\documentclass{article} 

\usepackage{amsmath,amsfonts,bm}

\def\1{\bm{1}}

\def\vq{{\bm{q}}}

\def\vv{{\bm{v}}}

\def\vx{{\bm{x}}}
\def\vy{{\bm{y}}}
\def\vz{{\bm{z}}}

\DeclareMathAlphabet{\mathsfit}{\encodingdefault}{\sfdefault}{m}{sl}
\SetMathAlphabet{\mathsfit}{bold}{\encodingdefault}{\sfdefault}{bx}{n}

\def\gD{{\mathcal{D}}}

\def\gH{{\mathcal{H}}}

\def\gM{{\mathcal{M}}}

\def\gX{{\mathcal{X}}}

\def\der{{\mathrm{d}}}

\usepackage[left=1in,top=1in,right=1in,bottom=1in,letterpaper]{geometry}
\usepackage{hyperref}
\usepackage{url}
\usepackage{microtype}
\usepackage{graphicx}
\usepackage{subfigure}
\usepackage{booktabs} 
\usepackage{amsmath}
\usepackage{amsfonts}
\usepackage{amsthm}
\usepackage{xspace}
\usepackage{soul}
\usepackage{mathtools}
\usepackage{floatrow}
\usepackage{breqn}
\newfloatcommand{capbtabbox}{table}[][\FBwidth]
\usepackage[T1]{fontenc}
\usepackage[font=small,labelfont=bf,tableposition=top]{caption}
\usepackage[bottom]{footmisc}
\usepackage{tablefootnote}
\usepackage[dvipsnames]{xcolor}

\DeclareCaptionLabelFormat{andtable}{#1~#2  \&  \tablename~\thetable}

\newtheorem{thm}{Theorem}[section]

\newtheorem{prop}[thm]{Proposition}
\newtheorem{lem}[thm]{Lemma}

\newtheorem{defn}[thm]{Definition}

\newlength\mylength
\newlength\mylengthLong
\title{Bridging Mean-Field Games and Normalizing Flows with Trajectory Regularization\thanks{This work is supported in part by an NSF Career Award (DMS-1752934)  and NSF DMS-2134168. }}

\author{Han Huang \thanks{H. Huang (huangh14@rpi.edu) is with the department of mathematics, Rensselaer Polytechnic Institute. }\qquad 
Jiajia Yu \thanks{J. Yu (yuj12@rpi.edu) is with the department of mathematics, Rensselaer Polytechnic Institute. } \qquad 
Jie Chen \thanks{J. Chen (chenjie@us.ibm.com) is with the MIT-IBM Watson AI Lab, IBM Research.}\qquad 
Rongjie Lai \thanks{Corresponding author. R Lai (lair@rpi.edu) is with the department of mathematics, Rensselaer Polytechnic Institute. }
}

\begin{document}

\date{}
\maketitle

\begin{abstract}
Mean-field games (MFGs) are a modeling framework for systems with a large number of interacting agents. They have applications in economics, finance, and game theory. Normalizing flows (NFs) are a family of deep generative models that compute data likelihoods by using an invertible mapping, which is typically parameterized by using neural networks. They are useful for density modeling and data generation. While active research has been conducted on both models, few noted the relationship between the two. In this work, we unravel the connections between MFGs and NFs by contextualizing the training of an NF as solving the MFG. This is achieved by reformulating the MFG problem in terms of agent trajectories and parameterizing a discretization of the resulting MFG with flow architectures. With this connection, we explore two research directions. First, we employ expressive NF architectures to accurately solve high-dimensional MFGs, sidestepping the curse of dimensionality in traditional numerical methods. Compared with other deep learning approaches, our trajectory-based formulation encodes the continuity equation in the neural network, resulting in a better approximation of the population dynamics. Second, we regularize the training of NFs with transport costs and show the effectiveness on controlling the model's Lipschitz bound, resulting in better generalization performance. We demonstrate numerical results through comprehensive experiments on a variety of synthetic and real-life datasets.
\end{abstract}

\section{Introduction}

Mean-field games (MFGs) are a powerful modeling framework for systems with a large number of interacting agents. They arise in the study of game theory~\cite{MFG_game_thry, Learning_MFG}, economics~\cite{MFG_econ1, MFG_econ2}, finance~\cite{MFG_fin1, MFG_fin2, MFG_fin3}, and industrial planning~\cite{MFG_industry1, MFG_industry2, MFG_industry3}. In general, the MFG formulation prescribes an identical objective for all agents and it seeks the optimal strategy over a time interval, which is analogous to solving for the Nash equilibrium in a finite $N$-player game. When the strategies concern spatial movement with a particular transport cost, the MFG is reduced to the optimal transport (OT) problem, which finds rich applications in signal processing and machine learning~\cite{Wass_image, Wass_signal, Wass_text, WGAN}.

A key theoretical underpinning for MFGs is the characterization of its optimality condition in terms of the Hamilton-Jacobi-Bellman (HJB) equation. Obtaining the Nash equilibrium of an MFG is thus reduced to solving a system of nonlinear partial differential equations (PDEs), for which a plethora of numerical methods have been developed based on the Eulerian framework~\cite{achdou2010mean, benamou2014augmented,benamou2017variational,benamou2000computational, jacobs2019solving,papadakis2014optimal,yu2021fast}. Nevertheless, the reliance on spatial discretization renders traditional approaches exponentially more expensive as the dimensionality grows. Hence, numerically solving MFGs in high dimensions remains a difficult problem, due to the curse of dimensionality. Meanwhile, machine learning-based approaches emerged lately. A recent method~\cite{ruthotto2020machine} successfully solves MFGs in quite high dimensions using deep neural networks. This method uses a Lagrangian-based approach to approximate the value function. However, besides penalizing the HJB equation, this method has to additionally solve the continuity equation that governs the evolution of the densities under the learned dynamics. 

Rather than approximating the value function in MFGs by using neural networks as discussed in~\cite{ruthotto2020machine}, we propose to approximate the trajectories instead. To this end, we establish a mathematical connection between the MFG framework and a popularly used family of generative models in deep learning---normalizing flows (NFs). Composed of a sequence of invertible mappings, NFs allow for the exact computation of the data likelihood. They can be used to model complex data distributions for density estimation and latent sampling~\cite{NF_survey}. Starting from the MFG problem, we reformulate it in terms of agent trajectories and canonically parameterize the discretization as an NF. As a result, the standard negative log-likelihood training of an NF is equivalent to a special case of the MFG without transport and interaction costs.

This connection between the MFGs and the NFs motivates us to explore two promising directions to improve the two. First, we approximate the trajectories of MFGs by employing expressive NF architectures, to accurately and efficiently solve high-dimensional MFG problems. The key novelty of our approach lies in its encoding of the continuity equation into the neural network, which mitigates the numerical errors incurred in the approximation and allows one to solve MFGs more accurately. We demonstrate the advantage of our approach in problems including OT, crowd motion, and multi-group path planning. Additionally, we conduct numerical analysis of a discretization scheme and explore the relationship between continuous and discretized MFGs. By characterizing the solution behaviors of discretized MFGs with the OT theory~\cite{OT_book}, we obtain an important insight, which sugggests that direct linear interpolation of the discretized problems in the temporal direction solves MFGs without the interaction term. Moreover, we adapt the universality theory from~\cite{universality_CL} to argue the effectiveness of using a particular kind of NFs---affine coupling flows---in parameterizing MFG trajectories. The resultant push-forward measure in the NF converges in distribution, corroborating accurate modeling of the population dynamics. 

Second, we improve NFs by introducing the transport cost to their training, resulting in \emph{trajectory-regularized flows}. Note that in high dimensions, there exists more than one flow that transforms one probability density to another. Therefore, the NF training problem is ill-posed and proper regularization is needed to incentivize the solution towards a more meaningful one. The OT theory suggests that this can be done via using the kinetic energy as the transport cost, which induces an OT plan that traverses the geodesic in the space of measures with the Wasserstein-2 metric~\cite{OT_book}. We show by using a variety of synthetic datasets that in the absence of the transport cost, one cannot expect existing NF models to approximate the OT trajectory.  Additionally, we find that the intermediate flows incur significantly less distortion when trained with transport costs. Regulating the kinetic energy, therefore, serves as an effective way of controlling the Lipschitz bound, which is intimately related to the model's generalization performance as well as its ability to counter adversarial samples~\cite{Lip_gen, Lip_adversarial}. Traditionally, the control over the Lipschitz bound for a neural network is implemented through explicit regularization such as the $l_2$ weight decay. To the best of our knowledge, this is the first work regularizing NFs by using transport costs and relating the regularization to the Lipschitz bound of the trained neural network, which proves to be more robust and effective than weight decay. We demonstrate the effectiveness of transport-based regularization on a variety of synthetic and real-life datasets and show improvements over popularly used NF models, such as RealNVP and the neural spline flow~\cite{RealNVP, NSF}.

\paragraph{Related Works}


Variational MFGs are a generalization of the dynamic formulation of OT~\cite{OT_book}, where the final density matching is encouraged but not enforced and transport costs other than the kinetic energy can be considered. Traditional methods for solving OT and variational MFG problems are well-studied in low dimensions~\cite{achdou2010mean, benamou2014augmented,benamou2017variational,benamou2000computational, jacobs2019solving,papadakis2014optimal,yu2021fast,cuturi2013sinkhorn}; and machine learning-based approaches emerged recently for solving high-dimensional MFGs. The work~\cite{ruthotto2020machine} is the first on solving high-dimensional deterministic MFGs by approximating the value function using deep neural networks. The work \cite{Tong_APACNet} solves high-dimensional MFGs in the stochastic setting using the primal-dual formulation of MFGs, where training is conducted similarly to that of generative adversarial networks~\cite{GAN,WGAN}.   





Our proposed trajectory-regularized flow can be built on any discrete NF. In general, existing NF models impose a certain structure on the flow mapping to enable fast evaluation of the log-determinant of the Jacobian. For a non-exhaustive list, RealNVP and NICE combine affine coupling flows to build simple yet flexible mappings~\cite{RealNVP, NICE}, while NSF and Flow$++$ use splines and a mixture of logistic CDFs as more sophisticated coupling functions~\cite{NSF, Flow++}. MAF and IAF are autoregressive flows with affine coupling functions~\cite{MAF, IAF}, while NAF and UMNN parameterize the coupling function with another neural network~\cite{NAF, UMNN}.

Not the focus of this work, continuous NF models consider the ODE limit of the discrete version, with significant changes to the computation of the parameter gradient and the change of variables. In these models, neural ODE outlines the theoretical framework~\cite{NODE} and FFJORD uses the Hutchinson estimator to approximate the trace term incurred in the density evaluation~\cite{FFJORD}. OT-Flow introduces the OT cost to continuous NFs, suggesting a parameterization of the HJB potential that allows for exact trace computations~\cite{OT_flow}. In addition, recent works found that the transport cost can enhance and stabilize the training of continuous NFs~\cite{RNODE}. For example, the OT theory suggests that the optimal agent trajectories are straight; this characterization allows the ODE solver in OT-Flow to integrate with larger time steps and reduce the computational cost~\cite{OT_flow}.

\paragraph{Roadmap} 

The rest of the paper is organized as follows. We provide the  mathematical background for MFGs and NFs in Section~\ref{sec:background} and elaborate the connections between the two in Section~\ref{sect_briding}, via a reformulation and discretization of the MFG objective. In Section~\ref{sec:theory}, we perform theoretical analysis to relate the discretized and continuous MFG problems, as well as establishing universality theorems for using NFs to solve MFGs. In Section~\ref{sec:expertments}, we provide numerical results that demonstrate the effectiveness and accuracy of solving high-dimensional MFGs by using NFs. Meanwhile, we illustrate the effect of regularized NF training, stress the ability of acquiring a lower Lipschitz bound, and show its superiority to weight decay in real-life datasets.


\section{Mathematical Background}
\label{sec:background}
\subsection{Mean-Field Games}

The study of MFGs considers classic $N$-player games at the limit of $N\to \infty$~\cite{MFG_varMFG,huang2006large,huang2007large}. Assuming homogeneity of the objectives and anonymity of the players, we set up an MFG with a continuum of non-cooperative rational agents distributed in a spatial domain $\Omega \subseteq \mathbb{R}^d$ and a time interval $[0, T]$. For any fixed $t\in [0,T]$, we denote the population density of agents by $p(\cdot, t) \in \mathcal{P}(\Omega)$, where $\mathcal{P}(\Omega)$ is the space of all probability densities over $\Omega$. For an agent starting at $\vx_0\in \Omega$, their position over time follows a trajectory $\vx: [0,T] \to \Omega$ governed by
\begin{equation} \label{agent_traj}
\left\{\begin{aligned}
    \der \vx(t) &= \vv(\vx(t),t)\der t, \quad \forall t\in [0,T]\\
    \vx(0) &= \vx_0,
\end{aligned}\right.
\end{equation}
where $\vv: \Omega \times [0,T] \to \Omega$ specifies an agent's action at a given time. For simplicity, there are no stochastic terms in~\eqref{agent_traj} in this work. Then, each agent's trajectory is completely determined by $\vv(\vx,t)$. To play the game over an interval $[t, T]$, each agent seeks to minimize their objective:
\begin{equation} \label{agent_cost}
\begin{aligned}
    J_{\vx_0,t }(\vv,p) &\coloneqq \int_t^T [L(\vx(s),\vv(\vx(s),s)) + I(\vx(s), p(\vx(s), s)) ] \der s + M(\vx(T), p(\vx(T), T))\\
    &\text{s.t. } \eqref{agent_traj} \text{ holds}.
\end{aligned}
\end{equation}
Each term in the objective denotes a particular type of cost. The running cost $L: \Omega \times \Omega \to \mathbb{R}$ is incurred due to each agent's action. A commonly used example for $L$ is the kinetic energy $L(\vx,\vv) = \|\vv\|^2$, which accounts for the total amount of movement along the trajectory. The running cost $I: \Omega \times \mathcal{P}(\Omega) \to \mathbb{R}$ is accumulated through each agent interacting with one another or with the environment. For example, this term can be an entropy that discourages the agents from grouping together, or a penalty for the agents to collide with an obstacle in the space~\cite{ruthotto2020machine}. The terminal cost $M: \Omega \times \mathcal{P}(\Omega) \to \mathbb{R}$ is computed from the final state of the agents. Typically, $M$ measures a discrepancy between the final distribution $p(\cdot, T)$ and a desirable distribution $p_1 \in \mathcal{P}(\Omega)$.

To solve the MFG, we define the value function $u: \Omega \times [0,T] \to \mathbb{R}$ as
\begin{equation} \label{value_fun}
\begin{aligned}
u(\vx_0,t) &\coloneqq \inf_{\vv} J_{\vx_0,t}(\vv,p),
\qquad  \text{s.t. } \eqref{agent_traj} \text{ holds}.
\end{aligned}
\end{equation}
One can show that $u(\vx,t)$ and $p(\vx,t)$ are solutions to the following HJB equation and continuity equations~\cite{MFG_varMFG,huang2006large,huang2007large}
\begin{equation} \label{HJB}
\left\{\begin{aligned}
    -\partial_t u(\vx,t) + H(\vx,\nabla_{\vx} u(\vx,t)) &= I(\vx,p(\vx,t)) \\
    \partial_t p(\vx,t) + \nabla_{\vx} \cdot (p(\vx,t)\vv(\vx,t)) &= 0 \\
    p(\vx,0) = p_0(\vx),  \quad \text{and} \quad   u(\vx,t)& = M(\vx, p(\vx,t)),
\end{aligned}\right.
\end{equation}
where $H: \Omega \times \Omega \to \mathbb{R}$ is the Hamiltonian
$  H(\vx, \vq) \coloneqq \sup_\vv -\langle \vq , \vv\rangle - L(\vx, \vv)$; 
$p_0 \in \mathcal{P}(\Omega)$ is the initial population density at $t=0$; and $\vv(\vx,t)$ denotes the optimal strategy for an agent at position $\vx$ and time $t$. It can be shown that the optimal strategy satisfies $\vv(\vx,t) = -\nabla_p H(\vx, \nabla_{\vx} u(\vx,t))$. We note that the continuity equations in the last two lines of~\eqref{HJB} are a special case of the Fokker-Planck equation with no diffusion terms. 



The setup~\eqref{agent_traj}--\eqref{agent_cost} concerns the strategy of an individual agent. Under suitable assumptions, the work~\cite{MFG_varMFG} developed a macroscopic formulation of the MFG that models the collective strategy of all agents. Suppose there exist functionals $\mathcal{I}, \mathcal{M}: \mathcal{P}(\Omega) \to \mathbb{R}$ such that
\begin{align*}
    I(\vx,p) = \frac{\delta \mathcal{I}}{\delta p}(\vx), \quad M(\vx,p) = \frac{\delta \mathcal{M}}{\delta p}(\vx),
\end{align*}
where $\frac{\delta}{\delta p}$ is the variational derivative. Then, the functions $p(\vx,t)$ and $\vv(\vx,t)$ satisfying~\eqref{HJB} coincide with the optimizers of the following variational problem: 
\begin{equation}\label{var_MFG}
\begin{split}
    \inf_{p,\vv} J(p,\vv) \coloneqq \int_0^T \int_\Omega L(\vx, \vv(\vx,t),p(\vx,t) \der\vx \der t &+ \int_0^T \mathcal{I} (p(\cdot, t))\der t + \mathcal{M}(p(\cdot, T))\\
    s.t. \quad \partial_t p(\vx,t) + \nabla_{\vx} \cdot (p(\vx,t)\vv(\vx,t)) &= 0, \quad \vx\in \Omega, t\in [0,T]\\
    p(\vx,0) &= p_0(\vx), \quad \vx\in \Omega.
\end{split}
\end{equation}
This formulation is termed the \emph{variational MFG} and it serves as the main optimization problem this work solves. 
It is known that traditional numerical PDE methods that solve~\eqref{HJB} are computationally intractable in dimensions higher than three, as the number of grid points grows exponentially. This challenge motivates us to solve the MFG problem with parameterizations of the agent trajectory based on NFs.

\subsection{Normalizing Flows}

NFs are a family of invertible neural networks that find applications in density estimation, data generation, and variational inference~\cite{NICE, RealNVP, NSF, Flow++, MAF, IAF, NAF}. In machine learning, NFs are considered deep generative models, which include variational autoencoders~\cite{VAE} and generative adversarial networks~\cite{GAN} as other well-known examples. The key advantage of NFs over these alternatives is the exact computation of the data likelihood, which is made possible by invertibility.

Mathematically, suppose we have a dataset $\gX= \{\vx_n\}_{n=1}^N \subseteq \mathbb{R}^d$ generated by an underlying data distribution $P_1$; that is, $\vx_n \overset{\mathrm{iid}}{\sim} P_1$. We define a \textit{flow} to be an invertible function $f_\theta:\mathbb{R}^d \to \mathbb{R}^d$ parameterized by $\theta \in \mathbb{R}^W$, and a \textit{normalizing flow} to be the composition of a sequence of flows: $F_\theta = f_{\theta_K} \circ f_{\theta_{K-1}} \circ ... \circ f_{\theta_2} \circ f_{\theta_1}$ parameterized by $\theta = (\theta_1, \theta_2, ..., \theta_K)$. To model the complex data distribution $P_1$, the idea is to use an NF to gradually transform a simple base distribution $P_0$ to $P_1$. Formally, the transformation of the base distribution is conducted through the push-forward operation $F_{\theta*}P_0$, where  $F_{\theta*}P(A) \coloneqq P(F_{\theta}^{-1}(A))$ for all measurable sets $A\subset\Omega$. The aim is that the transformed distribution resembles the data distribution. Thus, a commonly used loss function for training the NF is the KL divergence:
\begin{equation}\label{KL_NF}
    \min_\theta \quad  \mathcal{D}_{KL} (P_1 || F_{\theta*}P_0).
\end{equation}
If the measure $P$ admits density $p$ with respect to the Lebesgue measure, the density $F_{\theta*}p$ associated with the push-forward measure $F_{\theta*}P$ satisfies the change-of-variable formula~\cite{NF_survey}:
\begin{equation}\label{change_of_var}
    F_{\theta*}p(\vx) = p(F_\theta^{-1}(\vx))|\det \nabla F_\theta^{-1}(\vx)|,
\end{equation}
where $\nabla F_\theta^{-1}$ is the Jacobian of $F_\theta^{-1}$. This allows one to tractably compute the KL term in~\eqref{KL_NF}, which is equivalent to the negative log-likelihood of the data up to a constant (independent of $\theta$):
\begin{equation} 
\begin{split}
    \mathcal{D}_{KL} (P_1 || F_{\theta*}P_0) &= - \mathbb{E}_{\vx\sim P_1} [\log F_{\theta*}p_0(\vx)] + \text{const.} \\
    &= -\mathbb{E}_{\vx\sim P_1} \left[\log p_0(g_{\theta_1} \circ ... \circ ... g_{\theta_K}(\vx)) + \sum_{k=1}^K \log |\det \nabla g_{\theta_k}(\vy_{K-k})|\right] + \text{const.},
\end{split}
\label{NF_NLL}
\end{equation}
where $g_{\theta_i} \coloneqq f_{\theta_i}^{-1}$ are the flow inverses, and $\vy_j \coloneqq g_{\theta_{K-j}} \circ ... \circ g_{\theta_K}(\vx), \forall j = 1, 2, ..., K-1, \vy_0 \coloneqq \vx$.

We remark that there exists flexibility in choosing the base distribution $P_0$, as long as it admits a known density to evaluate on. In practice, a typical choice is the standard normal $\mathrm{N}(0,I)$. In addition, \eqref{NF_NLL} requires the log-determinant of the Jacobians, which take $O(d^3)$ time to compute in general. Existing NF architectures sidestep this issue by meticulously using functions that follow a (block) triangular structure~\cite{NF_survey}, so that the log-determinants can be computed in $O(d)$ time from the (block) diagonal elements.

One popular family of NF models is the coupling flow~\cite{NF_survey}, which composes individual flows in the following form
\begin{equation} 
\begin{aligned}
    f(\vx_1, \vx_2) = (h(\vx_1; \theta(\vx_2)), \vx_2) \coloneqq (\vy_1, \vy_2), 
\end{aligned}
\end{equation}
where $(\vx_1, \vx_2) \in \mathbb{R}^a \times \mathbb{R}^{d-a}$, $\theta: \mathbb{R}^{d-a} \to \mathbb{R}^d$ is the called the \emph{conditioner}, and $h: \mathbb{R}^d \to \mathbb{R}^d$ is the coupling function. Suppose $h$ is invertible, then the inverse $f^{-1}$ is:
\begin{equation} 
\begin{aligned}
    (\vx_1, \vx_2) = (h^{-1}(\vy_1; \theta(\vx_2)), \vy_2).
\end{aligned}
\end{equation}
In common scenarios, $h$ is defined element-wise, such that invertibility is ensured by strict monotonicity. For example, a simple yet expressive choice is the affine function $h(\vx;\theta) = \theta_1 \vx + \theta_2, \theta_1 \ne 0, \theta_2 \in \mathbb{R}$. This gives rise to the RealNVP flow~\cite{RealNVP}. Other popular coupling functions include splines, mixture of CDFs, and functions parameterized by neural networks, which are shown to be more expressive than the affine coupling function~\cite{NSF, Quad_spline, cubic_spline, Flow++, SOS, NAF}.


\section{Bridging MFGs and NFs}
\label{sect_briding}
In this section, we reformulate the MFG problem as a generalization of the NF training loss. This reformulation relates the two models and opens opportunities to improve both. The contribution is two-fold. On the one hand, we solve high-dimensional MFG problems by using an NF parameterization. The NF model encodes the discretized trajectory of the agents in its network architecture, which allows for efficient optimization through evaluation of the log-determinants in $O(d)$ cost. We show that the expressitivity of the NF architecture allows one to tractably solve high-dimensional MFG problems with small error. On the other hand, under the MFG framework we introduce a model-agnostic regularization to improve the training of NFs. In the MFG language, existing NF models are trained by using only the terminal cost (i.e., the KL loss). With the introduction of the transport cost, the intermediate flows become better behaved and more robust, owing to a smaller Lipschitz bound. It turns out that the learned NF improves the matching of the densities, as evident in various synthetic and real-life examples. 

\subsection{Trajectory-Based Formulation of MFGs}

Rather than solving for the density $p(\vx,t)$ and the action $\vv(\vx,t)$, we reparameterize the problem~\eqref{var_MFG} to directly work with agent trajectories, from which $p(\vx,t)$ and $\vv(\vx,t)$ can be derived. Let $P(\cdot,t)$ be the measure that admits $p(\cdot,t)$ as its density for all $t\in [0,T]$. Define the agent trajectory as $F: \mathbb{R}^{d} \times \mathbb{R} \to \mathbb{R}^d$, where $F(\vx,t)$ is the position of the agent starting at $\vx$ and having traveled for time $t$. It satisfies the following ordinary differential equation:
\begin{equation} \label{NF_traj}
\left\{\begin{aligned}
    \partial_t F(\vx,t) &=  \vv(F(\vx,t), t), \quad \vx\in \Omega, t\in[0,T] \\
    F(\vx,0) &= \vx, \quad \vx\in \Omega.
\end{aligned}\right.
\end{equation}
The evolution of the population density is determined by the movement of agents. Thus, $P(\cdot,t)$ is simply the push-forward of $P_0$ under $F(\cdot, t)$; namely, $P(\vx, t) = F(\cdot, t)_* P_0(\vx)$, whose associated density satisfies
\begin{equation} \label{density_push-forward}
\begin{aligned}
    p(\vx,t) &= \mathrm{d}(F(\cdot, t)_* P_0)(\vx),
\end{aligned}
\end{equation}
where $\mathrm{d}(F(\cdot, t)_* p_0(\vx))$ is a Radon-Nikodym derivative. Note that this definition of $p(\vx,t)$ automatically satisfies the continuity equation as well as the initial condition for $p(\vx,t)$ in~\eqref{var_MFG}. Unless specified otherwise, we assume $T=1$ and use the $L_2$ transport cost
 $   L(\vx,\vv(\vx,t),p(\vx,t)) = \lambda_L p(\vx,t) \|\vv(\vx,t)\|_2^2$ from now on. Here, $\lambda_L \geq 0$ is treated as a model (hyper-)parameter. Furthermore, we can apply a change of variables on~\eqref{density_push-forward} to rewrite the integral involving $L(\vx, \vv,p)$:
\begin{align*}
    \int_0^1 \int_\Omega L(\vx, \vv(\vx,t), p(\vx,t)) \der\vx \der t &=\lambda_L \int_0^1 \int_\Omega  p(\vx,t) \|\vv(\vx,t)\|_2^2 \der\vx \der t\\
    &= \lambda_L \int_0^1 \int_\Omega \|\vv(\vx,t)\|_2^2 d(F(\cdot, t)_* P_0)(\vx) \der t\\
    &= \lambda_L\int_0^1 \int_\Omega  p_0(\vx) \|\partial_t F(\vx,t)\|_2^2 \der\vx \der t.
\end{align*}
Similarly, we can rewrite the interaction cost $\mathcal{I}(p(\cdot, t))$ and the terminal cost $\mathcal{M}(p(\cdot, T))$ as:
\begin{equation} \label{I_M_traj}
\begin{aligned}
    \mathcal{I}(p(\cdot, t)) &= \mathcal{I}(F(\cdot, t)_*p_0),\\
    \mathcal{M}(p(\cdot, T)) &= \mathcal{M}(F(\cdot, T)_*p_0).
\end{aligned}
\end{equation}
%
%
Therefore, \eqref{var_MFG} becomes:
\begin{equation} \label{var_MFG_traj}
\begin{aligned}
    \inf_{F} \quad \lambda_L \int_0^1 \int_\Omega p_0(\vx) \|\partial_t F(\vx,t)\|_2^2 \der\vx \der t &+ \int_0^1 \mathcal{I}(F(\cdot, t)_*p_0)\der t + \mathcal{M}(F(\cdot, 1)_*p_0) \coloneqq \mathcal{L} (F)\\
    s.t. \quad  &F(\vx,0) = \vx.
\end{aligned}
\end{equation}


It is worth noting that~\eqref{var_MFG_traj} removes the continuity-equation constraint. The reformulated problem optimizes over agent trajectories $F$, which allows one to parameterize with an NF. In addition, the remaining constraint $F(\cdot, 0) = \mathrm{Id}$ will be automatically integrated into the NF, such that it can be trained in an unconstrained manner.

\subsection{Discretization of MFG Trajectories with NFs}

Starting from the base density $p_0$, each flow in the NF advances the density one step forward. Hence, the NF model specifies a discretized evolution of the agent trajectories. Formally, we discretize the trajectory $F(\vx,t)$ in time with regular grid points $ t_i \coloneqq i\cdot \Delta t, \forall i = 0,1,...K$. Thus, the grid spacing is $\Delta t \coloneqq \frac{1}{K}.$ Denote $F(\vx,t_i) \coloneqq F_i(\vx)$.  By approximating $\partial_t F(\vx,t)$ with the forward difference, the transport cost becomes:
\begin{equation}\label{eqn:dis}
\begin{split}
    \lambda_L \int_0^1 \int_\Omega  p_0(\vx) \|\partial_t F(\vx,t)\|_2^2 \der\vx \der t &= \lambda_L  \int_0^1 \mathbb{E}_{z\sim P_0} [\|\partial_t F(\vz,t)\|_2^2] \der t\\
    &\sim \lambda_L K \cdot \mathbb{E}_{z\sim P_0} \left[\sum_{i=0}^{K-1}\|F_{i+1}(\vz) - F_i(\vz)\|_2^2\right].
\end{split}
\end{equation}
Moreover, suppose that the terminal cost is computed as the KL divergence as in standard NFs:
\begin{equation} \label{var_MFG_NF_terminal}
\begin{aligned}
    \mathcal{M}(F(\cdot, T)_*p_0) &= \lambda_{\mathcal{M}} D_{KL} (P_1 || F(\cdot, T)_*P_0) = \lambda_{\mathcal{M}} D_{KL} (P_1 || F_{K_*}P_0).
\end{aligned}
\end{equation}
For simplicity of exposition, for now we omit the interaction term (that is, $\mathcal{I} \equiv 0$). We can discretize the interaction cost $\mathcal{I}(F(\cdot, t)_*P_0)$, similarly to~\eqref{eqn:dis}, given its exact form. Subsequently, we will give numerical examples for the case $\mathcal{I} \not \equiv 0$ (see Sections~\ref{subsec:CrowdMotion} and~\ref{subsec:pathplanning}). 

For any $F$, the discretized objective value converges to the continuous version at $O(\frac{1}{K})$. However, we show later that the optimal discretized and continous objective values in fact agree for any fixed number of grid points used. Hence, there introduce no additional errors by discretizing the original MFG problem, in the absence of interaction costs.

To solve the discretized MFG, we define the flow maps for $0\leq t_i \leq t_j \leq 1$ as  $\Phi_{t_{i}}^{t_{j}} (\vx_0) = \vx(t_{j})$, where
\begin{equation}\label{flow_map}
\left\{\begin{aligned}
    \der\vx(t) &= \vv(\vx(t),t)\der t, \quad t\in [t_i, t_j]\\
    \vx(t_i) &= \vx_0.
\end{aligned}\right.
\end{equation}
Here, we leverage the semi-group property $\Phi_{t_{b}}^{t_{c}} \circ \Phi_{t_{a}}^{t_{b}} = \Phi_{t_{a}}^{t_{c}}$ to decompose the agent trajectories into a sequence of flow maps.  It follows that $F_k = \Phi_{t_{k-1}}^{t_k} \circ \Phi_{t_{k-2}}^{t_{k-1}} \circ ... \circ \Phi_{t_0}^{t_1}, \forall k = 1, 2, ... K$, and we parameterize each $ \Phi_{t_i-1}^{t_i}: \mathbb{R}^d \to \mathbb{R}^d$ by a flow model $f_{\theta_i}$. Denote $F_{\theta^k} \coloneqq f_{\theta_k} \circ ... \circ f_{\theta_1}$, $F_{\theta^0} \coloneqq Id$, and $F_\theta \coloneqq F_{\theta^K}$. The discretized  MFG problem with flow-parameterized trajectories is:
\begin{equation} \label{var_MFG_NF_final}
\begin{aligned}
    \inf_{\theta} \quad \lambda_L K\cdot \mathbb{E}_{\vz\sim P_0} &\left[ \sum_{k=0}^{K-1}\|F_{\theta^{k+1}}(\vz) - F_{\theta^k}(\vz)\|_2^2\right] + \lambda_{\mathcal{M}} \mathcal{D}_{KL} (P_1||F_{\theta_*}P_0),
\end{aligned}
\end{equation}
where $\lambda_{\mathcal{M}}, \lambda_L \geq 0$ are weights.


By comparing~\eqref{var_MFG_NF_final} to the NF optimization problem~\eqref{NF_NLL}, we see that training a standard NF is equivalent to solving a special MFG with no transport and interaction costs. In the past, many highly flexible flows were developed that are capable of solving this special MFG~\cite{NAF, NSF, RealNVP, glow, Flow++}. As a result, the trajectory-based formulation opens the door to employ expressive NFs to accurately solve MFGs. In addition, note that the density evolution can be directly obtained from the agent trajectories, which are the outputs of the individual flows in the NF model. Thus, this formulation sidesteps the additional effort required to solve the Fokker-Planck equation~\eqref{HJB}, resulting in a more efficient solution of the population dynamics. 

We also note that the computation of~\eqref{var_MFG_NF_final} involves both a forward (normalizing) and an inverse (generating) pass of the NF: the forward pass computes $\mathcal{D}_{KL} (P_1||F_{\theta_*}P_0)$ via~\eqref{NF_NLL} and the inverse pass yields the transport cost. For NFs that have comparable run time for passes in both directions (e.g. coupling flows), adding the transport cost only doubles the overall computation time. For NFs where the inverse pass is much slower than the forward pass (or vice versa), such as autoregressive flows~\cite{NF_survey}, the additional computational effort associated with the transport cost may be prohibitive. To remedy this, we propose an alternative, equivalent formulation that incurs less overhead. 

\paragraph{Alternative Formulation}

Recall that the forward agent trajectory is defined in~\eqref{NF_traj}, such that $p(\vx, t) = F(\cdot, t)_* p_0(\vx)$. Equivalently, we can consider the backward trajectory $G: \mathbb{R}^d \times \mathbb{R} \to \mathbb{R}^d$:
\begin{equation} \label{NF_traj_rev}
\begin{aligned}
    \partial_t G(\vx,1-t) &= \vv(G(\vx,1-t), t), \quad \vx\in \Omega, t\in[0,T] \\
    G(\vx,0) &= \vx, \quad \vx\in \Omega.
\end{aligned}
\end{equation}
Therefore, $G(\cdot, t)_* P_1 = P(\cdot, 1-t), \forall t\in [0,1]$. As a result, the objective~\eqref{var_MFG_NF_final} becomes:
\begin{equation} \label{NF_traj_rev_2}
\begin{aligned}
    \inf_{\theta} \quad \lambda_L K\cdot \mathbb{E}_{\vx\sim P_1} &\left[ \sum_{k=0}^{K-1}\|G_{\theta^k}(\vx) - G_{\theta^{k+1}}(\vx)\|_2^2\right] + \lambda_{\mathcal{M}} \mathcal{D}_{KL} (P_1||F_{\theta*}P_0), \\
\end{aligned}
\end{equation}
where $G_{\theta^k}(\vx) = g_{\theta_k} \circ g_{\theta_{k-1}} \circ ... \circ g_{\theta_1}(\vx), \forall k = 1, 2,..., K$ are the flow inverses, and $G_{\theta^0}(\vx) \coloneqq \vx$. Not that the evaluation of $\mathcal{D}_{KL} (P_1||F_{\theta*}P_0)$ requires the inverse of $F_\theta$, which is expressed in the $g_{\theta_i}$'s defined in~\eqref{NF_NLL}.  We defer the detailed derivation to Appendix~\ref{append:alternative}.

The transport cost in the equivalent formulation~\eqref{NF_traj_rev_2} can be computed via a single forward pass on the training data, which can be done together with the negative log-likelihood computation. As a result, the addition of the transport cost incurs minimal overhead compared to the standard NF training.

\subsection{Regularizing NF Training with MFG Costs}
\label{subsec:RegularizingNFs}
The aforementioned method for solving high-dimensional MFGs using NFs leads to the view of NFs as a special type of MFG problems.  Comparing the MFG formulation~\eqref{var_MFG_NF_final} with the standard NF training~\eqref{NF_NLL}, one can interpret the MFG transport cost as a regularization term in the NF training with strength $\lambda_L / \lambda_{\mathcal{M}}$; we call this regularization a \emph{trajectory regularization}. It is straightforward to implement it on top of any existing discrete NF models (not necessarily restricted to coupling flows). Since the transport cost aggregates squared differences between the input and the output of each flow, the regularization encourages smoothness in the intermediate transforms.

From the perspective of inverse problems, there exist many flows (transport plans) that can transform the base $P_0$ to $P_1$. Therefore, regularization helps ameliorate the inevitably ill-posed nature of NF training. Through OT theory, the $L_2$ regularized evolution of densities is characterized as the geodesic connecting $P_0$ and $P_1$ in the space of measures equipped with the Wasserstein-2 metric~\cite{W2_metric}. Since $P_0$ is always chosen to be absolutely continuous for NFs, the OT between the two densities will be a unique and bijective mapping known as the Monge map, $T(\vx)$. The same result holds for the discretized MFG problem, shown later in Theorem~\ref{cont_disc_MFG_soln}. Therefore, the NF training problem is no longer ill-posed. Furthermore, the optimal trajectory is then a linear interpolation between each point $\vx$ and its destination $T(\vx)$: $F^*(\vx,t) = (1-t)\vx + tT(\vx)$. 

From the perspective of machine learning, the transport cost provides an effective framework to control the flow's Lipschitz bound, which closely correlates with the robustness and the generalization capability~\cite{Lip_norm_bound, Lip_adversarial}. Note that by enforcing smoothness in the time domain, the trajectory regularization also encourages smoothness in the spatial domain, owing to the triangle inequality:
\begin{align}
\label{eqn:LipRegu}
    \|F_K(\vx) - F_K(\vy)\|^2_2 &= \left\| \sum_{i = 0}^{K-1} (F_{i+1}(\vx) - F_i(\vx)) - \sum_{i = 0}^{K-1} (F_{i+1}(\vy) - F_i(\vy)) + (\vx - \vy) \right\|^2 \nonumber\\
    &\leq 2\|\vx-\vy\|^2_2 +  K \sum_{i=0}^{K-1} \|F_{i+1}(\vx) - F_i(\vx)\|_2^2 + K \sum_{i=0}^{K-1} \|F_{i+1}(\vy) - F_i(\vy)\|_2^2,
\end{align}
where $F_0 = Id, F_i = f_i \circ ... f_1, i=1,...,K$, and $F_K$ is the NF model. By regularizing with the transport cost as in~\eqref{var_MFG_NF_final}, the right-hand side can be controlled. 

As observed by the authors of~\cite{Lip_gen}, explicit regularization methods such as $l_2$ weight decay has little impact on the model's generalization gap. In contrast, the trajectory regularization is implicit in that it imposes no direct penalty on the parameters. We find it to be much more effective than weight decay in reducing the NF's Lipschitz bound. We will demonstrate in the numerical experiments that a regularized evolution can reduce overfitting and improve density estimation.

\section{Theoretical Analysis}
\label{sec:theory}
We first characterize the solutions of the MFG problem and its discretized analogue, when the interaction cost is absent and the transport cost is the kinetic energy. We prove that the discretized objective converges to the continuous version linearly, and that the minimum values agree between the two problems, for any number of discretized points. Through reduction to an OT problem, we show uniqueness of the optimizer for discretized MFG problems, thereby alleviating ill-posedness of the NF training objective.

In addition, we show that with mild assumptions on the MFG solution, certain classes of normalizing flows can be used to approximate the optimal trajectories of any MFG problem in $L_2$ locally. Our results build on the work~\cite{universality_CL}, which establishes universality results for a broad class of coupling flows. In addition, we state that under the learned mapping, the  evolution of the probability measure converges to the ground truth in distribution, at every discretized time step.
We sketch the intuitions here and defer the proofs to the appendix.

\subsection{Relating Continuous and Discretized MFGs}

Our analysis focuses on the MFG problem without the interaction cost. We recall the discretized MFG problem with $\mathcal{I} \equiv 0$ here for convenience:
\begin{equation} \label{MFG_K}
\begin{aligned}
    &\inf_{F = \{F_i\}_{i=0}^K} \quad \lambda_L K\cdot \mathbb{E}_{z\sim P_0} \left[ \sum_{k=0}^{K-1}\|F_{i+1}(\vz) - F_{i}(\vz)\|_2^2\right] +  \mathcal{M}(F_{K*}p_0) \coloneqq \mathcal{L}_K(F) \\
    &s.t. \quad  F_0(\vz) = \vz.
\end{aligned}
\end{equation}

First, we note that in the absence of the interaction cost, the discretized objective converges to the continuous analogue in $O(\frac{1}{K})$ when evaluated on any $F$.

\begin{thm}
\label{thm:disc_error}
Assume $\mathcal{I} \equiv 0$. The discretized MFG objective value converges to the continuous version at $O(\frac{1}{K})$ for any $F$, when both objectives are finite:
\begin{align*}
    |\mathcal{L}_K(F) - \mathcal{L}(F)| = O\left(\frac{1}{K}\right).
\end{align*}

\end{thm}

The proof is based on a straightfoward Taylor expansion, similar to that used in finite difference error analysis. In general, one can use a higher order scheme to discretize~\eqref{var_MFG_traj} in time, yielding a higher order convergence rate. For example, the composite Simpson's rule is used in our numerical experiments. 

Next, we consider the optimization and establish two lemmas that characterize the optimal trajectories of both the discretized and the continuous MFG problems.

\begin{lem}\label{cont_MFG_soln}
Assume $P_0$ is absolutely continuous with density $p_0$. Let $F^*(\vz,t)$ be an optimizer for~\eqref{var_MFG_traj} with $\mathcal{I} \equiv 0$. Define $F^*(\vz,1) \coloneqq  F_1(\vz)$. Then, $F^*(\vz,t) = (1-t)\vz + tF_1(\vz)$. 

\end{lem}

Lemma~\ref{cont_MFG_soln} can be shown by recasting the MFG as a suitable OT problem. Then, the straight trajectories follow from the OT theory. A similar result holds for the discretized MFG problem.

\begin{lem}\label{disc_MFG_soln}

Assume $P_0$ is absolutely continuous with density $p_0$. Let $\{F^*_i\}_{i=0}^K$ be an optimizer for~\eqref{MFG_K} with $F_K^* \coloneqq F_1$. Then, $F_i^* = (1-\frac{i}{K})\vz + \frac{i}{K}F_1, i=0, 1, ..., K$.

\end{lem}

Lemma~\ref{disc_MFG_soln}  can be shown by recasting the problem in a similar fashion, followed by solving its KKT system. Combining Lemma~\ref{cont_MFG_soln} and Lemma~\ref{disc_MFG_soln}, we can show that the optimal value of the discretized problem agrees with its continuous counterpart.


\begin{thm}\label{cont_disc_MFG_soln}
Assume $\mathcal{I} \equiv 0$. Both the continuous~\eqref{var_MFG_traj} and the discretized~\eqref{MFG_K} MFG problems admit unique optimizers and their optimal values agree. In addition, the trajectory obtained from linearly interpolating the optimizers  of the discretized MFG problem, $F(\vz,t) = (1-t)\vz + tF_K^*(\vz)$, is the optimizer of the continuous MFG problem.
\end{thm}

Theorem~\ref{cont_disc_MFG_soln} has two important implications. First, when we use an NF to parameterize and solve the MFG problem with no interaction costs, the optimizer incurs no discretization error. Therefore, it is not necessary to use higher order discretization schemes for the problem. Second, the discretized MFG problem, which corresponds to the training objective~\eqref{var_MFG_NF_final} modulo parameterization, admits a unique solution. In other words, the regularized NF training is well-posed.

\subsection{Universality of NFs for Solving MFGs}

The universality results from~\cite{universality_CL} indicate that NFs built with affine coupling transforms and sufficiently expressive conditioners are universal approximators for $\gD^2$ locally in $L_p, p\in [1,\infty)$. Here, $\gD^2 = \{ f: Dom(f) \to Im(f) \subseteq \mathbb{R}^d\}$, where $f$ is $C^2$ diffeomorphic, and $Dom(f) \subseteq \mathbb{R}^d$ is open and diffeomorphic to $\mathbb{R}^d$. As a result, affine coupling flows are well-suited for approximating the optimal MFG trajectories. To be more precise, we consider the following sets introduced in~\cite{universality_CL}. 

\begin{defn}
Define the set of single coordinate affine coupling flows with conditioner class $\gH$ as $\mathbf{\gH\textbf{-}ACF} = \{\Psi_{d-1, s,t} : s,t\in \gH \}$, where $\Psi_{k,s,t}(\vx_{\leq k}, \vx_{>k}) = (\vx_{\leq k}, \vx_{>k} \odot \exp{(s(\vx_{\leq k})}) + t(\vx_{\leq k}))$ is an affine coupling flow. 
Further, define the set of invertible neural networks $\mathbf{INN_{\gH\textbf{-}ACF}} = \{W_1 \circ g_1 \circ ... \circ W_n \circ g_n : n\in \mathbb{N}, g_i\in \mathbf{\gH\textbf{-}ACF}, W_i(\vx) = A_ix + b_i, A_i \text{ invertible} \}$, which is a family of normalizing flows with flow transforms from the class of invertible functions $\mathbf{\gH\textbf{-}ACF}$.   
\end{defn}

Note that the affine coupling flows used in our numerical experiments belong to $\mathbf{INN_{\gH\textbf{-}ACF}}$, where $\gH$ is the class of ReLU networks. It turns out that the universality of $\mathbf{INN_{\gH\textbf{-}ACF}}$ provides an ideal characterization for the approximate solutions of~\eqref{var_MFG_traj}.

\begin{thm}\label{thm_optimizer_conv}Let $\gH$ be sup-universal for $C^\infty_c(\mathbb{R}^{d-1})$ and contain piecewise $C^1$ functions. 
For a fixed $K\in \mathbb{N}$, let the optimizer of~\eqref{var_MFG_traj} be $F^*(\vz,t)$ and its evaluation on the grid points $t_i = i\Delta t$ be $F^*(\vz,t_i) = \Phi_{t_{i-1}}^{t_i} \circ ... \Phi_{t_0}^{t_1} , i=1,...,K$. Assume $F^*(\vz,t_i) \in \gD^2, i=1,...,K$. Then, for all $\epsilon > 0$, $p\in [1,\infty)$, and compact sets $A$, there exists $\{F_{\theta^{i}}\}_{i=1}^K \subset \mathbf{INN_{\gH\textbf{-}ACF}}$ such that
\begin{align*}
    \|F^*(\cdot,t_i) - F_{\theta^{i}}\|_{p,A} < \epsilon, \quad  i=1,..., K.
\end{align*}

\end{thm}

This result can be directly obtained from the conclusion in~\cite{universality_CL}. 
We remark that the conditions on $\gH$ can be satisfied by fully connected networks equipped with the ReLU activation. In addition, we show that the push-forward measures from the approximated optimal trajectories converge to the ground truth in distribution. This ensures the validity of the learned evolution of the agent population.

\begin{thm}\label{thm_conv_measure}
Let $P_0$ be absolutely continuous with a bounded density $p_0$, and $F(\vz,t)$ be the optimizer for~\eqref{var_MFG_traj} such that $F(\vz,t_i) \in \gD^2, i=1,...,K$. Then, for each $k=1,2,..., K$, there exists a sequence of normalizing flows $\{F^n_{\theta_k}\}_{n=1}^\infty \coloneqq \{f^n_{\theta_k} \circ \cdots \circ f^n_{\theta_1}\}_{n=1}^\infty \subset \mathbf{INN_{\gH\textbf{-}ACF}}$ such that 
\begin{align*}
    F^n_{\theta_k*}P_0 \xrightarrow{d} F(\cdot, t_k)_*P_0, \hspace{2mm} \text{ as } n\to \infty,
\end{align*}
where $\xrightarrow{d}$ denotes convergence in distribution.
\end{thm}
The proof is based on the Lipschitz-bounded definition of convergence in distribution. See the proof in the appendix.

\section{Numerical Experiments}
\label{sec:expertments}
This section contains numerical results into two parts. The first part showcases NFs being an effective parameterization for solving high-dimensional MFG problems. The second part showcases the use of MFG transport costs to improve the robustness of NF models for generative modeling. 

\subsection{Solving MFGs with NFs}
In the first part, we follow the settings in~\cite{ruthotto2020machine} to construct two MFG problems. They are designed so that the behavior of the MFG solutions is invariant to dimensionality, when projected onto the first two components. This property allows for qualitative evaluation through visualization. As the dimensionality increases from 2 to 100, the problem becomes more challenging in terms of the neural network's capacity and the computational complexity. We further extend our numerical experiments to multi-group path planning. In these experiments, we use spline flows as it is commonly believed that they are more expressive than coupling flows.

\subsubsection{OT in High Dimensional Spaces}

In this experiment, we devise a dynamically formulated OT problem that can be solved efficiently with classic numerical methods in $d \leq 3$~\cite{achdou2010mean, benamou2014augmented,benamou2017variational,benamou2000computational, jacobs2019solving,papadakis2014optimal,yu2021fast,cuturi2013sinkhorn}. The optimization problem is the following:
\begin{equation}\label{OT}
\begin{split}
    \inf_{p, \vv} \int_0^1 \int_\Omega p(\vx,t) &\|\vv(\vx,t)\|^2_2 \der\vx \der t\\
    s.t. \quad \partial_t p(\vx,t) + \nabla_{\vx} \cdot (p(\vx,t)\vv(\vx,t)) &= 0, \quad \vx\in \Omega, t\in [0,T]\\
    p(\vx,0) &= p_0(\vx), \quad \vx\in \Omega\\
    p(\vx,1) &= p_1(\vx), \quad \vx\in \Omega.
\end{split}
\end{equation}
Let $\mathrm{N}(\vx; \mu, \Sigma)$ denote the density function of a multivariate Gaussian with mean $\mu \in \mathbb{R}^d$ and covariance $\Sigma \in \mathbb{R}^{d\times d}$. We set the initial density to be $p_0(\vx) = \mathrm{N}(\vx; 0, 0.3I) $ and the terminal density to be $p_1(\vx) = \frac{1}{8}\sum_{i=1}^8 \mathrm{N}(\vx; \mu_i, 0.3I)$, where $\mu_i = 4\cos(\frac{\pi}{4}i)e_1 + 4\sin(\frac{\pi}{4}i)e_2, \forall i=1, ..., 8$, and $e_1, e_2$ are the first two standard basis vectors. 

Contrary to MFGs, the OT problem imposes a strict constraint on the matching of the terminal density. To handle this, we cast the OT as a variational MFG problem by penalizing the mismatch of the terminal density using the KL divergence, and we parameterize the MFG with the coupling version of the neural spline flow (NSF-CL)~\cite{NSF}, which employs flexible rational-quadratic coupling functions to learn expressive transformations. Overall, we find NSF-CL to be a robust architecture for solving different kinds of MFGs.

\begin{figure}[h]
\centering
\begin{minipage}{0.25\linewidth}
\includegraphics[width=1\linewidth]{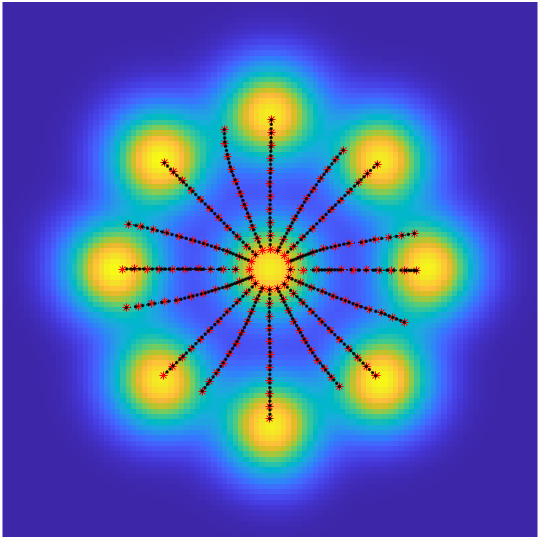}\\
\end{minipage}\hfill
\begin{minipage}{0.25\linewidth}
\includegraphics[width=1\linewidth]{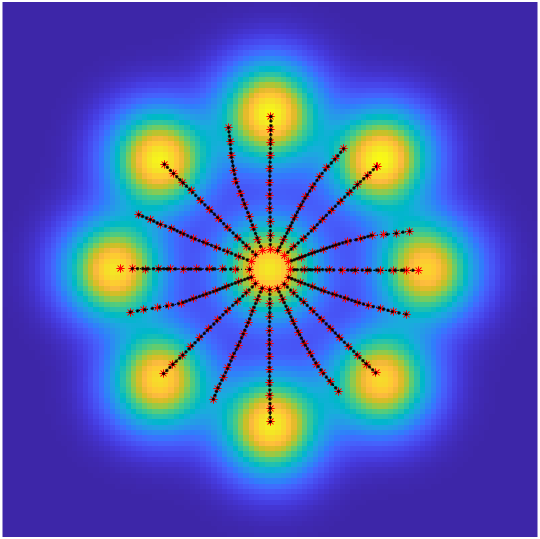}\\
\end{minipage}\hfill
\begin{minipage}{0.25\linewidth}
\includegraphics[width=1\linewidth]{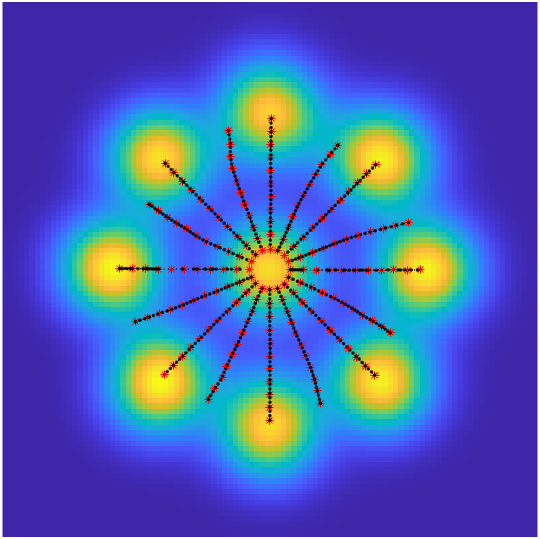}\\
\end{minipage}\hfill
\begin{minipage}{0.25\linewidth}
\includegraphics[width=1\linewidth]{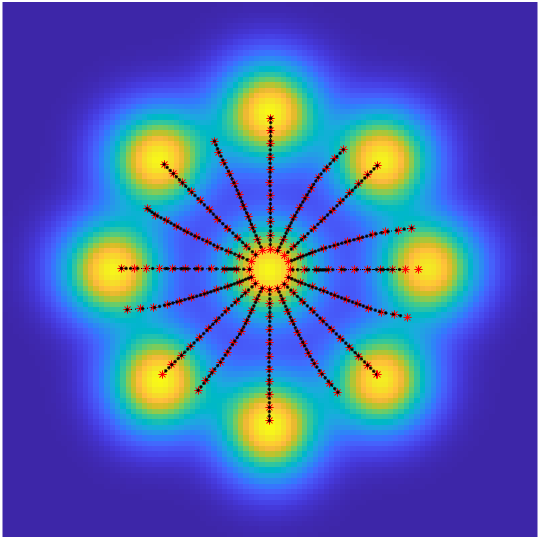}\\
\end{minipage}\hfill
\vskip -10pt
\caption{Left to right: samples of 16 computed trajectories in 2D, 10D, 50D, and 100D.}
\label{fig:OT_trajectory}
\end{figure}

\begin{figure}[h]
\centering
\begin{minipage}{0.16\linewidth}
\includegraphics[width=1\linewidth]{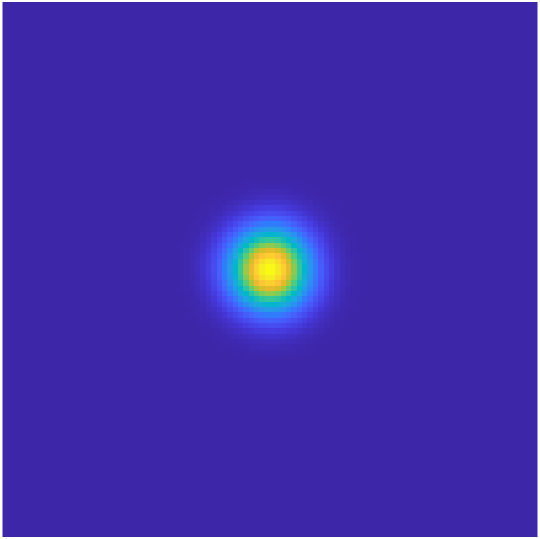}\\
\centering {$P_0$}
\end{minipage}\hfill
\begin{minipage}{0.16\linewidth}
\includegraphics[width=1\linewidth]{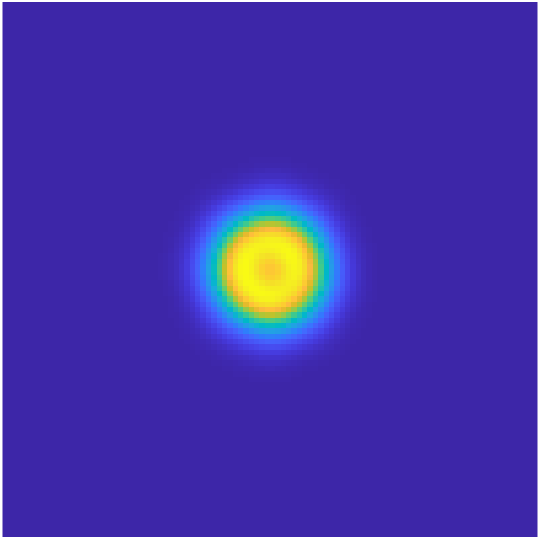}\\
\centering {$F_{1*}(P_0)$}
\end{minipage}\hfill
\begin{minipage}{0.16\linewidth}
\includegraphics[width=1\linewidth]{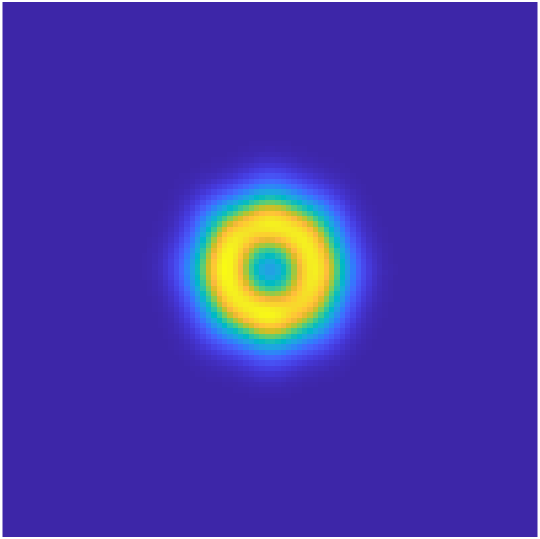}\\
\centering {$F_{2*}(P_0)$}
\end{minipage}\hfill
\begin{minipage}{0.16\linewidth}
\includegraphics[width=1\linewidth]{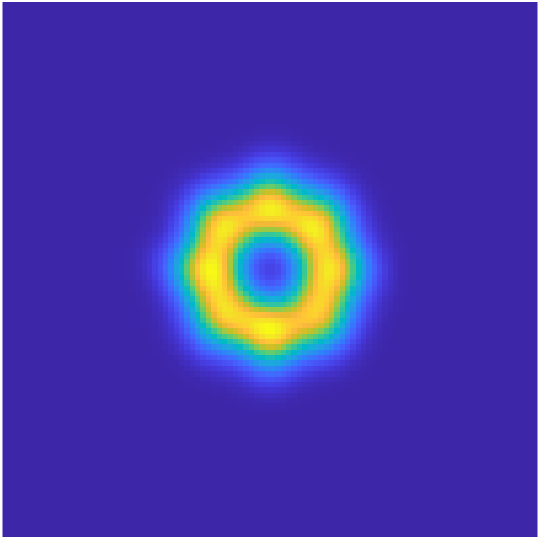}\\
\centering {$F_{3*}(P_0)$}
\end{minipage}\hfill
\begin{minipage}{0.16\linewidth}
\includegraphics[width=1\linewidth]{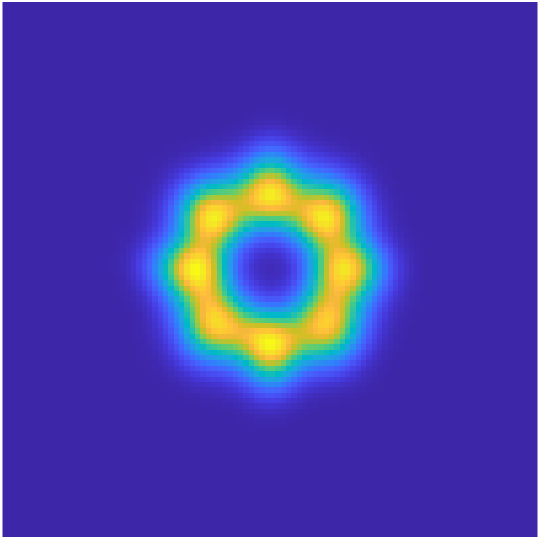}\\
\centering {$F_{4*}(P_0)$}
\end{minipage}\hfill
\begin{minipage}{0.16\linewidth}
\includegraphics[width=1\linewidth]{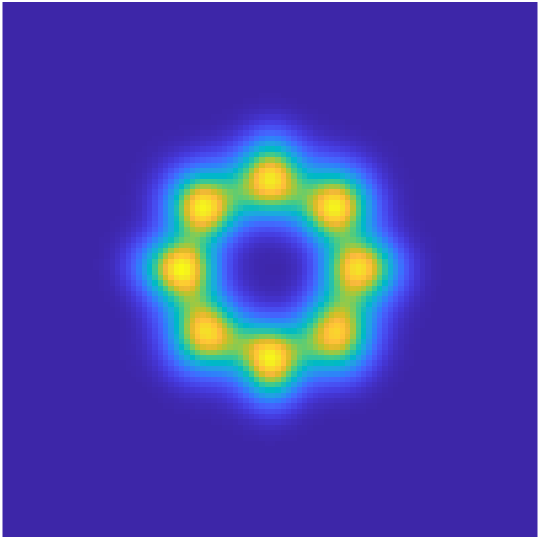}\\
\centering {$F_{5*}(P_0)$}
\end{minipage}\hfill
\begin{minipage}{0.16\linewidth}
\includegraphics[width=1\linewidth]{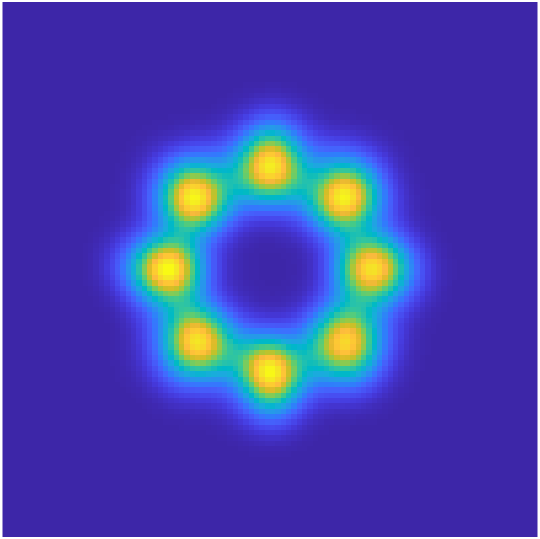}\\
\centering {$F_{6*}(P_0)$}
\end{minipage}\hfill
\begin{minipage}{0.16\linewidth}
\includegraphics[width=1\linewidth]{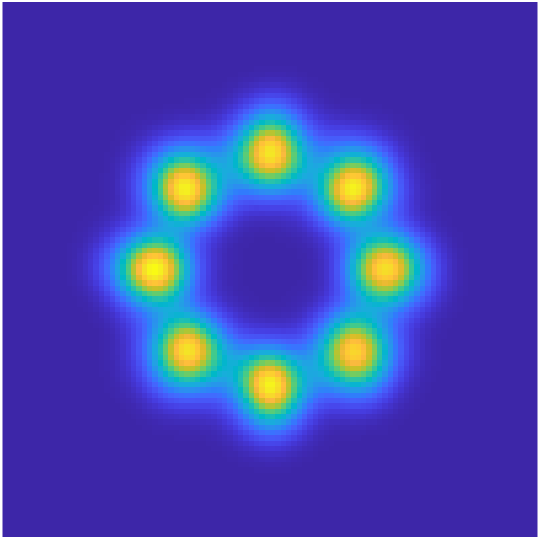}\\
\centering {$F_{7*}(P_0)$}
\end{minipage}\hfill
\begin{minipage}{0.16\linewidth}
\includegraphics[width=1\linewidth]{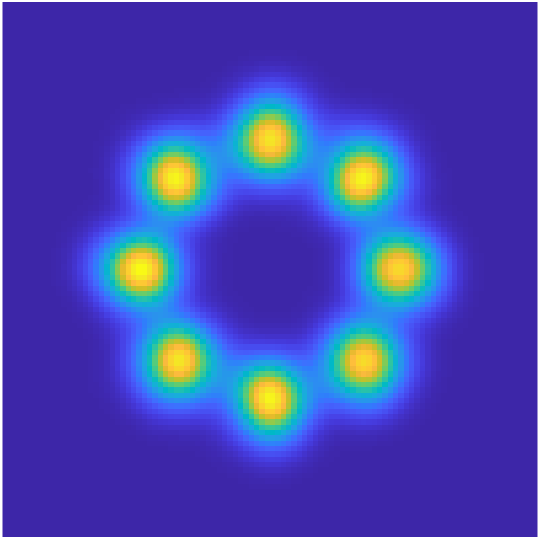}\\
\centering {$F_{8*}(P_0)$}
\end{minipage}\hfill
\begin{minipage}{0.16\linewidth}
\includegraphics[width=1\linewidth]{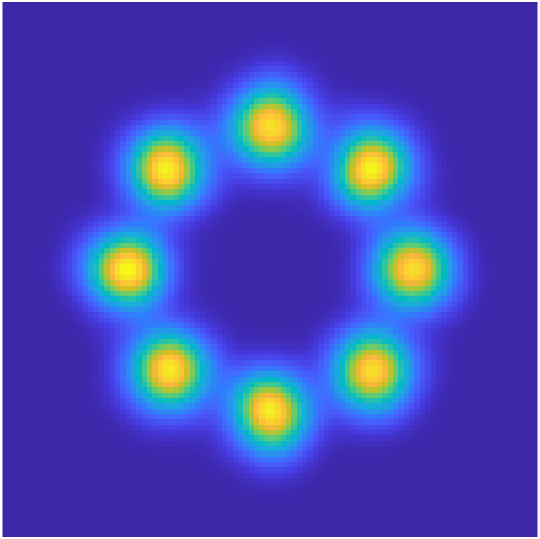}\\
\centering {$F_{9*}(P_0)$}
\end{minipage}\hfill
\begin{minipage}{0.16\linewidth}
\includegraphics[width=1\linewidth]{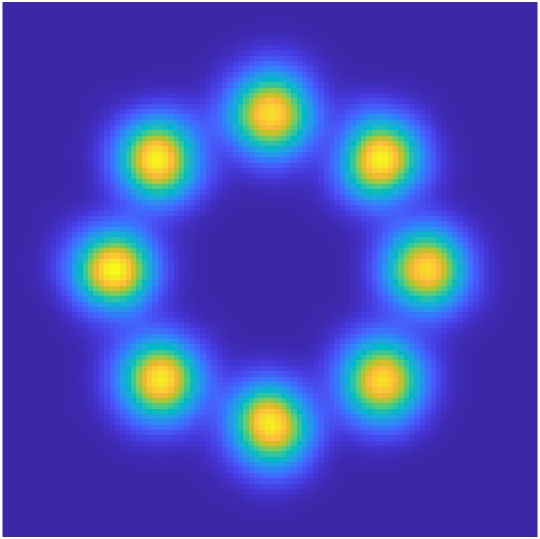}\\
\centering {$F_{10*}(P_0)$}
\end{minipage}\hfill
\begin{minipage}{0.16\linewidth}
\includegraphics[width=1\linewidth]{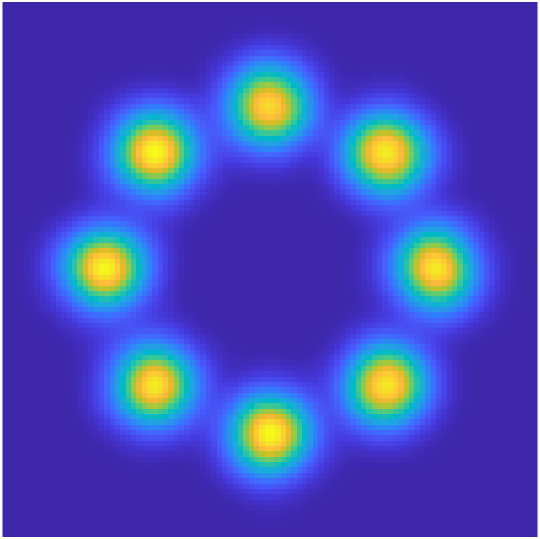}\\
\centering {$P_1$}
\end{minipage}\hfill
\caption{Evolution of the density for the OT problem~\eqref{OT} in 50D.}
\label{OT_density_evo_50D}
\end{figure}

An illustration of the agent trajectories obtained from our method is given in Figure~\ref{fig:OT_trajectory} and the evolution of the population density is displayed in Figure~\ref{OT_density_evo_50D}. We see that the NF learns to  equally partition the initial Gaussian density and transport each partition to a separate Gaussian in an approximately straight trajectory. Moreover, it does so consistently as the problem dimension grows. More experimental results can be found in the appendix.

\begin{table}[h]
\begin{center}
\begin{small}
\begin{sc}
\begin{tabular}{lccccccccccc}
\toprule
Method & Dimension & $L$ & $\gM $ & Time/iter & Number of iters\\
\midrule
 Eulerian & 2 & 9.8 & 0.1725 & - & - \\
 MFGNet  & 2 & 9.9 & 0.1402 & 2.038 & 1e3$^*$ \\
 Ours  & 2 & 10.1 & 0.0663 & 0.549 & 1e5 \\
\midrule
 MFGNet  & 10 & 10.1 & 0.1616 & 8.256 & 1e3$^*$\\
 Ours  & 10 & 10.1 & 0.0641 & 0.561 & 1e5\\
 \midrule
 MFGNet  & 50 & 10.1 & 0.1396 & 81.764 & 1e3$^*$\\
 Ours  & 50 & 10.1 & 0.0641 & 0.584 & 1e5\\
 \midrule
 MFGNet  & 100 & 10.1 & 0.1616 & 301.043 & 1e3$^*$\\
 Ours  & 100 & 10.1 & 0.0622 & 0.625 & 1e5\\
\bottomrule
\end{tabular}
\end{sc}
\end{small}
\end{center}
\caption{Comparison of methods for solving the OT problem~\eqref{OT}. $L$: transport cost; $\mathcal{M} = D_{KL}$: terminal cost. $^*$The first 500 iterations are trained with less data for a warm start.}
\label{OT_table}
\end{table}

We compute the MFG cost term and compare our results to a traditional Eulerian method~\cite{yu2021fast} as well as the MFGNet~\cite{ruthotto2020machine}, another deep learning approach based on continuous NFs. The Eulerian method is provably convergent, but its reliance on a discrete grid renders it infeasible for dimensions higher than three. For a fair comparison, we report the unweighted cost and treat the KL divergence $D_{KL}(F(\cdot, T)_*P_0 || P_1)$ as the terminal cost.

Table~\ref{OT_table} demonstrates that the NF-based parameterization can solve MFG problems up to 100 dimensions while maintaining a consistent transport cost. Compared to the MFGNet, our approach yields a similar transport cost but a much lower terminal cost, entailing a more accurate density matching at the final time. We attribute the improvement in terminal matching to the trajectory-based parameterization of the MFG. By encoding the continuity equation~\eqref{HJB} in the flow architecture, our approach sidesteps the errors incurred in the numerical schemes used in the MFGNet, when approximating the solution. 

We remark that our model is optimized with ADAM~\cite{ADAM}  while the MFGNet is trained with BFGS. The superlinear convergence of BFGS allows the MFGNet to converge with fewer iterations, at the price of higher costs per update. In addition, our implementation is based on Pytorch, which effectively leverages GPU parallelism to achieve an optimal runtime scaling, up to 100 dimensions. While being slower than MFGNet in lower dimensions, our approach uses only one third of the time in $d=100$, since its runtime is nearly constant in different dimensions.


\subsubsection{Crowd Motion}
\label{subsec:CrowdMotion}
In this experiment, we consider an MFG problem with nonzero interaction costs. Following the general setting in~\eqref{var_MFG_traj}, we use a two-part interaction cost $\mathcal{I}$:
\begin{equation}\label{crowd_motion_MFG}
\begin{split}
    \mathcal{I}_P(F(\cdot, t)_*p_0) &\coloneqq \lambda_P\int_\Omega Q(\vx) \der (F(\cdot, t)_*P_0)(\vx) = \lambda_P\int_\Omega Q(F(\vx,t)) p_0(\vx)\der \vx\\
    \mathcal{I}_E(F(\cdot, t)_*p_0) &\coloneqq \lambda_E\int_\Omega (\log p_0(\vx) - \log |\det \nabla  F(\vx,t)|) p_0(\vx) \der x\\
    \mathcal{I} (F(\cdot, t)_*p_0) &\coloneqq \lambda_{\mathcal{I}}( \mathcal{I}_P(F(\cdot, t)_*p_0) +  \mathcal{I}_E(F(\cdot, t)_*P_0)).
\end{split}
\end{equation}
Here, $Q(\vx)$ acts as an obstacle that incurs a cost for an agent that passes through it. It is set to be the density of a bivariate Gaussian centered at the origin, with a magnitude of 50:
\begin{equation}\label{Q}
\begin{split}
    Q(\vx) &\coloneqq 50 \cdot \mathrm{N} (\vx; 0, \text{diag}(1,0.5)).
\end{split}
\end{equation}

In dimensions higher than two, we evaluate the first two components of $F(\vx,t)$ on $Q$ to compute $\mathcal{I}_P$. The initial and the terminal densities are chosen as Gaussians: $p_0(\vx) = \mathrm{N}(\vx;3e_2, 0.3I), p_1(\vx) = \mathrm{N}(\vx;-3e_2, 0.3I)$; and the terminal cost is identical to that in the OT experiment~\eqref{OT}.  Intuitively, this problem aims at transporting an initial Gaussian density to a different location, while avoiding an obstacle placed at the origin. Moreover, the optimal trajectories are invariant to the dimensionality on the first two components, allowing us to visualize the dynamics in high dimensions.

Similar to the transport cost $L$, we perform a discretization on~\eqref{crowd_motion_MFG}, based on the NF parameterization. Note that $\mathcal{I}$ is integrated in time in the MFG problem~\eqref{var_MFG_traj}. Thus, discretizing the integral with the right-point rule yields
\begin{align}
    \int_0^1 \mathcal{I} (F(\cdot, t)_*p_0) \der t &= \mathbb{E}_{\vx\sim P_0} \left\{\lambda_{\mathcal{I}} \int_0^1 [\lambda_P Q(F(\vx,t)) + \lambda_E (\log p_0(\vx) - \log |\det \nabla  F(\vx,t)|] \der t \right\}\\
    &\approx \frac{\lambda_{\mathcal{I}}}{K} \mathbb{E}_{\vx\sim P_0} \left\{ \sum_{i=1}^K [\lambda_P Q(F_i(\vx)) + \lambda_E (\log p_0(\vx) - \log |\det \nabla  F_i(\vx_{i-1})|)] \right\},
\end{align}
where $F_i = f_i \circ f_{i-1} \circ ... \circ f_1, \forall i = 1, ..., K$ and $\vx_0 \coloneqq \vx, \vx_j \coloneqq f_j(\vx_{j-1}), \forall j=1, ..., K-1$. Similar to before, our discretization is consistent with the $O(\frac{1}{K})$ global error. In practice, we use fourth-order forward difference for the $\partial_t F(\vz,t)$ term in $L$ and approximate all integrals with Simpson's rule. The overall numerical scheme converges to the continuous MFG with $O(\frac{1}{K^4})$ global error.

\begin{figure}[h]
\centering
\begin{minipage}{0.25\linewidth}
\includegraphics[width=1\linewidth]{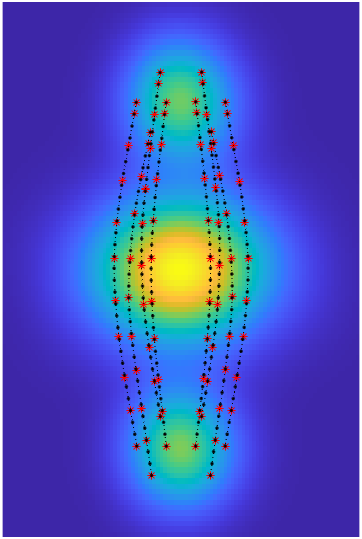}\\
\end{minipage}\hfill
\begin{minipage}{0.25\linewidth}
\includegraphics[width=1\linewidth]{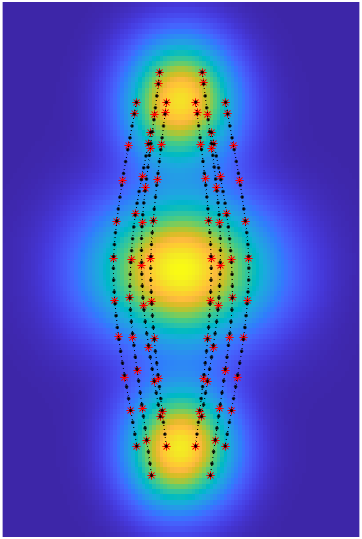}\\
\end{minipage}\hfill
\begin{minipage}{0.25\linewidth}
\includegraphics[width=1\linewidth]{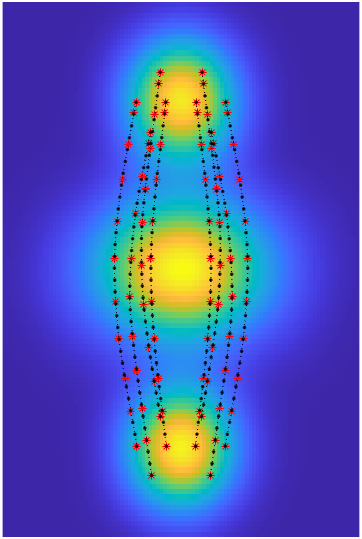}\\
\end{minipage}\hfill
\begin{minipage}{0.25\linewidth}
\includegraphics[width=1\linewidth]{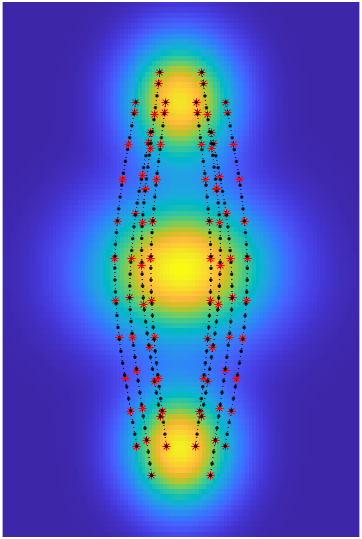}\\
\end{minipage}\hfill
\caption{Left to right: samples of the computed trajectories in 2D, 10D, 50D, and 100D. The top and bottom densities are $p_0$ and $p_1$, respectively, while the middle density denotes the obstacle.}
\label{crowd_traj_allD}
\end{figure}

\begin{table}[h]
\caption{Comparison of methods for solving the crowd motion problem~\eqref{crowd_motion_MFG}. $L$: transport cost; $\mathcal{I}$: interaction cost; $\mathcal{M}$: terminal cost. $^*$The first 500 iterations are trained with less data for a warm start.}
\label{crowd_motion_table}
\begin{center}
\begin{small}
\begin{sc}
\begin{tabular}{lccccccccccc}
\toprule
Method & Dimension & $L$ & $\mathcal{I}$ & $\mathcal{M}$ & Time/iter & Number of iters\\
\midrule
 Eulerian  & 2 & 31.8 & 2.27 & 0.1190 & - & - \\
 MFGNet  & 2 & 33.0 & 2.29 & 0.1562 & 4.122 & 1e3$^*$ \\
 Ours  & 2 & 32.8 & 2.19 & 0.0417 & 0.559 & 1e5\\
\midrule
 MFGNet  & 10 & 33.0 & 2.29 & 0.1502 & 17.205 & 1e3$^*$ \\
 Ours  & 10 & 32.9 & 2.13 & 0.0436 & 0.568 & 1e5\\
 \midrule
 MFGNet  & 50 & 33.0 & 1.91 & 0.1440 & 134.938 & 1e3$^*$ \\
 Ours  & 50 & 32.9 & 1.82 & 0.0381 & 0.581 & 1e5\\
 \midrule
 MFGNet  & 100 & 33.0 & 1.49 & 0.2000 & 241.727 & 1e3$^*$ \\
 Ours  & 100 & 33.0 & 1.36 & 0.0464 & 0.625 & 1e5\\
\bottomrule
\end{tabular}
\end{sc}
\end{small}
\end{center}
\end{table}

\begin{figure}[h]
\centering
\begin{minipage}{0.33\linewidth}
\centering
\includegraphics[width=.9\linewidth]{density_traj_NSF_CL_10D_N=1M_identicalW.png}
\end{minipage}\hfill
\begin{minipage}{0.33\linewidth}
\centering
\includegraphics[width=.9\linewidth]{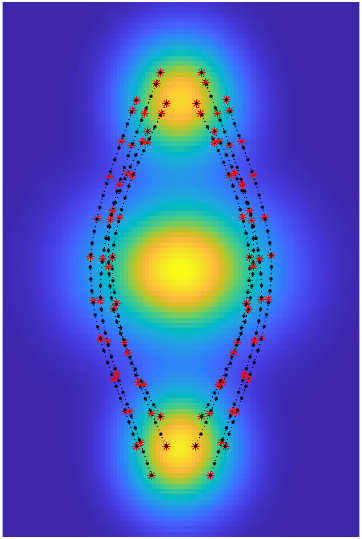}
\end{minipage}\hfill
\begin{minipage}{0.33\linewidth}
\centering
\includegraphics[width=.9\linewidth]{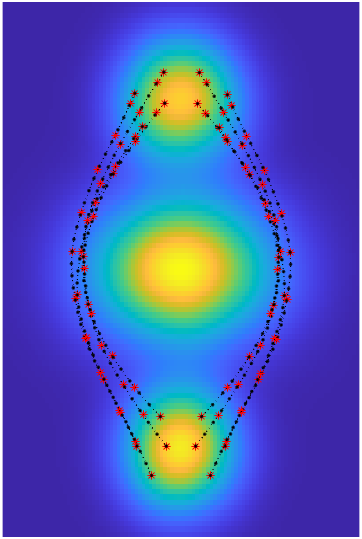}
\end{minipage}\hfill
\caption{Computed trajectories for the crowd motion problem~\eqref{crowd_motion_MFG} in 10D, projected onto the first two dimensions. From left to right: minor ($\lambda_{\mathcal{I}}=0.2$), moderate ($\lambda_{\mathcal{I}}=0.5$), and strong ($\lambda_{\mathcal{I}}=1$) penalty on conflicts with the obstacle. The left figure corresponds to the case reported in Table~\ref{crowd_motion_table}.}
\label{crowd_traj}
\end{figure}

Similar to the OT experiment, we use the same NSF-CL flow~\cite{NSF} to solve the crowd motion problem from 2 to 100 dimensions and compare the computed cost terms to those of MFGNet in Table~\ref{crowd_motion_table}. Compared to the baseline, our method yields comparable transport cost but lower interaction and terminal costs, suggesting a more accurate characterization of the agents' behavior in avoiding the conflicts with the obstacle, as well as a more authentic matching of the final population distribution. Similar to the MFGNet, our computed $L, \mathcal{M}$ costs are approximately invariant with respect to the dimensionality, as further validated by the visualizations of the sampled trajectories shown in Figure~\ref{crowd_traj_allD}.
The agents plan their paths to circumvent the obstacle; moreover, trajectories of symmetric initial points are symmetric as well, which is expected from the problem setup. In Figure~\ref{crowd_traj}, we also show the optimal trajectories at different levels of aversion preference, further validating the robustness of the method. From left to right, the agents are farther away from the obstacle, while maintaining approximately symmetric trajectories and accurate terminal matching. The associated costs for each scenario are given in the appendix, together with more numerical results.

\subsubsection{Multi-Group Path Planning}
\label{subsec:pathplanning}

In the above crowd motion experiment, we demonstrated that the interaction cost can be used to devise the optimal trajectories for a single population, in the presence of obstacles. Here, we take one step further and investigate path planning when multiple populations are moving simultaneously. Therein, the interaction cost is used to discourage  conflicts among different populations. 

We consider a generalized MFG problem that accommodates multiple densities; this problem is also explored in~\cite{MFG_city, MFG_drones_GAN}. We consider $N_p$ populations $\{p^i(\vx,t)\}_{i=1}^{N_p}$ and their associated trajectories $\{F^i(\vx,t)\}_{i=1}^{N_p}$, each of which is parameterized by an NF. Every population maintains an independent transport cost $ \lambda_L \int_0^1 \int_\Omega p_0(\vx) \|\partial_t F^i(\vx,t)\|_2^2 \der\vx \der t$ and an independent terminal cost $\mathcal{M}_i = \lambda_{\mathcal{M}} \mathcal{D}_{KL} (P^i_1||F^i(\cdot, 1)_*P^i_0)$, where $P^i_1, P^i_0$ are the initial and terminal distributions for the $i$-th population, respectively.

To keep the setting simple yet illustrative, we consider an inter-group interaction cost that encourages different populations to avoid each other in their paths. In a real-life scenario, this can be thought of as path planning for multiple groups of drones, and the desired outcome is for them to arrive at their desired locations through short paths while minimizing collisions. Formally, the interaction cost is modeled by a Gaussian kernel $\mathcal{I} = \sum_{j\ne i} \int_{\mathbb{R}^d}\int_{\mathbb{R}^d} e^{-\frac{1}{2} \|\vx-\vy\|^2_2} \der P^i(\vx,t) \der P^j(\vy,t)$, which decays exponentially as the radial distance between two populations increases.

\begin{figure}[H]
\centering
\begin{minipage}{0.33\linewidth}
\includegraphics[width=.9\linewidth]{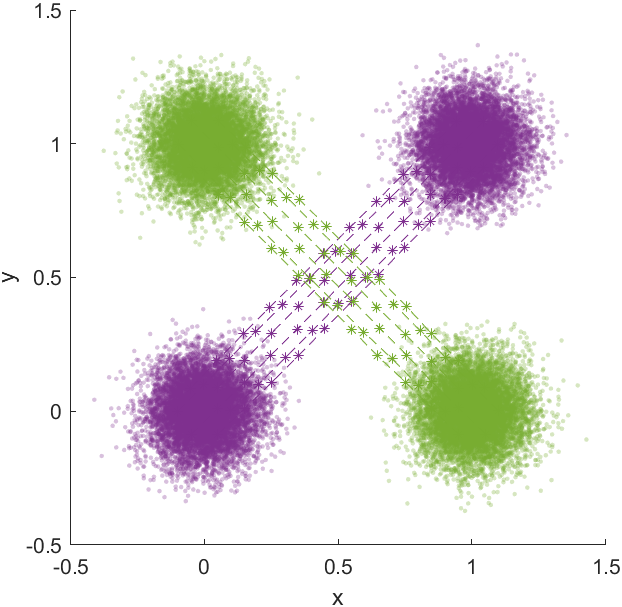}\\
\end{minipage}\hfill
\begin{minipage}{0.33\linewidth}
\includegraphics[width=.9\linewidth]{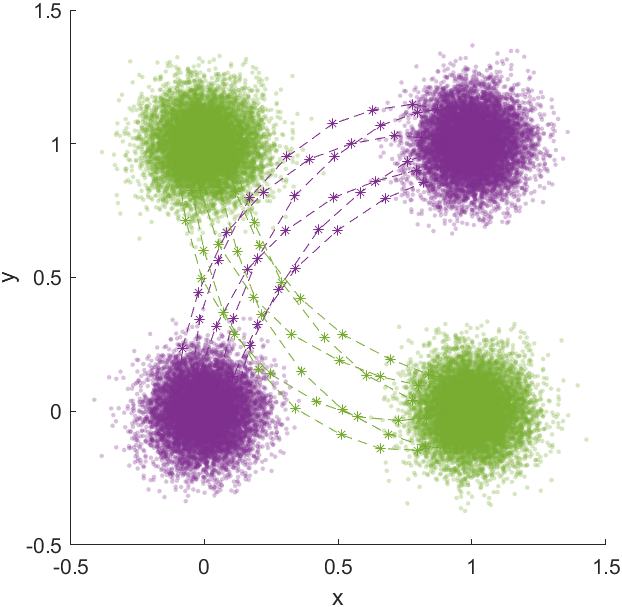}\\
\end{minipage}\hfill
\begin{minipage}{0.33\linewidth}
\includegraphics[width=.9\linewidth]{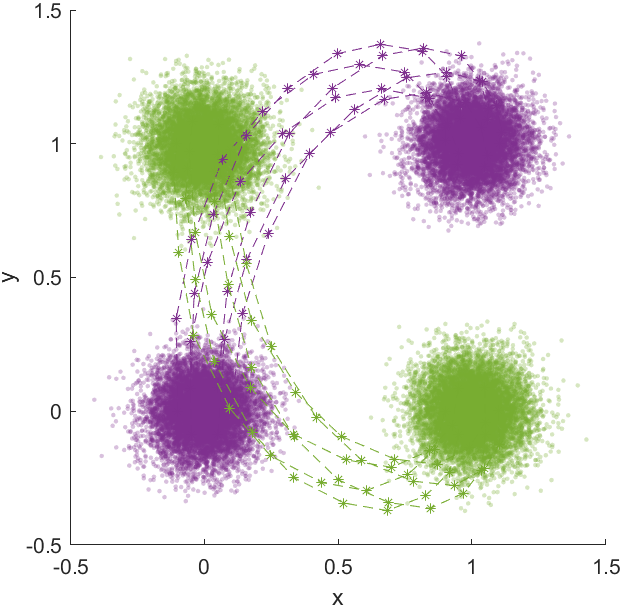}\\
\end{minipage}\hfill
\begin{minipage}{0.33\linewidth}
\includegraphics[width=.95\linewidth]{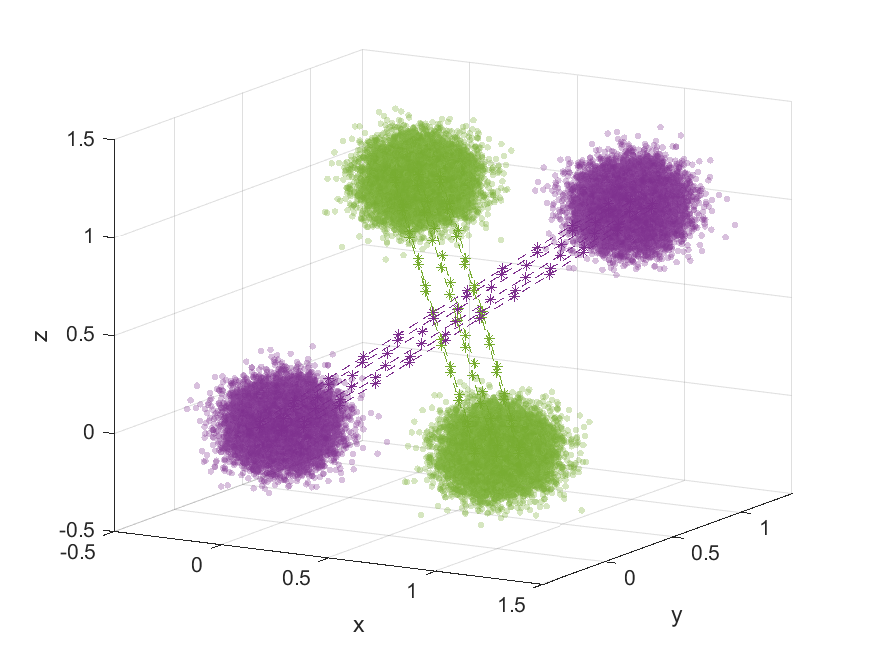}\\
\end{minipage}\hfill
\begin{minipage}{0.33\linewidth}
\includegraphics[width=.95\linewidth]{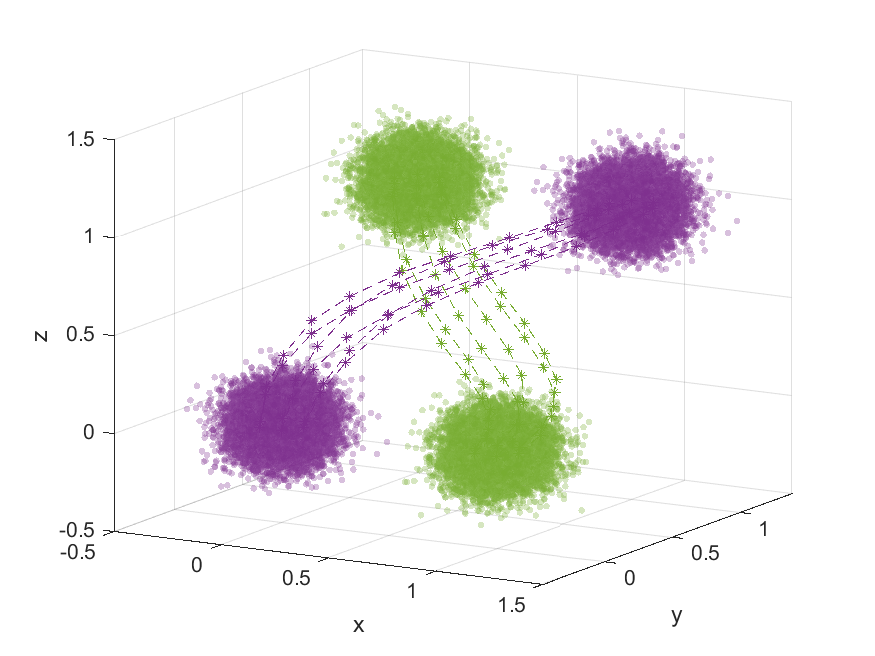}\\
\end{minipage}\hfill
\begin{minipage}{0.33\linewidth}
\includegraphics[width=.95\linewidth]{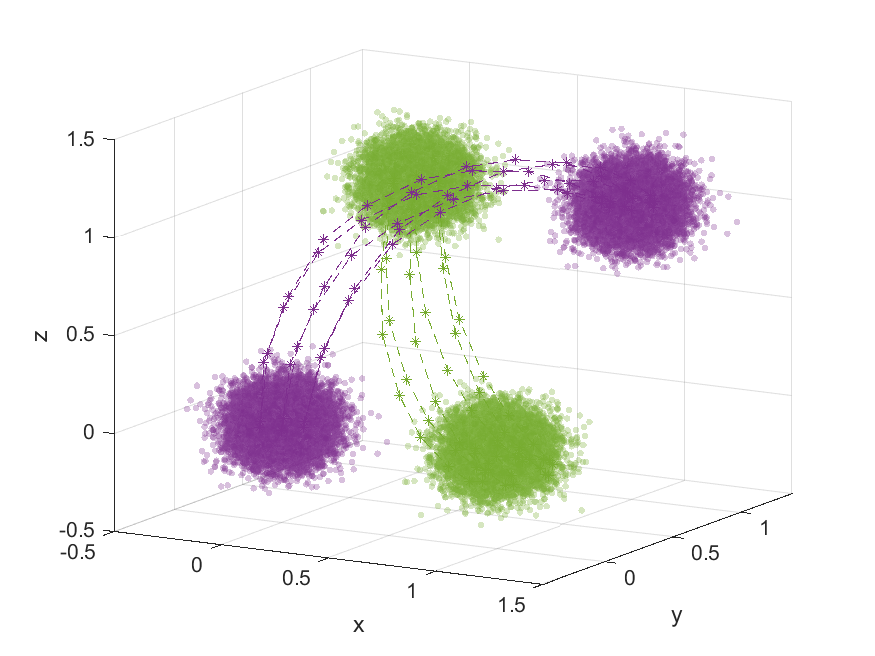}\\
\end{minipage}\hfill
\begin{minipage}{0.33\linewidth}
\includegraphics[width=.95\linewidth]{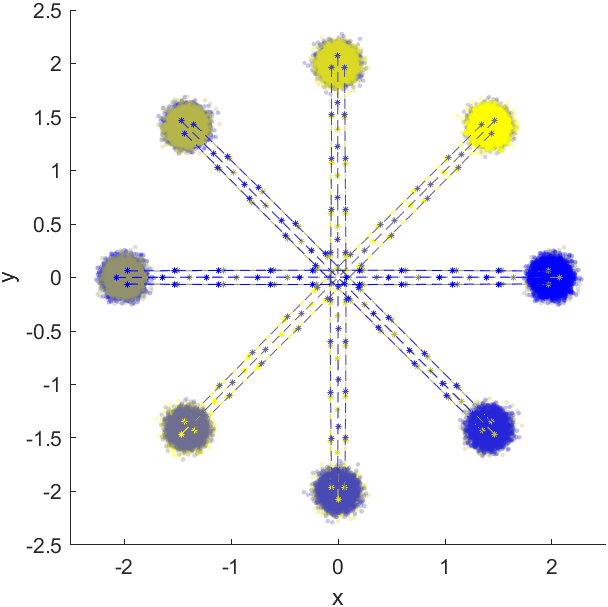}\\
\end{minipage}\hfill
\begin{minipage}{0.33\linewidth}
\includegraphics[width=.95\linewidth]{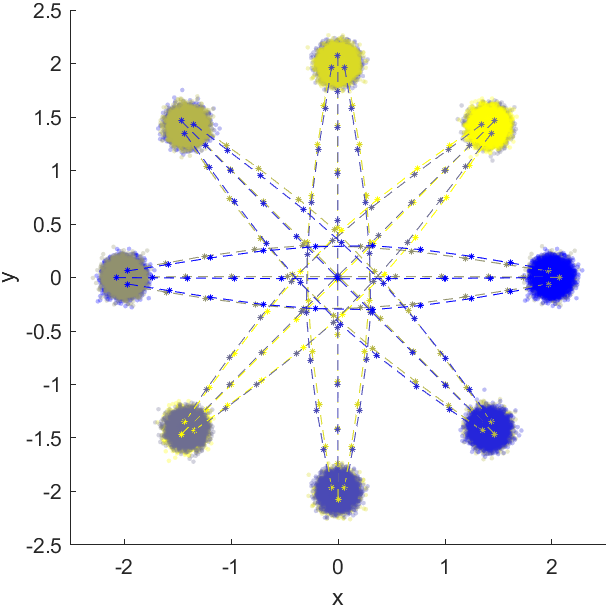}\\
\end{minipage}\hfill
\begin{minipage}{0.33\linewidth}
\includegraphics[width=.95\linewidth]{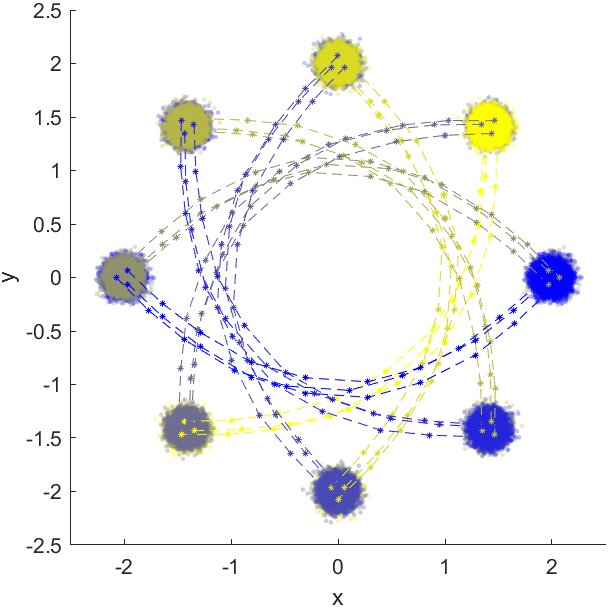}\\
\end{minipage}\hfill
%
\caption{Computed trajectories for the multi-group path planning problem. From left to right: zero, moderate, and strong penalty on conflicts between populations. Colors indicate different populations. From top to bottom: two populations in 2D, two populations in 3D, and eight populations in 2D.}
\label{drone_traj}
\end{figure}

In Figure~\ref{drone_traj}, we show results for problems in two and three dimensions. 
For the 2D experiment with two populations, we choose $p^1_0(\vx,t) = \mathrm{N}((0,0), 0.01I), p^1_1(\vx,t) = \mathrm{N}((1,1), 0.01I)$ as the initial and terminal density for the first population, and $p^2_0(\vx,t) = \mathrm{N}((1,0), 0.01I), p^2_1(\vx,t) = \mathrm{N}((0,1), 0.01I)$ for the second population. For the 3D experiment with two populations, we have $P^1_0(\vx,t) = \mathrm{N}((0,0,0), 0.01I), p^1_1(\vx,t) = \mathrm{N}((1,1,1), 0.01I), p^2_0(\vx,t) = \mathrm{N}((1,0,0), 0.01I), p^1_2(\vx,t) = \mathrm{N}((0,1,1),0.01I)$. For the 2D experiment with eight populations, the densities are $p_i(\vx) = \mathrm{N}(\vx; \mu_i, 0.05I)$, where $\mu_i = 4\cos(\frac{\pi}{4}i)e_1 + 4\sin(\frac{\pi}{4}i)e_2, \forall i=1, ..., 8$. Each population uses the density mirrored across the origin as its destination. 

In these problems, the initial and the terminal densities are configured to create a full collision if the interaction term is absent. As the weights on the interaction cost increase, the populations coordinate a delicate dance to avert each other during travel. Note that the computed trajectories are approximately symmetric, which is expected according to the symmetric nature of the problem setup. We summarize the computed unweighted costs in Table~\ref{multi_group_cost}. As higher weights are placed on the interaction cost, the populations spend more effort maneuvering and hence the transport cost increases.

\begin{table}[h]
\caption{Results of solving the multi-group path planning problem by using different interaction wrights. $\lambda_{\mathcal{I}}$: weight for $\mathcal{I}$; $L$: transport cost; $\mathcal{I}$: interaction cost; $\mathcal{M}$: terminal cost.}
\label{multi_group_cost}
\begin{center}
\begin{small}
\begin{sc}
\begin{tabular}{lccc|cccc|cccc}
\toprule
 & \multicolumn{3}{c|}{2D, 2 populations} & \multicolumn{4}{c|}{3D, 2 populations} & \multicolumn{4}{c}{2D, 8 populations}\\
 \midrule
$\lambda_{\mathcal{I}}$ & $L$ & $\mathcal{I}$ & $\mathcal{M}$ & $\lambda_{\mathcal{I}}$ & $L$ & $\mathcal{I}$ & $\mathcal{M}$ & $\lambda_{\mathcal{I}}$ & $L$ & $\mathcal{I}$ & $\mathcal{M}$\\
\midrule
0   & 3.9755 & 0.8992  & 0.0038  & 0 & 5.9899 & 0.8897 & 0.0023 & 0 & 120.8807 & 12.5090 & 0.0067\\
3   & 5.5874 & 0.6556  & 0.0034 & 3 & 7.2835 & 0.6574 & 0.0017 & 1 & 125.7670 &  10.6740 & 0.0086 \\
5   & 9.1048 & 0.4690  & 0.0099  & 5 & 10.8889 & 0.4677 & 0.0047 & 3 & 164.5073 & 6.0951 & 0.0089 \\
\bottomrule
\end{tabular}
\end{sc}
\end{small}
\end{center}
\end{table}

\subsection{Improving NF Training with MFGs}

The standard training of NFs only considers minimizing the KL divergence (the terminal cost in the MFG), without transport and interaction costs. Therefore, intermediate distributions may be meaningless artifacts of the flow parameterization and they often come with enormous distortion. By viewing NFs as a parameterization of the discretized MFG, it is natural to introduce the transport cost into NF training. We showed in Theorem~\ref{cont_disc_MFG_soln} that the training problem admits a unique mapping as its solution. In other words, the problem is no longer ill-posed. 

More interestingly, the unique transport-efficient mapping tends to be well-behaved, in terms of the Lipschitz bound illustrated in~\eqref{eqn:LipRegu}. This behavior is closely related to the model's complexity and the robustness of the neural network~\cite{Lip_gen}. As a result, the weight on the trajectory regularization serves as an effective means of controlling the Lipschitz bound. With appropriate choices of the weight, we observe in a variety of tabular and image datasets that the trained model has a better generalization performance.

\subsubsection{Synthetic Datasets}

\begin{figure}[h]
\centering

\begin{minipage}{0.142\linewidth}
\includegraphics[width=1\linewidth]{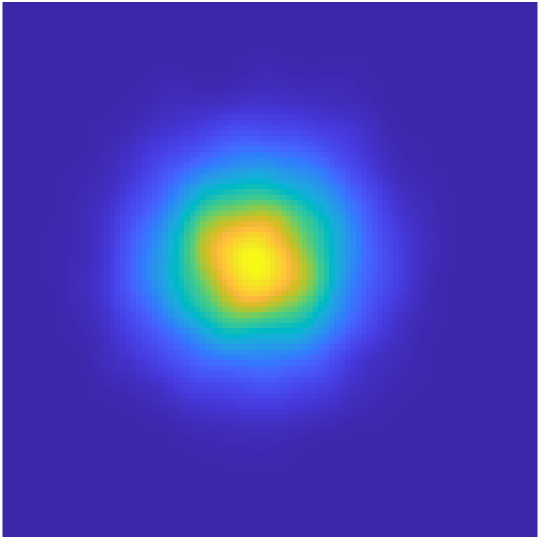}\\
\end{minipage}\hfill
\begin{minipage}{0.142\linewidth}
\includegraphics[width=1\linewidth]{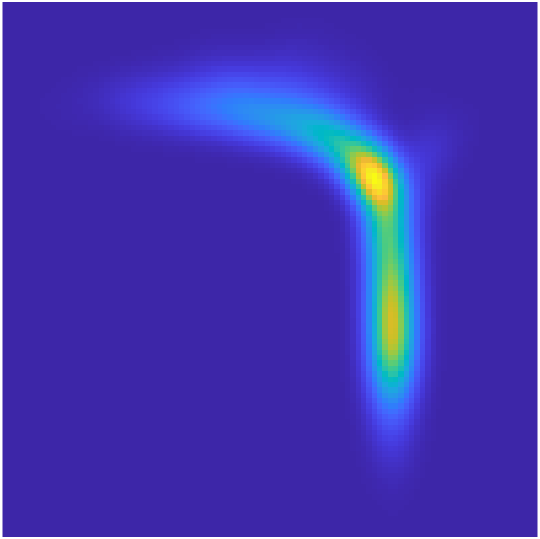}\\
\end{minipage}\hfill
\begin{minipage}{0.142\linewidth}
\includegraphics[width=1\linewidth]{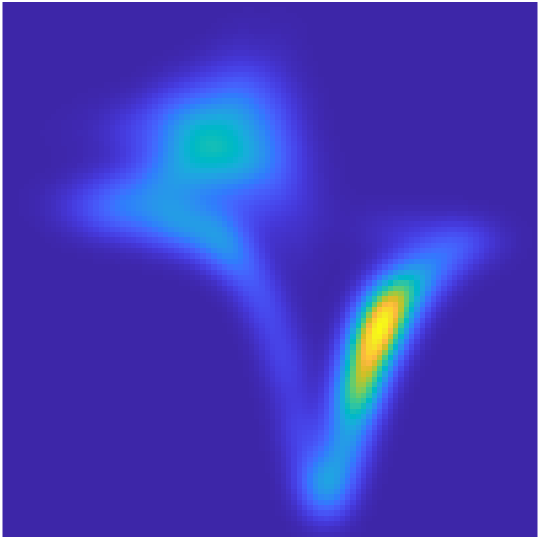}\\
\end{minipage}\hfill
\begin{minipage}{0.142\linewidth}
\includegraphics[width=1\linewidth]{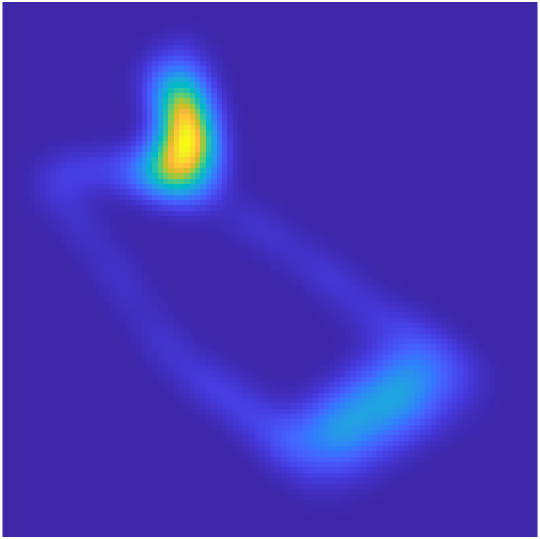}\\
\end{minipage}\hfill
\begin{minipage}{0.142\linewidth}
\includegraphics[width=1\linewidth]{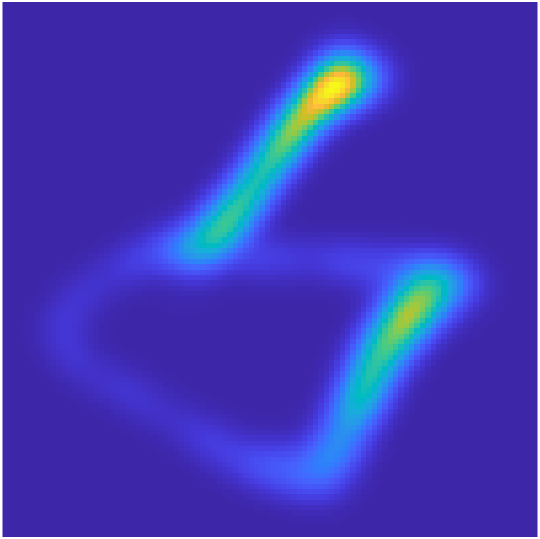}\\
\end{minipage}\hfill
\begin{minipage}{0.142\linewidth}
\includegraphics[width=1\linewidth]{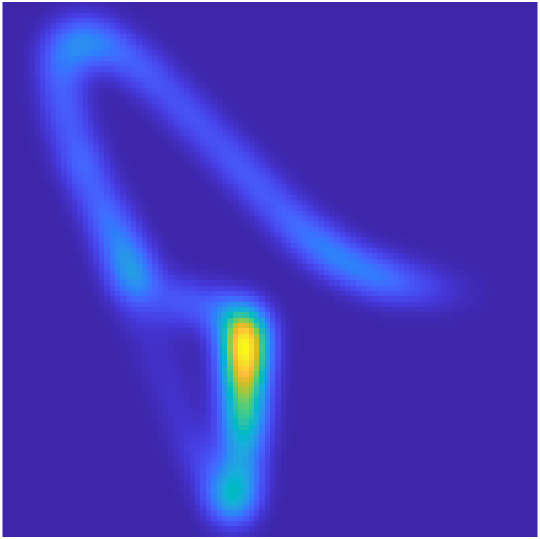}\\
\end{minipage}\hfill
\begin{minipage}{0.142\linewidth}
\includegraphics[width=1\linewidth]{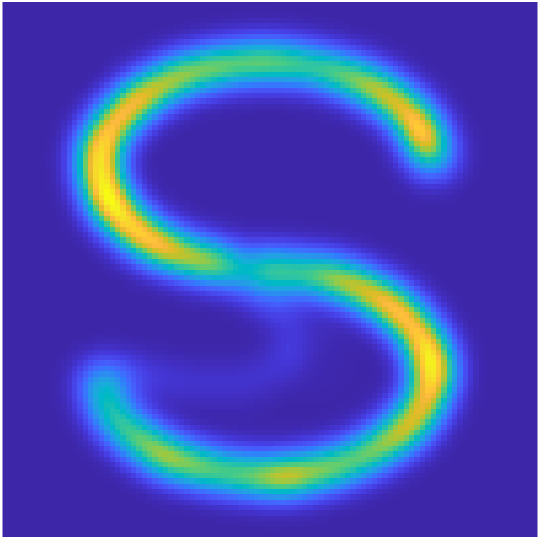}\\
\end{minipage}\hfill

\begin{minipage}{0.142\linewidth}
\includegraphics[width=1\linewidth]{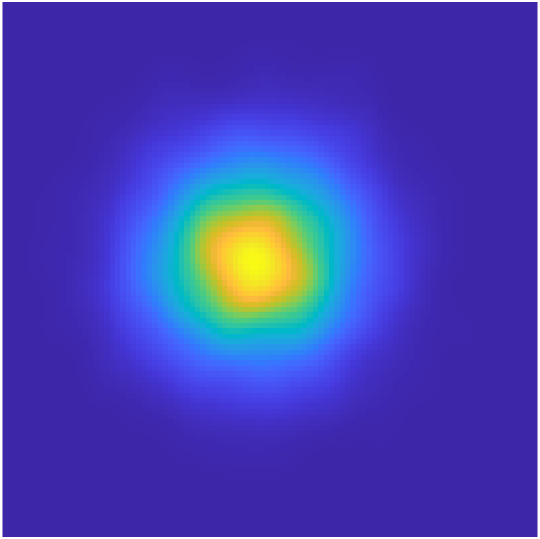}\\
\end{minipage}\hfill
\begin{minipage}{0.142\linewidth}
\includegraphics[width=1\linewidth]{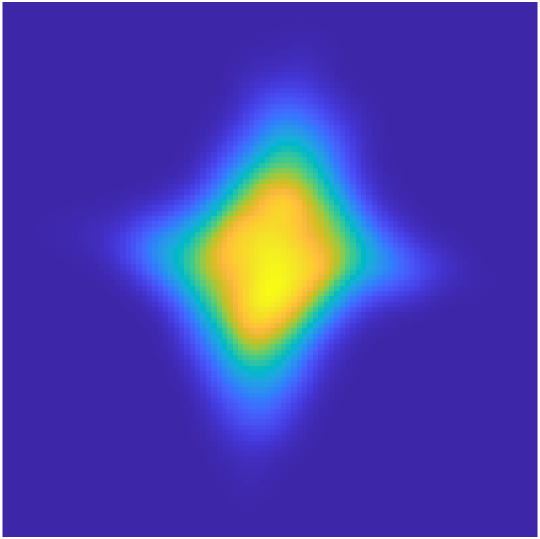}\\
\end{minipage}\hfill
\begin{minipage}{0.142\linewidth}
\includegraphics[width=1\linewidth]{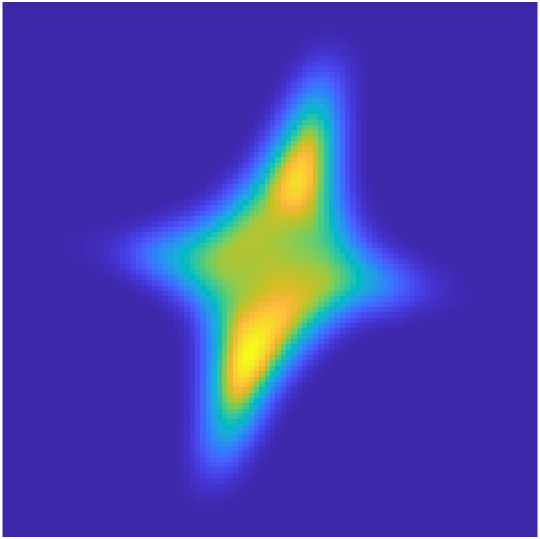}\\
\end{minipage}\hfill
\begin{minipage}{0.142\linewidth}
\includegraphics[width=1\linewidth]{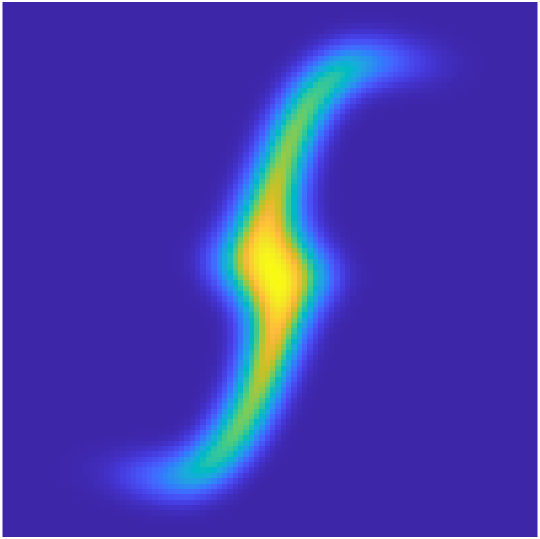}\\
\end{minipage}\hfill
\begin{minipage}{0.142\linewidth}
\includegraphics[width=1\linewidth]{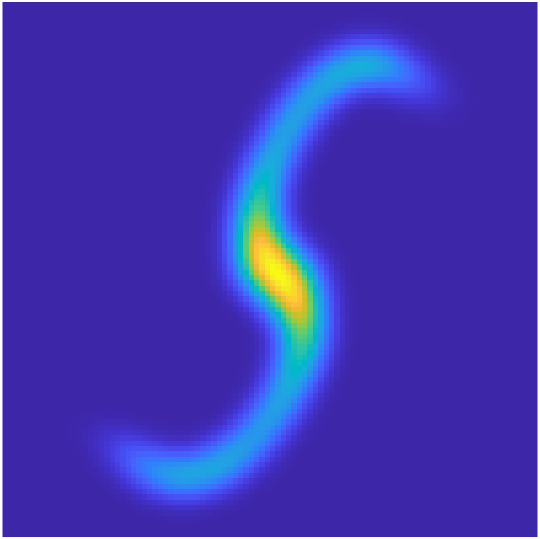}\\
\end{minipage}\hfill
\begin{minipage}{0.142\linewidth}
\includegraphics[width=1\linewidth]{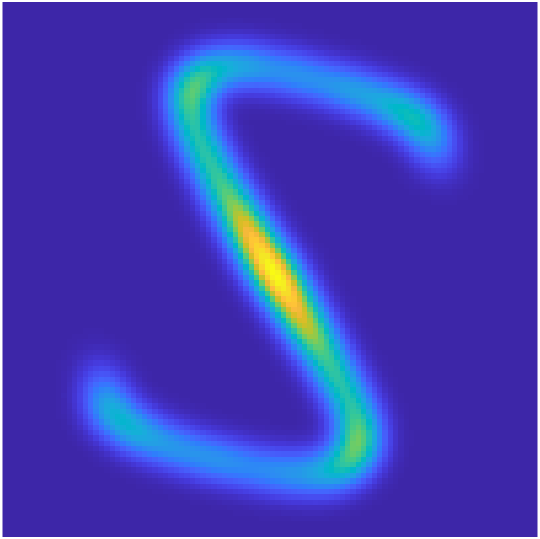}\\
\end{minipage}\hfill
\begin{minipage}{0.142\linewidth}
\includegraphics[width=1\linewidth]{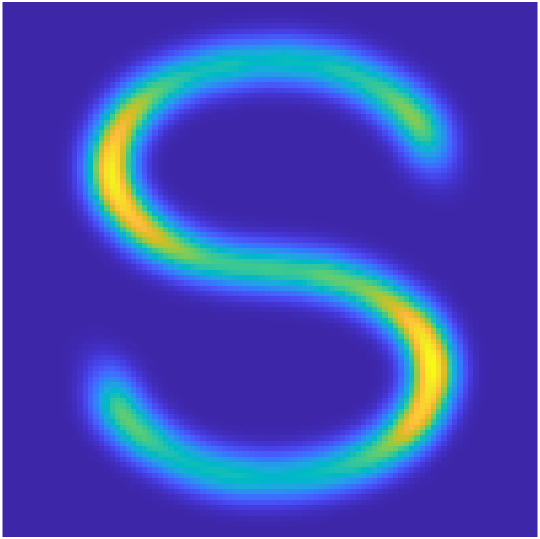}\\
\end{minipage}\hfill

\begin{minipage}{0.49\linewidth}
\centering
\includegraphics[width=.9\linewidth]{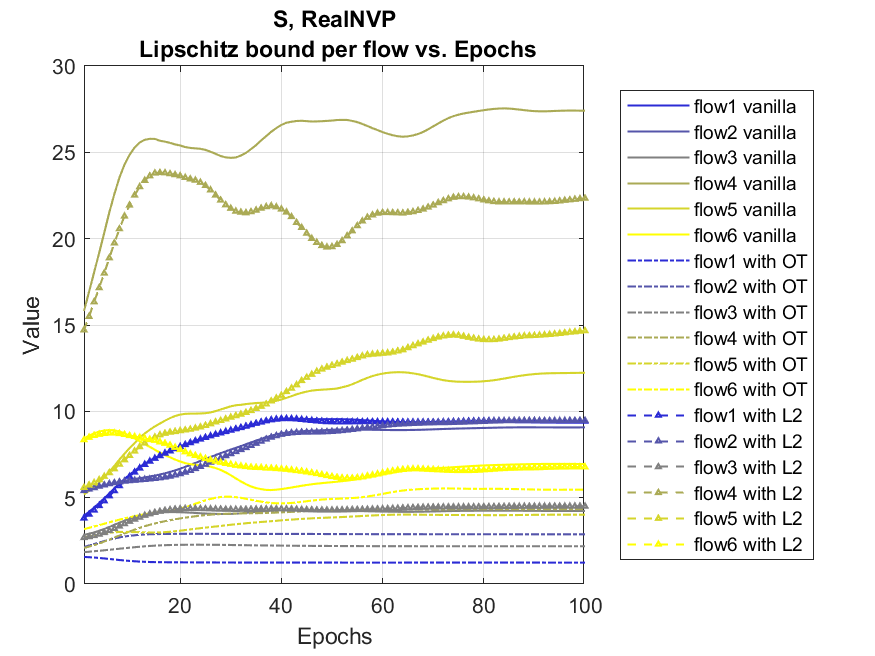}\\
\end{minipage}\hfill
\begin{minipage}{0.49\linewidth}
\centering
\includegraphics[width=.9\linewidth]{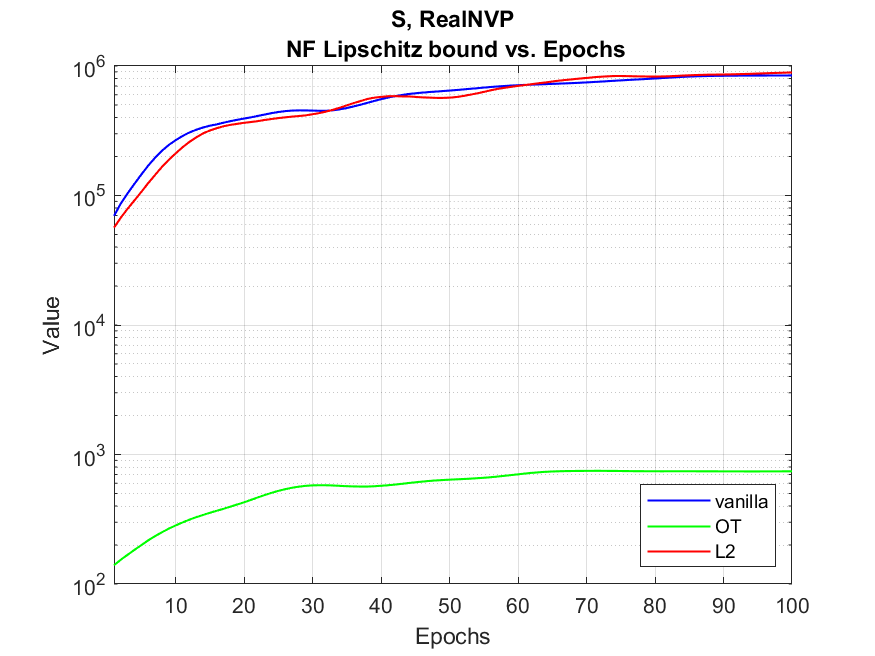}\\
\end{minipage}\hfill
\caption{Top row: output of each intermediate RealNVP flow between a single Gaussian and the S-shape density, trained without using the transport cost. Middle row: same but trained with using the transport cost. Bottom left: the Lipschitz bound for each flow over the training epochs. Bottom right: the Lipschitz bound for the entire flow. The weights of the transport regularization are chosen so that the negative log-likelihood is not severely obstructed.}
\label{RNVP_S}
\end{figure}

In this experiment, we compare the intermediate flows trained on 2D synthetic datasets, with and without the transport cost. 
%
%
%
%
The introduction of the transport cost encourages NFs to evolve the data to their desired destinations with the least amount of distortion. As a result, the intermediate flows are more sensible and visually appealing, as illustrated in Figure~\ref{RNVP_S} (see more examples in the appendix; Figures~\ref{RNVP_Swiss}, \ref{syn_fig_RealNVP}, \ref{syn_fig_NSF}). By comparing the NFs trained with and without the transport cost, we can see that in the absence of the regularization, existing NF models cannot be expected to learn trajectories that are approximately optimal under the $L_2$ cost. Instead, the learned trajectories are artifacts of the flow parameterization and they are sensitive to initialization. 

Visually, NFs trained with the transport cost incur significantly less distortion, so it is not surprising that the Lipschitz bound for each flow and that for the whole model are much smaller. In Figure~\ref{RNVP_S}, we compare the Lipschitz bounds between the standard NF, the NF trained with $l_2$ weight decay, and the NF trained with trajectory regularization. It is evident that while weight decay has a small impact on the Lipschitz bound, the trajectory regularization yields more noticeable improvement.

\subsubsection{Tabular Datasets}\label{sect_tabular_datasets}

\begin{table}[h]
\caption{The log-likelihoods obtained by using different models. Bolded results highlight the improvements of the OT-based models over their counterparts (e.g., RNVP + OT versus RNVP). The other NF models are included for reference. The numbers after the $\pm$ symbol are twice the standard errors.}
\label{Tabular_table}
\begin{center}
\begin{small}
\begin{sc}
\resizebox{\textwidth}{!}{\begin{tabular}{lccccccccccc}
\toprule
Model & POWER & GAS & HEPMASS & MINIBOONE & BSDS300\\
\midrule
 MAF (10)  & 0.24 $\pm$ 0.01
 & 10.08 $\pm$ 0.02
 & -17.73 $\pm$ 0.02
 & -12.24 $\pm$ 0.45
 & 154.93 $\pm$ 0.28\\
 MAF MoG  & 0.30 $\pm$ 0.01
 & 9.59 $\pm$ 0.02
 & -17.39 $\pm$ 0.02
 & -11.68 $\pm$ 0.44
 & 156.36 $\pm$ 0.28\\
 NAF  & 0.60 $\pm$ 0.02
 & 11.96 $\pm$ 0.33
 & -15.32 $\pm$ 0.23
 & -9.01 $\pm$ 0.01
 & 157.43 $\pm$ 0.30\\
 SOS  & 0.60 $\pm$ 0.01
 & 11.99 $\pm$ 0.41
 & -15.15 $\pm$ 0.10
 & -8.90 $\pm$ 0.11
 & 157.48 $\pm$ 0.41\\
 Quad. Spline  & 0.64 $\pm$ 0.01
 & 12.80 $\pm$ 0.02
 & -15.35 $\pm$ 0.02
 & -9.35 $\pm$ 0.44
 & 157.65 $\pm$ 0.28\\
 RNVP & 0.335 $\pm$ 0.013
 & 11.017 $\pm$ 0.022
 & -17.983 $\pm$ 0.021
 & -11.540 $\pm$ 0.469 & 153.398 $\pm$ 0.283\\
 NSF-AR & 0.650 $\pm$ 0.013
 & 13.001 $\pm$ 0.018
 & -14.145 $\pm$ 0.025
 & -9.662 $\pm$ 0.501
 & 157.546 $\pm$ 0.282\\
 \midrule
 RNVP + OT  & \textbf{0.374 $\pm$ 0.013} & \textbf{11.063 $\pm$ 0.022} & \textbf{-17.879 $\pm$ 0.0215}
 & \textbf{-11.167 $\pm$ 0.466}
 & \textbf{153.406 $\pm$ 0.283}\\
NSF-AR + OT  & \textbf{0.654 $\pm$ 0.012}
 &\textbf{ 13.018 $\pm$ 0.018}
 & \textbf{-13.989 $\pm$ 0.026}
 & \textbf{-9.395 $\pm$ 0.494}
 & \textbf{157.729 $\pm$ 0.280}
 \\
\bottomrule
\end{tabular}}
\end{sc}
\end{small}
\end{center}
\vspace{-.2cm}
\end{table}

\begin{figure}[h]
\centering
\begin{minipage}{0.33\linewidth}
\includegraphics[width=1\linewidth]{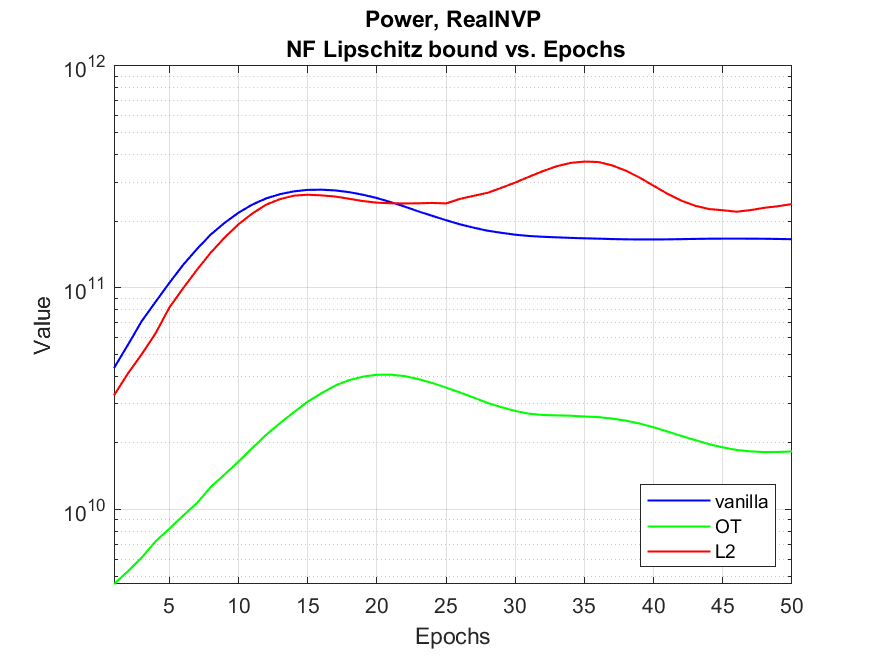}\\
\end{minipage}\hfill
\begin{minipage}{0.33\linewidth}
\includegraphics[width=1\linewidth]{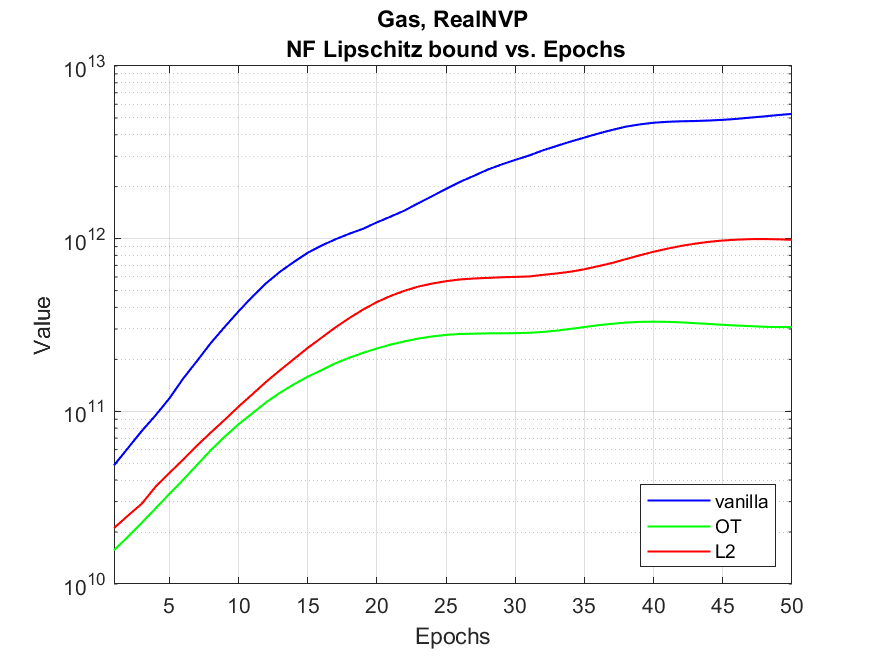}\\
\end{minipage}\hfill
\begin{minipage}{0.33\linewidth}
\includegraphics[width=1\linewidth]{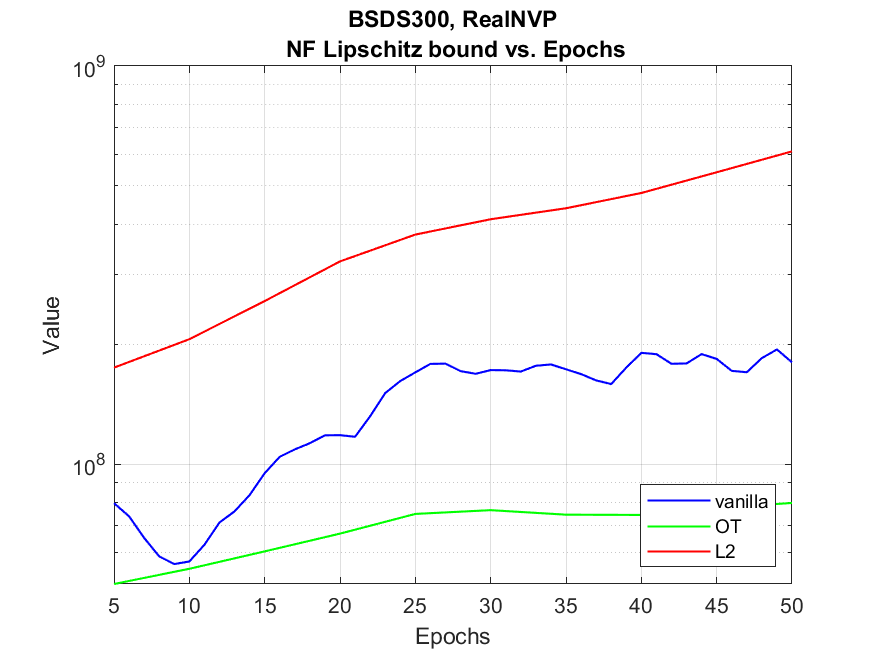}\\
\end{minipage}\hfill
%
\caption{The Lipschitz bound over trianing epochs, for different tabular datasets (left to right: Power, Gas, and BSDS300).}
\label{RNVP_tabular_Lip_bound}
\end{figure}

To demonstrate the benefit of training with transport costs, we replicate the density estimation experiments on popular benchmarks from the UCI repository and BSDS300~\cite{Tabular_dataset, BSDS300}, using the same setup as in~\cite{MAF}. Our approach requires no changes to the baseline architecture and hence any existing discrete NF is applicable. Here, we use two popular NF architectures, RealNVP and  NSF-AR, which are representatives of coupling flows and autoregressive flows, respectively~\cite{RealNVP, NSF}. The former uses simple (affine) coupling functions while the latter uses sophisticated (spline) coupling functions.

Table~\ref{Tabular_table} indicates that NF models trained with OT (trajectory regularization) are able to achieve better performance on the test set across all tabular datasets. The Lipschitz bounds for RealNVP trained on the datasets Power, Gas, and BSDS300 are visualized in Figure~\ref{RNVP_tabular_Lip_bound} (for other dadtasets, see the appendix). Clearly, while $l_2$ regularization has a mixed impact on the model's Lipschitz bound, the trajectory regularization reduces the bound effectively and improves the model's generalization performance.

We can interpret the discretized points on the OT trajectories as landmarks that should be visited by each flow in the interim. Our training setup is thus analogous to DSN~\cite{DSN}, which guides the learning of intermediate layers with intermediate classifiers. Consistent with their observation, the transport cost provides a natural mechanism for distributing losses throughout the network to improve convergence of the initial flows. 



\subsubsection{Image Datasets}\label{sect_image_datasets}

To further demonstrate the usefulness of our transport-efficient model, we add transport costs to Glow~\cite{glow}, a prominent NF architecture well-suited for generating high-quality images. Similar to tabular data, we observe improvements in density estimation via training Glow with the transport cost, in a number of datasets.

\begin{table}[ht]
\caption{The bits per dimension obtained by using different models ( Glow versus Glow with transport cost). The numbers after the $\pm$ symbol are twice the standard errors. }
\label{glow_table}
\begin{center}
\begin{small}
\begin{sc}
\resizebox{\textwidth}{!}{\begin{tabular}{lccccc}
\toprule
Model & MNIST & FMNIST & CIFAR-10 & SVHN & EUROSAT\\
\midrule
Glow & 1.182 $\pm$ 0.005 & 2.944 $\pm$ 0.017 & 3.489 $\pm$ 0.012 & 2.156 $\pm$ 0.005 & 2.670 $\pm$ 0.006
\\
 \midrule
 Glow + OT & \textbf{1.142 $\pm$ 0.005} & \textbf{2.901 $\pm$ 0.017} & \textbf{3.465 $\pm$ 0.012} & \textbf{2.146 $\pm$ 0.005} & \textbf{2.658 $\pm$ 0.006}\\
\bottomrule
\end{tabular}}
\end{sc}
\end{small}
\end{center}
\vspace{-.2cm}
\end{table}

 In Table~\ref{glow_table}, we select five popular image datasets to compare the density estimation results. The selected datasets include both greyscale and color images, with sizes varying between $28 \times 28$ and $64 \times 64$. The Glow implementation is down-scaled to fit on a single GPU; and the baseline results are similar to those reported by other works that follow this approach; e.g. \cite{OOD_gen}. In addition, we omit the 1x1 convolution transforms~\cite{glow} in the transport cost, as they play a role similar to permutations. Overall, we observe improvements in the esitmated density across all five image datasets, with similar regularization weights.






\section{Conclusion}

In this work, we unravel connections between the training of an NF and the solving of an MFG. Based on this insight, we explore two research  directions. First, we employ expressive flow transforms to accurately and efficiently solve high-dimensional MFG problems. We also approximation analysis for the proposed method. Compared to existing deep learning approaches, our approach encodes the continuity equation into the NF architecture and achieves a more accurate terminal matching. In addition, our implementation makes effective use of GPU parallelism to achieve fast computations in high dimensions. Second, we introduce transport costs into the standard NF training and demonstrate the effectiveness in controlling the network's Lipschitz bound. As a result, NFs trained with appropriate trajectory regularization are shown to achieve better generalization performance in density modeling across tabular and image datasets.

\bibliographystyle{plainnat}
\bibliography{references}

\newpage

\appendix

\DeclareRobustCommand{\vect}[1]{\bm{#1}}
\pdfstringdefDisableCommands{%
  \renewcommand{\vect}[1]{#1}%
}

\section{Derivations Regarding $F(\vect{x},t)$ and $\Phi_{t_{i}}^{t_{i+1}}(\vect{x})$}

In this section, we provide additional details for some statements made in~\eqref{sect_briding}, regarding the agent trajectories $F(\vx,t)$ and $\Phi_{t_{i}}^{t^{i+1}}(\vx)$. First, we will derive the ODE system defining $F(\vx,t)$ that is central to our MFG reformulation. Second, we will show how $F_k(\vx) \coloneqq F(\vx,t_k)$ can be decomposed into a series of flow maps $\Phi_{t_{i}}^{t_{i+1}}(\vx)$ that we subsequently parameterize with neural networks.

At the base level, we have in~\eqref{agent_traj} the trajectory $\vx: [0,T] \to \Omega$ for a single agent starting at $\vx_0\in \Omega$. From there, we can define the trajectories for all agents $F: \mathbb{R}^d \times [0,T] \to \mathbb{R}^d$ as $F(\vx,t) = \vx(t)$, where 
\begin{equation}
\begin{split}
    \der\vx(t) &= \vv(\vx(t),t)\der t, \quad \forall t\in [0,T]\\
    \vx(0) &= \vx.
\end{split}
\end{equation}
We see that $\partial_t F(\vx,t) = \der\vx(t) = \vv(\vx(t),t) = v(F(\vx,t),t), \forall t\in[0,T]$, and $F(\vx,0) = \vx(0) = \vx$. This gives the ODE system governing $F(\vx,t)$.

Next, we prove the additive property of $\Phi_{t_{i}}^{t_{j}}$.

\begin{prop}
For any $0 \leq t_a \leq t_b \leq t_c \leq T$, we have $\Phi_{t_{b}}^{t_{c}} \circ \Phi_{t_{a}}^{t_{b}} = \Phi_{t_{a}}^{t_{c}}$.
\end{prop}

\begin{proof}

By directly integrating the ODE governing $\vx(t)$,
\begin{align*}
    \Phi_{t_{a}}^{t_{b}}(\vx_0) = X(t_b) &= \vx_0 + \int_{t_{a}}^{t_{b}} \vv(\vx(s),s)\der s \\
    \implies \Phi_{t_{b}}^{t_{c}} ( \Phi_{t_{a}}^{t_{b}}(\vx_0)) &= (\vx_0 + \int_{t_{a}}^{t_{b}} \vv(\vx(s),s)\der s) + \int_{t_{b}}^{t_{c}} \vv(\vx(s),s)\der s\\
    &= \vx_0 + \int_{t_{a}}^{t_{c}} \vv(\vx(s),s)\der s\\
    &= \Phi_{t_{a}}^{t_{c}}(\vx_0).
\end{align*}

\end{proof}

Now, recalling the definition of $\Phi_{t_{i}}^{t_{j}}(\vx)$ in~\eqref{flow_map}, we see that $F_k(\vx) \coloneqq F(x,t_k) = \Phi_{0}^{t_{k}}(\vx) = \Phi_{t_{0}}^{t_{k}}(\vx)$, which by the additive property can be decomposed iteratively: $F_k = \Phi_{t_{0}}^{t_{k}} = \Phi_{t_{k-1}}^{t_{k}} \circ \Phi_{t_{0}}^{t_{k-1}} = \Phi_{t_{k-1}}^{t_{k}} \circ \Phi_{t_{k-2}}^{t_{k-1}} \circ \cdots \circ \Phi_{t_{0}}^{t_{1}}$.

\section{Derivation of the Alternative Formulation}
\label{append:alternative}
Consider reverse agent trajectories $\mathcal{G}: \mathbb{R}^d \times \mathbb{R} \to \mathbb{R}^d$ satisfying:
\begin{equation} \label{NF_traj_rev_3}
\begin{split}
    \partial_t \mathcal{G}(x,1-t) &= v(\mathcal{G}(x,1-t), t), \quad \vx\in \Omega, t\in[0,T] \\
    \mathcal{G}(x,0) &= x, \quad \vx\in \Omega\\
    \mathcal{G}(\cdot, t)_* P_1 &= P(\cdot, 1-t), \forall t\in [0,1].
\end{split}
\end{equation}
Applying the change of measure, the transport cost becomes:
\begin{align*}
    \int_0^1 \int_\Omega \langle \vv(\vx,t), \vv(\vx,t) \rangle p(\vx,t) \der\vx \der t &= \int_0^1 \int_\Omega \langle v(\mathcal{G}(x,1-t),t), v(\mathcal{G}(x,1-t),t) \rangle p_1(\vx) \der\vx \der t\\
    &= \int_0^1 \int_\Omega \langle \partial_t \mathcal{G}(x,1-t), \partial \mathcal{G}(x,1-t)\rangle p_1(\vx)\der x \der t\\
    &=  \mathbb{E}_{x\sim P_1} \left[\int_0^1\| \partial_t \mathcal{G}(x,1-t)\|^2_2\right] \der t.
\end{align*}
Defining a regular grid in time: $ t_k \coloneqq k \cdot \Delta t, \forall k = 0,1,...K$, $\Delta t \coloneqq \frac{1}{K}$, we obtain the approximation:
\begin{align*}
   \partial_t \mathcal{G}(x,1-t_k) \approx \frac{G_{K-k}(\vx) - G_{K-k-1}(\vx)}{\Delta t}, \forall k=0,1, .., K-1.
\end{align*}
Therefore, the discretized transport cost is:
\begin{align*}
    \mathbb{E}_{x\sim P_1} \left[\int_0^1\| \partial_t \mathcal{G}(x,1-t)\|^2_2\right] \der t &= \mathbb{E}_{x\sim P_1} \left[\sum_{k=0}^{K-1}\left\| \frac{G_{K-k}(\vx) - G_{K-k-1}(\vx)}{\Delta t} \right\|^2_2 \Delta t\right] \\
    &= K \cdot \mathbb{E}_{x\sim P_1} \left[\sum_{k=0}^{K-1}\| G_{k+1}(\vx) - G_{k}(\vx) \|^2_2\right].
\end{align*}

\section{Proofs}\label{appendix_proofs}
In this section, we provide the proofs of the statements mentioned in Section~\ref{sec:theory}. 
\subsection{Theorem~\ref{thm:disc_error}}\label{disc_error}

The proof essentially follows the error analysis used in finite difference.  Define a regular grid in time with grid points $ t_i \coloneqq i\cdot \Delta t, \forall i = 0,1,...K$ and grid spacing $\Delta t \coloneqq \frac{1}{K}$. Let $\hat{L}(t) = \int_{\mathbb{R}^d} \|\partial_t F(\vz,t)\|_2^2 p_0(\vz) \der \vz \coloneqq \|\partial_t F(\cdot,t)\|^2_{P_0}$. We start with the continuous objective
\begin{equation}\label{MFG_disc_err}
\begin{split}
    \int_0^1 \int_{\mathbb{R}^d} \|\partial_t F(\vz,t)\|_2^2 p_0(\vz) \det \vz\der t &\coloneqq \int_0^1\|\partial_t F(\cdot,t)\|_{P_0}^2  \der t\\
    &= \sum_{i=0}^{K-1} \int_{t_i}^{t_{i+1}} \hat{L}(t) \der t\\
    &= \sum_{i=0}^{K-1} \int_{t_i}^{t_{i+1}} \left[\hat{L}(t_i) + (t-t_i)\partial_t \hat{L}(t_i) + \frac{1}{2}(t-t_i)^2\partial_{tt}\hat{L}(t_i) + \cdots\right] \der t\\
    &= \sum_{i=0}^{K-1} \hat{L}(t_i)\Delta t + \partial_t \hat{L}(t_i)\Delta t^2  + \frac{1}{6}\partial_{tt}\hat{L}(t_i)\Delta t^3 + O(\Delta t^4).
\end{split}
\end{equation}
To obtain $\hat{L}(t_i)$, we use the forward difference:
\begin{align*}
    \partial_t F(x,t_i) &= \frac{F_{i+1}(\vx)-F_i(\vx)}{\Delta t} - \frac{1}{2} \Delta t\partial_{tt}F_i(\vx) + O(\Delta t^2)\\
    \implies \hat{L}(t_i) &= \left\|\frac{F_{i+1}-F_i}{\Delta t} - \frac{1}{2} \Delta t\partial_{tt}F_i + O(\Delta t^2)\right\|^2_{P_0}\\
    &= \left\|\frac{F_{i+1} - F_i}{\Delta t}\right\|^2_{P_0} - \left\langle \frac{F_{i+1} - F_i}{\Delta t}, \partial_{tt} F_i \right\rangle_{P_0} \Delta t + O(\Delta t^2),
\end{align*}
where $F_i(\vx) \coloneqq F(x,t_i)$. Therefore, if we approximate $\hat{L}(t_i) \approx \|\frac{F_{i+1} - F_i}{\Delta t}\|^2_{P_0}$, the error of the integral is:
\begin{equation}\label{MFG_disc_err_2}
\begin{split}
    &\sum_{i=0}^{K-1} \left[\left\|\frac{F_{i+1} - F_i}{\Delta t}\right\|^2_{P_0} - \left\|\frac{F_{i+1} - F_i}{\Delta t}\right\|^2_{P_0} + \left\langle \frac{F_{i+1} - F_i}{\Delta t}, \partial_{tt} F_i \right\rangle_{P_0} \Delta t + O(\Delta t^2)\right]\Delta t
     + \partial_t \hat{L}(t_i)\Delta t^2 + O(\Delta t^3)\\
    &= \sum_{i=0}^{K-1} \left[\left\langle \frac{F_{i+1} - F_i}{\Delta t}, \partial_{tt} F_i \right\rangle_{P_0} + \partial_t \hat{L}(t_i) \right] \Delta t^2 + O(\Delta t^3)\\
    &= \sum_{i=0}^{K-1} [\langle \partial_t F_i + O(\Delta t), \partial_{tt} F_i \rangle_{P_0} + \partial_t \hat{L}(t_i) ] \Delta t^2 + O(\Delta t^3)\\
    &= \sum_{i=0}^{K-1} [\langle \partial_t F_i, \partial_{tt} F_i \rangle_{P_0} + \partial_t \hat{L}(t_i) ] \Delta t^2 + O(\Delta t^3)\\
    &= O(\Delta t),
\end{split}
\end{equation}
where the last equality is a result of summing over $K$ grid points.

\subsection{Lemma~\ref{cont_MFG_soln}}

First, note that the MFG problem~\eqref{var_MFG_traj} admits the same optimizer and objective value if we add $F(\vz,1) = F_1(\vz)$ as an additional constraint:
\begin{equation} \label{MFG_2}
\begin{split}
    &\inf_{F} \quad  \lambda_L\int_0^1 \int_{\mathbb{R}^n} p_0(\vz) \|\partial_t F(\vz,t)\|_2^2 \der z \der t  + \mathcal{M}(F(\cdot, 1)_*p_0)\\
    &s.t. \quad  F(\vz,0) = x, F(\vz,1) = F_1(\vz).
\end{split}
\end{equation}
Note that $\mathcal{M}(F(\cdot, 1)_*p_0)$ is a constant with respect to $F$. Therefore, the above problem has the same solution as:
\begin{equation} \label{MFG_3}
\begin{split}
    &\inf_{F} \quad  \lambda_L\int_0^1 \int_{\mathbb{R}^n} p_0(\vz) \|\partial_t F(\vz,t)\|_2^2 \der z \der t\\
    &s.t. \quad  F(\vz,0) = x, F(\vz,1) = F_1(\vz),
\end{split}
\end{equation}
which in turn is equivalent to the following OT problem:
\begin{equation} \label{OT_1}
\begin{split}
    &\inf_{p, \vv} \quad  \lambda_L\int_0^1 \int_{\mathbb{R}^n} p(\vx,t) \|\vv(\vx,t)\|_2^2 \der\vx \der t\\
    &s.t. \quad  p(\vx,0) = p_0(\vx), p(\vx,1) = F_{1*}p_0(\vx),
\end{split}
\end{equation}
where $ F_{1*}p_0(\vx) = p_0(F^{-1}_1(\vx))|\det \nabla F^{-1}_1(\vx)|$. By the OT theory, we know that there exists a unique Monge map $T(\vx)$ such that $p(\vx,1) = T_*p_0(\vx)$~\cite{OT_book}. Therefore, $T = F_1$ and the optimizer $F^*$ of~\eqref{MFG_3} satisfies $F^*(\vz,t) = (1-t)z + tT(\vz) = (1-t)\vz + tF_1(\vz)$. As noted before, $F^*$ minimizes the continuous MFG problem~\eqref{var_MFG_traj} as well.

\subsection{Lemma~\ref{disc_MFG_soln}}





Applying a similar logic to the proof of Lemma~\ref{cont_MFG_soln}, we know that the optimizer of the discretized MFG problem~\eqref{MFG_K} coincides with that of the following problem
\begin{align*}
    &\inf_{\{F_i\}_{i=0}^K} \quad \lambda_L \sum_{i=0}^{K-1} \mathbb{E}_{z\sim P_0} [\|F_{i+1}(\vz) - F_i(\vz)\|^2_2] \coloneqq  \lambda_L \sum_{i=0}^{K-1} \|F_{i+1} - F_i\|^2_{P_0}\\
    &s.t. \quad  F_0(\vz) = \vz, F_K(\vz) = F_1(\vz),
\end{align*}
where $F_1$ is an optimizer of the discretized MFG problem for $F_K$, and $ \|f\|_{P_0}^2 \coloneqq \int_{\mathbb{R}^d} f(\vx)^2 p_0(\vx)\der x$. Applying the KKT conditions on the above problem, we obtain: 
\begin{align*}
    &\partial_{F_i} \lambda_L \sum_{i=0}^{K-1} \|F_{i+1} - F_i\|^2_{P_0} = 0\\
    \implies &F^*_i = \frac{1}{2} (F^*_{i+1} + F^*_{i-1}),\quad \forall i = 1,2,..., K-1.
\end{align*}
Along with $F_0^* = Id, F_K^* = F_1$, we see that the solution to this recurrence relation is $F_i^* = z(1-\frac{i}{K}) + F_1\frac{i}{K}, i=0, 1, ..., K$.


\subsection{Theorem~\ref{cont_disc_MFG_soln}}

From Lemma~\ref{cont_MFG_soln}, we know the optimizer of the continuous MFG problem~\eqref{var_MFG_traj} has the form
\begin{align*}
    F^*(\vz,t) &= z(1-t) + T(\vx) t.
\end{align*}
We can substitute this into the continuous MFG problem and optimize over $T$ instead of $F$, yielding
\begin{align*}
&\inf_{T} \quad  \lambda_L\int_0^1 \int_{\mathbb{R}^d} p_0(\vz) \|\partial_t [z(1-t) + T(\vz)t)]\|_2^2 \der z \der t  + \mathcal{M}(T_*p_0)\\
= &\inf_{T} \quad   \lambda_L\int_0^1 \int_{\mathbb{R}^d} p_0(\vz) \|T(\vz) - z\|_2^2 \der z \der t  + \mathcal{M}(T_*p_0) \\
= &\inf_{T} \quad   \lambda_L \int_{\mathbb{R}^d} p_0(\vz) \|T(\vz) - z\|_2^2 \der z  + \mathcal{M}(T_*p_0).
\end{align*}

Second, we know from~\eqref{disc_MFG_soln} the optimizer $\{F_i^*\}_{i=0}^K$ of the discretized MFG problem with $K$ grid points has the form $F_i^* = z(1-\frac{i}{K}) + T_K\frac{i}{K}, i=0, 1, ..., K$. As a result, we can rewrite the discretized MFG problem as
\begin{align*}
&\inf_{T_K} \quad \lambda_L K\cdot \mathbb{E}_{z\sim P_0} \left[ \sum_{k=0}^{K-1}\left\|z\left(1-\frac{k+1}{K}\right) + T_K(\vz)\frac{k+1}{K} - z\left(1-\frac{k}{K}\right) - T_K(\vz)\frac{k}{K}\right\|_2^2\right] +  \mathcal{M}(T_{K_*}p_0)\\
= &\inf_{T_K} \quad \lambda_L K\cdot \sum_{k=0}^{K-1}\int_{\mathbb{R}^d} p_0(\vz) \left\|-\frac{1}{K} z + \frac{1}{K} T_K(\vz)\right\|^2_2 \der z + \mathcal{M}(T_{K_*}p_0)\\
= &\inf_{T_K} \quad \lambda_L \cdot \int_{\mathbb{R}^d} p_0(\vz) \| T_K(\vz) - z\|^2_2 \der z + \mathcal{M}(T_{K_*}p_0),
\end{align*}
which is identical to the problem over $T$ for the continuous MFG. Therefore, their objective values agree; furthermore, the minimizers of the two transformed problems are also the same: $T = T_K$. As a result, the mapping at the terminal time agrees between the continuous and the discretized problem: $F^*(\vz, 1) = F_K^*(\vz)$. It then follows that the linear interpolation between the initial and terminal mappings in the discretized problem, $z(1-t) + F_K^*(\vz)t$, is precisely the unique optimizer for the continuous problem.


\subsection{Theorem~\ref{thm_conv_measure}}
The proof can be obtained by combining Theorem~\ref{thm_optimizer_conv} with the following lemma. 

\begin{lem}\label{lem_conv_measure}
Let $P_0$ be absolutely continuous with a bounded density $p_0$. Suppose for any compact $A\subset \mathbb{R^d}$, there is a sequence of mappings $\{F_A^n\}_{n=1}^\infty$ such that  $\|F^n_A - F\|_{1,A} \to 0$. Then, there exists $\{F^n\}_{n=1}^\infty$ such that $F^n_{*} P_0 \xrightarrow{d} F_*P_0$, where $\xrightarrow{d}$ denotes convergence in distribution.
\end{lem}
\begin{proof}
Define 
\begin{align*}
    \mathbf{BL_1} = \{\eta: \eta \text{ bounded and Lipschitz with } Lip(\eta) + \|\eta\|_\infty \leq 1\}. 
\end{align*}
According to~\cite{analysis_book}, to show $ F^n_{*}P_0 \xrightarrow{d} F_*P_0$, it suffices to prove that:
\begin{align*}
    \sup_{\eta\in \mathbf{BL_1}} |\mathbb{E}_{F^n_{*}P_0}[\eta] - \mathbb{E}_{F_*{P_0}}[\eta]| < \frac{1}{n}, \quad \forall n \in \mathbb{N}.
\end{align*}
For each $n\in \mathbb{N}$, pick a compact $A_n \subset \mathbb{R}^d$ such that $P_0(A_n^c) = P_0(\mathbb{R}^d) - P_0(A_n) = 1 - P_0(A_n) < \frac{1}{2n}$. By assumption, we can find a sequence $\{F^k_n\}_{k=1}^\infty$ satisfying $\displaystyle \|F^k_n - F\|_{1,A_n} < \frac{1}{k\|p_0\|_{\infty,A_n}}$. Construct the sequence $\{F^n\}_{n=1}^\infty \coloneqq \{F^n_n\}_{n=1}^\infty$, which satisfies $\displaystyle \|F^n - F\|_{1,A_n} < \frac{1}{n\|p_0\|_{\infty,A_n}}, \forall n\in \mathbb{N}$.

For any $\eta\in \mathbf{BL_1}$, we have:
\begin{align*}
    |\mathbb{E}_{F^n_{*}P_0}[\eta] - \mathbb{E}_{F_*{P_0}}[\eta]| &= \left|\int_{\mathbb{R}^d} \eta \der F^n_{*}P_0 - \eta \der F_*P_0\right|\\
    &= \left|\int_{\mathbb{R}^d} \eta  \circ F^n \der P_0 - \eta  \circ F \der P_0\right|\\
    &= \left|\int_{A_n^c} \eta  \circ F^n \der P_0 + \int_{A_n} \eta  \circ F^n \der P_0 - \int_{A_n^c} \eta  \circ F \der P_0 - \int_{A_n} \eta  \circ F \der P_0\right|\\
    &\leq \left|\int_{A_n^c} \eta  \circ F^n \der P_0\right| + \left| \int_{A_n^c} \eta  \circ F \der P_0\right| +  \left|\int_{A_n} \eta  \circ F^n \der P_0 - \eta  \circ F \der P_0\right|\\
    &\leq  \int_{A_n^c} |\eta \circ F^n| \der P_0 + \int_{A_n^c} |\eta  \circ F| \der P_0 + \int_{A_n} |\eta (F^n(\vx)) - \eta (F(\vx))|\der P_0\\
    &\leq \|\eta \|_\infty P_0(A_n^c) + \|\eta \|_\infty P_0(A_n^c) + Lip(\eta ) \int_A \|F^n(\vx) - F(\vx)\|_1 \der P_0\\
    &\leq 2\|\eta \|_\infty P_0(A_n^c) + Lip(\eta ) \|p_0\|_{\infty,A_n} \int_{A_n} \|F^n(\vx) - F(\vx)\|_1 \der x\\
    &< \|\eta \|_\infty \frac{1}{n} + Lip(\eta ) \frac{1}{n}\\
    &= (\|\eta\|_\infty + Lip(\eta ))\frac{1}{n} \leq \frac{1}{n}.
\end{align*}

\end{proof}

\section{MFG Experiments}

We use NSF-CL with identical hyperparameters to solve both the OT and the crowd motion problems across different dimensions. The NF model contains 10 flows, each consisting of an alternating block with rational-quadratic splines as the coupling function, followed by a linear transformation. The conditioner has 256 hidden features, 8 bins, 2 transform blocks, ReLU activation, and a dropout probability of 0.25. 

In both the OT and the crowd motion experiments, we use Jeffery's divergence for terminal matching, which is the symmetrized KL divergence: $\mathcal{M}(F(\cdot, T)_*P_0) = D_{KL}(P_1 || F(\cdot, T)_*P_0) + D_{KL}(F(\cdot, T)_*P_0 || P_1)$. This terminal cost is different from the typical negative log-likelihood loss in NF training, but it is still tractable because $p_0, p_1$ are known. The choice of the terminal cost is not unique. For example, one can use either forward or  backward KL divergence. However, we find empirically that Jeffery's divergence yields the best results. 

In the Monte Carlo approximation of the expectation, we resample a batch of data from the initial and the terminal distributions to compute the various loss terms at each iteration. Although~\cite{ruthotto2020machine} remarked that this approach yields slightly worse performance than resampling every 20 iterations, we find it to be good enough for our method.

The computed dynamics for the OT and the crowd motion experiments are approximately identical with respect to the dimensionality, so we select only one $d\in \{2,10,50,100\}$ to show in the main text. For completeness, we include the remaining plots here.

\subsection{OT}

We show the density evolution and the sample trajectories for $d\in \{10,100\}$ here; for the case $d=50$, see Figure~\ref{OT_density_evo_50D} in the main text. Visually, the evolutions are approximately invariant with respect to the dimensionality.

\begin{figure}[H]
\centering
\begin{minipage}{0.16\linewidth}
\includegraphics[width=1\linewidth]{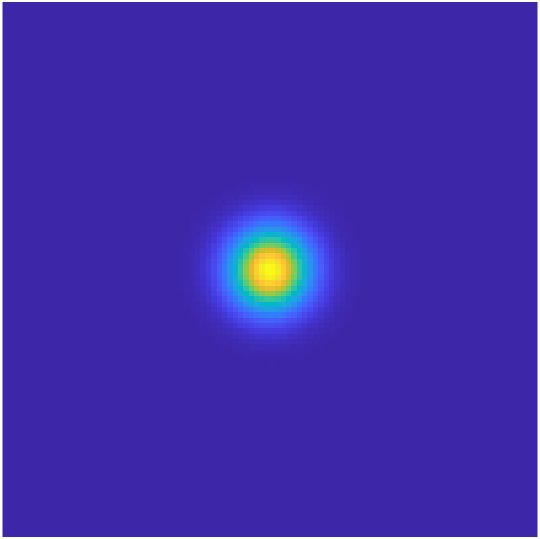}\\
\centering {$P_0$}
\end{minipage}\hfill
\begin{minipage}{0.16\linewidth}
\includegraphics[width=1\linewidth]{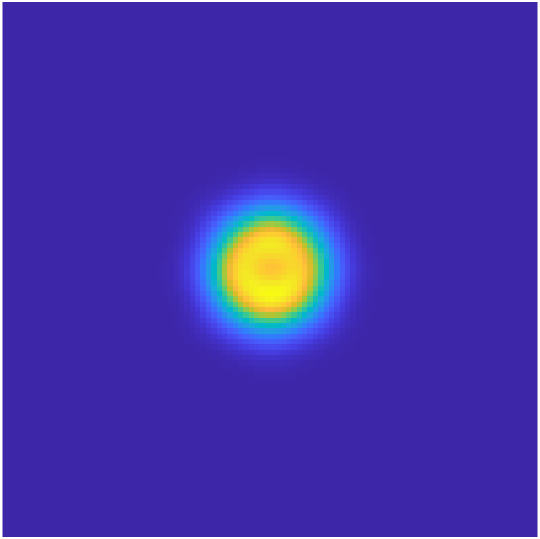}\\
\centering {$F_{1*}(P_0)$}
\end{minipage}\hfill
\begin{minipage}{0.16\linewidth}
\includegraphics[width=1\linewidth]{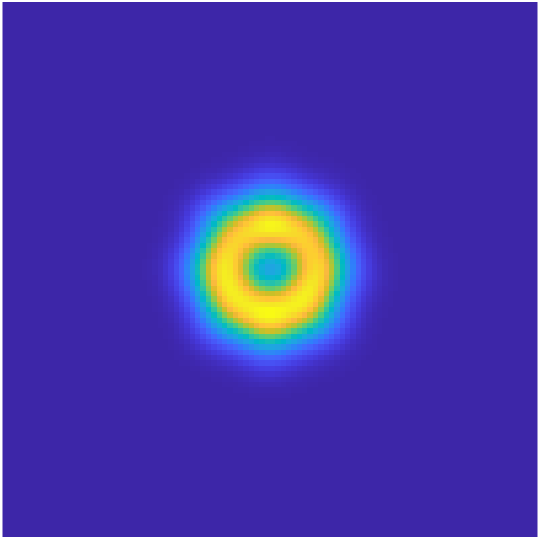}\\
\centering {$F_{2 *}(P_0)$}
\end{minipage}\hfill
\begin{minipage}{0.16\linewidth}
\includegraphics[width=1\linewidth]{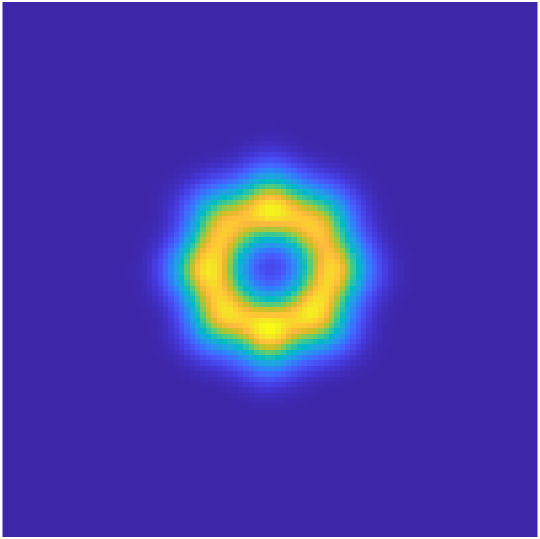}\\
\centering {$F_{3 *}(P_0)$}
\end{minipage}\hfill
\begin{minipage}{0.16\linewidth}
\includegraphics[width=1\linewidth]{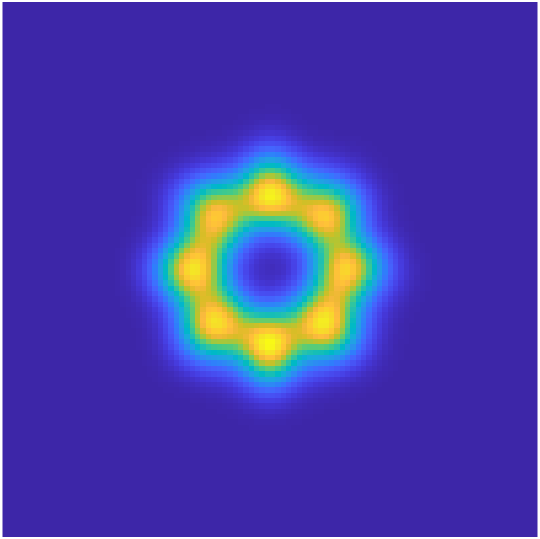}\\
\centering {$F_{4 *}(P_0)$}
\end{minipage}\hfill
\begin{minipage}{0.16\linewidth}
\includegraphics[width=1\linewidth]{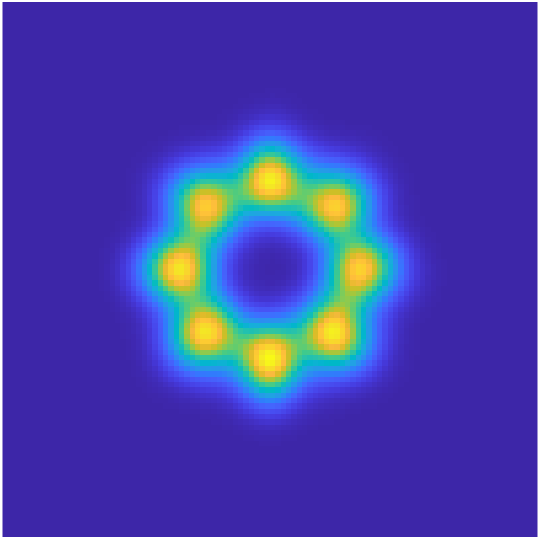}\\
\centering {$F_{5_*}(P_0)$}
\end{minipage}\hfill
\begin{minipage}{0.16\linewidth}
\includegraphics[width=1\linewidth]{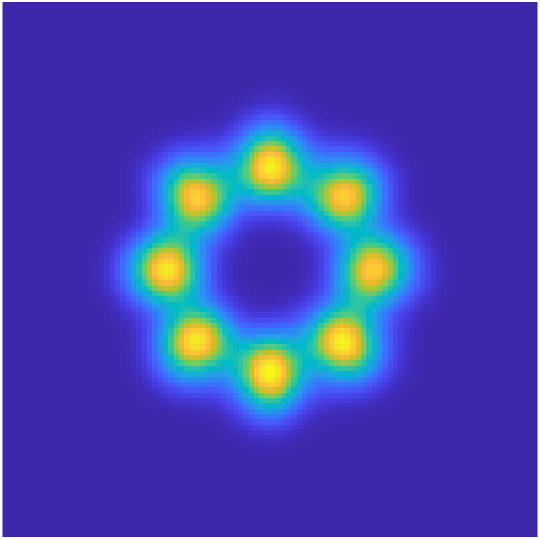}\\
\centering {$F_{6_*}(P_0)$}
\end{minipage}\hfill
\begin{minipage}{0.16\linewidth}
\includegraphics[width=1\linewidth]{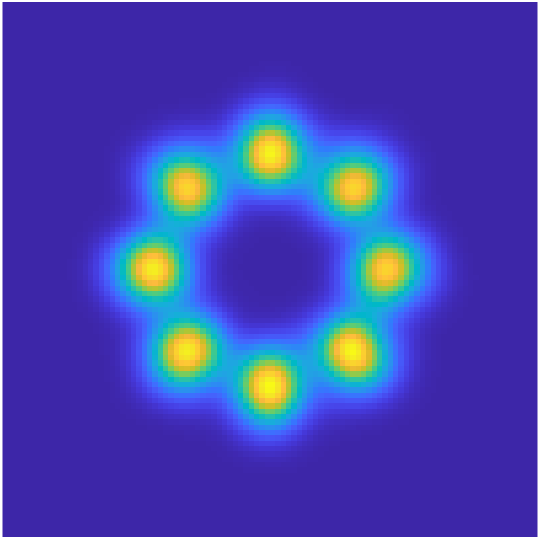}\\
\centering {$F_{7_*}(P_0)$}
\end{minipage}\hfill
\begin{minipage}{0.16\linewidth}
\includegraphics[width=1\linewidth]{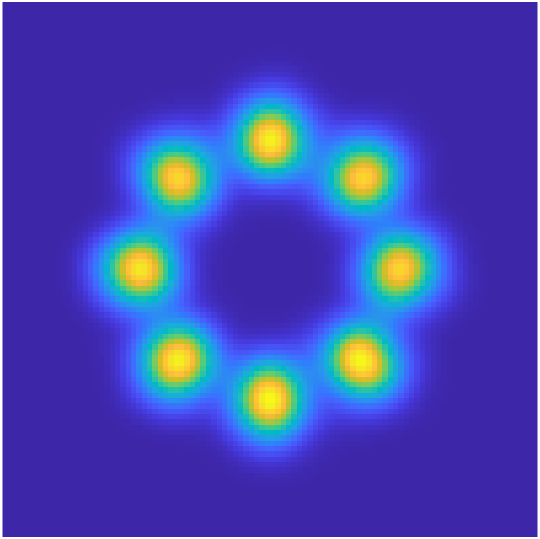}\\
\centering {$F_{8_*}(P_0)$}
\end{minipage}\hfill
\begin{minipage}{0.16\linewidth}
\includegraphics[width=1\linewidth]{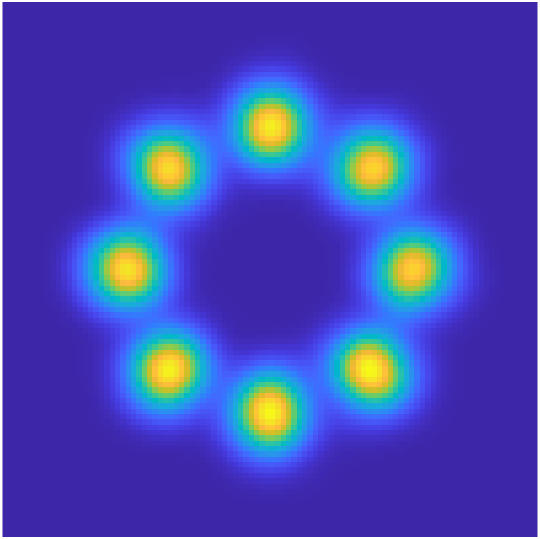}\\
\centering {$F_{9_*}(P_0)$}
\end{minipage}\hfill
\begin{minipage}{0.16\linewidth}
\includegraphics[width=1\linewidth]{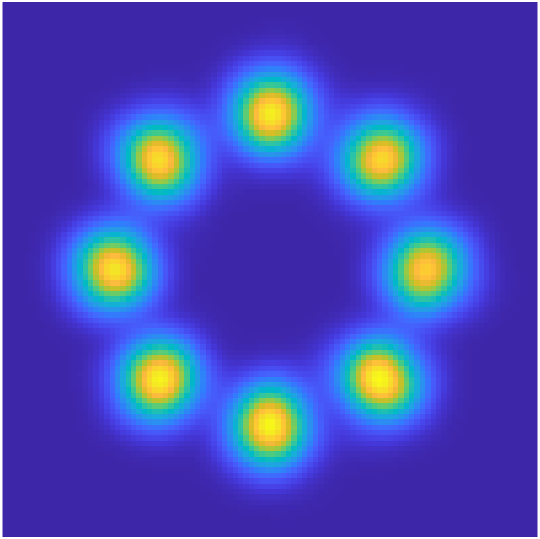}\\
\centering {$F_{10*}(P_0)$}
\end{minipage}\hfill
\begin{minipage}{0.16\linewidth}
\includegraphics[width=1\linewidth]{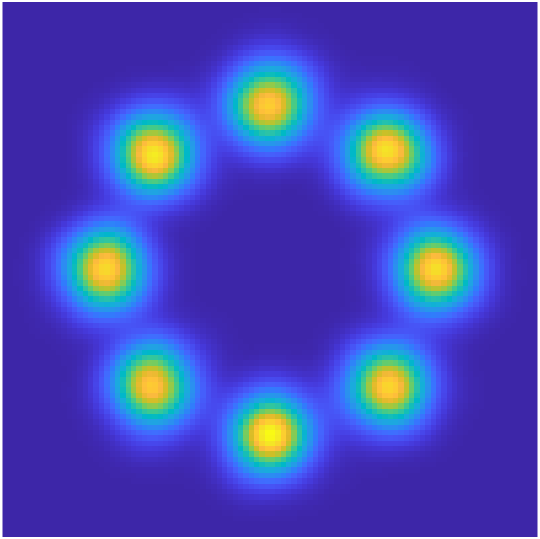}\\
\centering {$P_1$}
\end{minipage}\hfill
\caption{Evolution of the density for the OT problem~\eqref{OT} in 10D.}
\label{fig:OT_trajectory_10D}
\end{figure}

\begin{figure}[H]
\centering
\begin{minipage}{0.16\linewidth}
\includegraphics[width=1\linewidth]{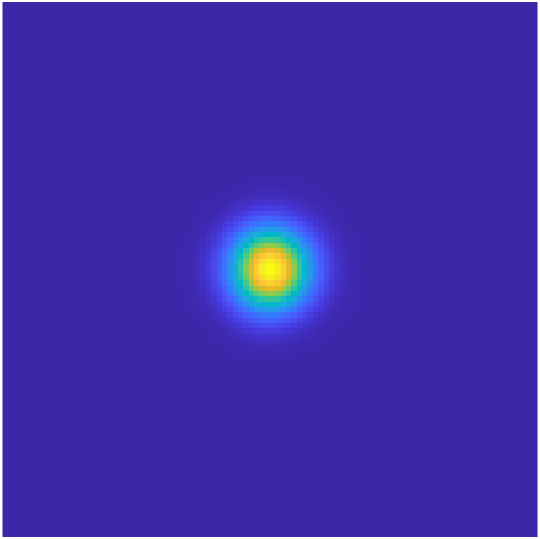}\\
\centering {$P_0$}
\end{minipage}\hfill
\begin{minipage}{0.16\linewidth}
\includegraphics[width=1\linewidth]{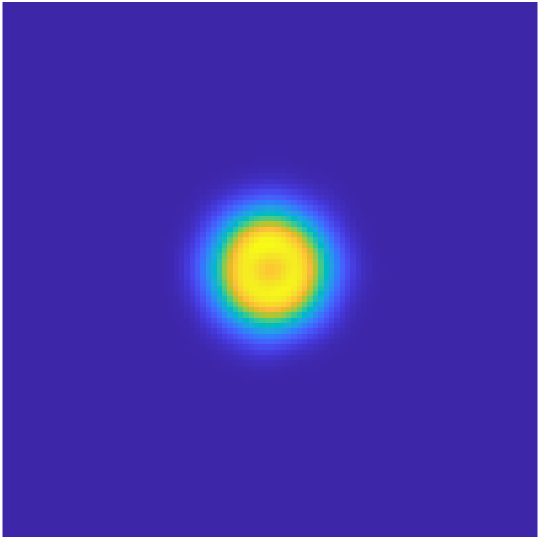}\\
\centering {$F_{1*}(P_0)$}
\end{minipage}\hfill
\begin{minipage}{0.16\linewidth}
\includegraphics[width=1\linewidth]{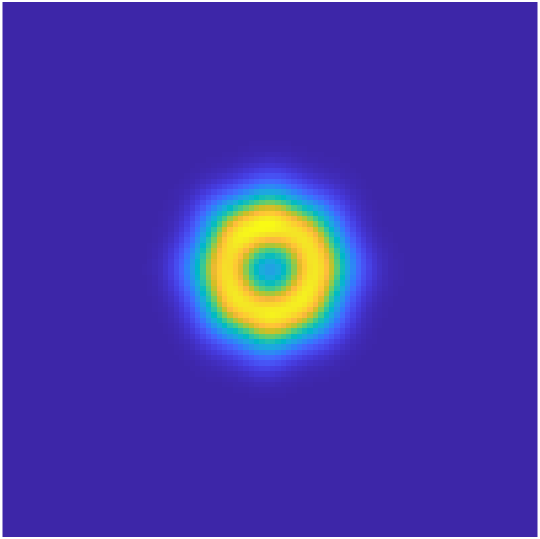}\\
\centering {$F_{2*}(P_0)$}
\end{minipage}\hfill
\begin{minipage}{0.16\linewidth}
\includegraphics[width=1\linewidth]{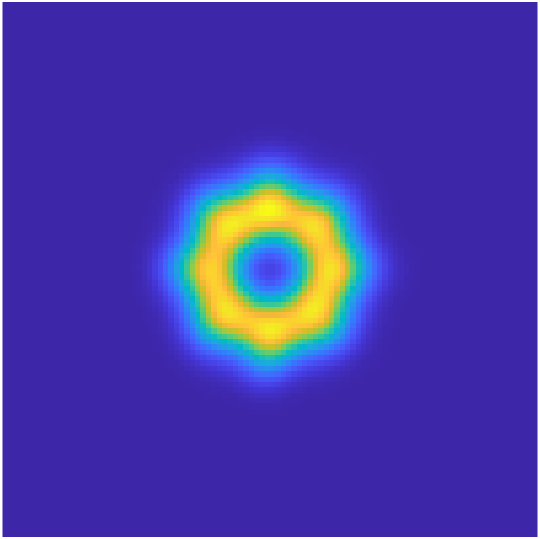}\\
\centering {$F_{3*}(P_0)$}
\end{minipage}\hfill
\begin{minipage}{0.16\linewidth}
\includegraphics[width=1\linewidth]{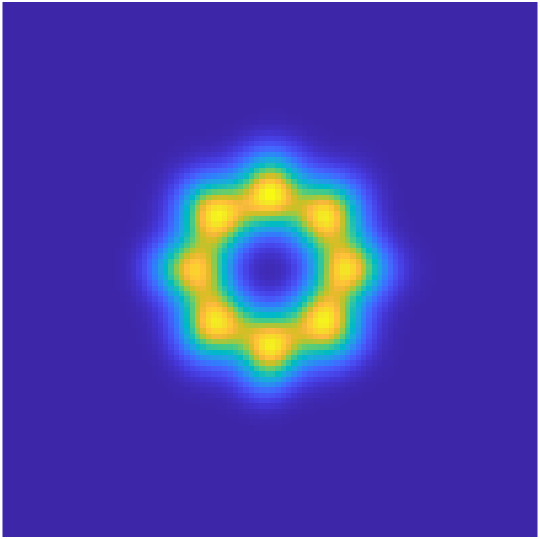}\\
\centering {$F_{4*}(P_0)$}
\end{minipage}\hfill
\begin{minipage}{0.16\linewidth}
\includegraphics[width=1\linewidth]{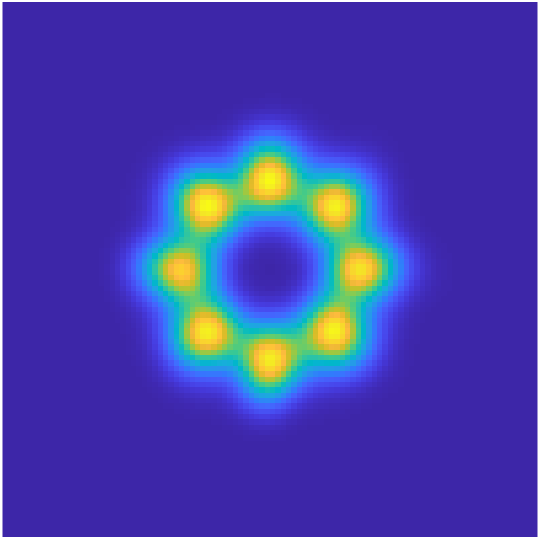}\\
\centering {$F_{5*}(P_0)$}
\end{minipage}\hfill
\begin{minipage}{0.16\linewidth}
\includegraphics[width=1\linewidth]{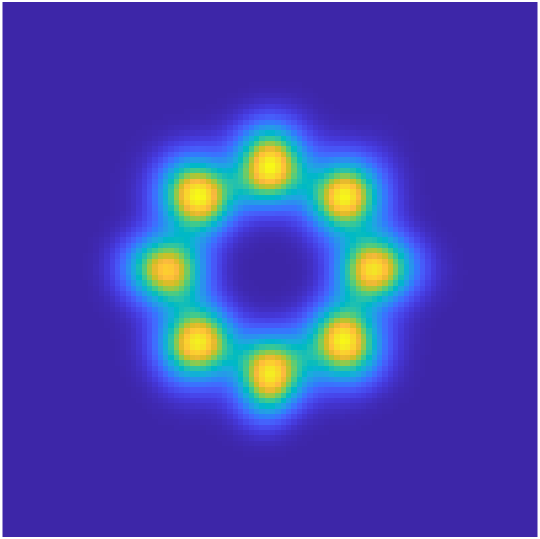}\\
\centering {$F_{6*}(P_0)$}
\end{minipage}\hfill
\begin{minipage}{0.16\linewidth}
\includegraphics[width=1\linewidth]{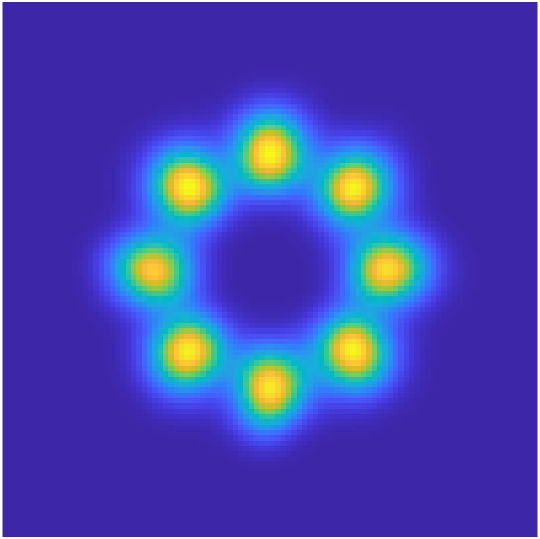}\\
\centering {$F_{7 *}(P_0)$}
\end{minipage}\hfill
\begin{minipage}{0.16\linewidth}
\includegraphics[width=1\linewidth]{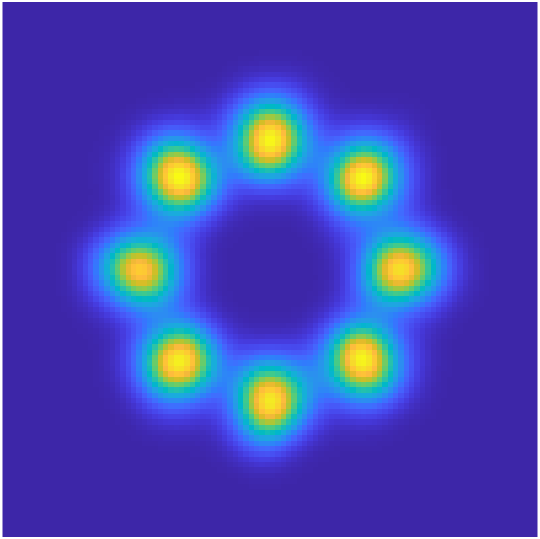}\\
\centering {$F_{8*}(P_0)$}
\end{minipage}\hfill
\begin{minipage}{0.16\linewidth}
\includegraphics[width=1\linewidth]{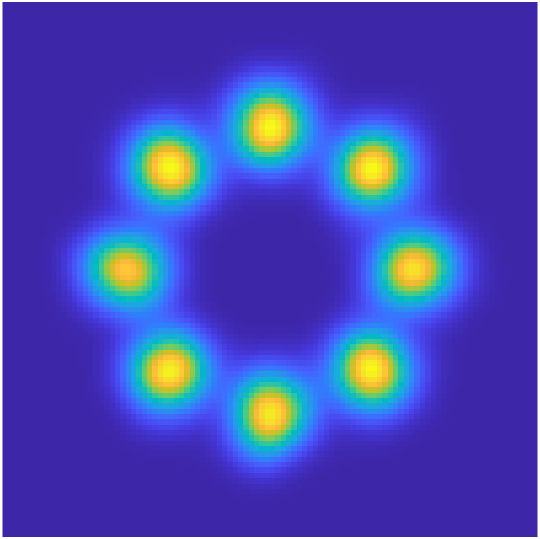}\\
\centering {$F_{9*}(P_0)$}
\end{minipage}\hfill
\begin{minipage}{0.16\linewidth}
\includegraphics[width=1\linewidth]{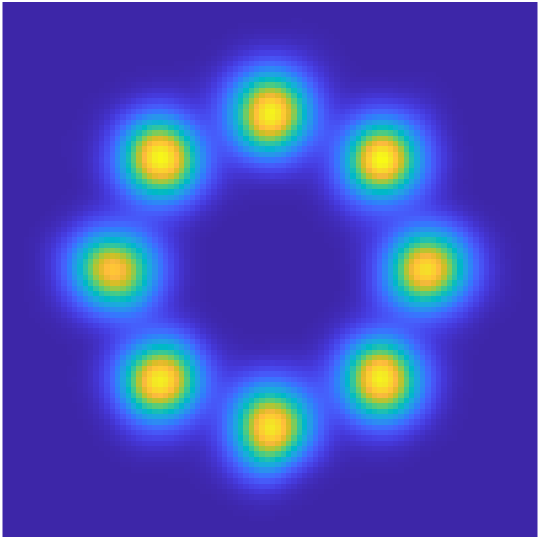}\\
\centering {$F_{10*}(P_0)$}
\end{minipage}\hfill
\begin{minipage}{0.16\linewidth}
\includegraphics[width=1\linewidth]{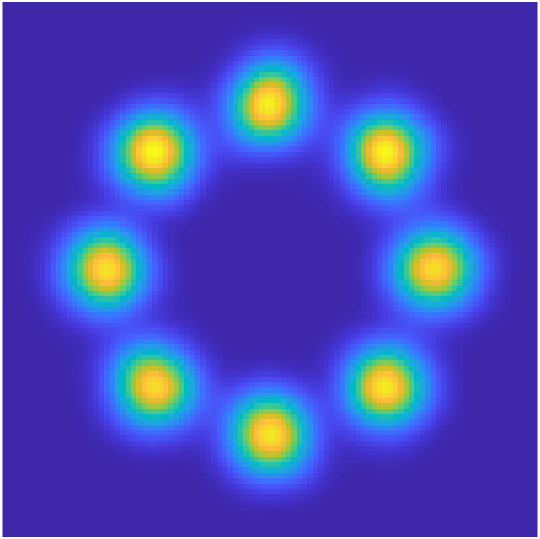}\\
\centering {$P_1$}
\end{minipage}\hfill

\caption{Evolution of the density for the OT problem~\eqref{OT} in 100D.}
\label{fig:OT_trajectory_100D}
\end{figure}

\subsection{Crowd Motion}

\subsubsection{Density Evolution in Different Dimensions}

We show the density evolution for $d\in \{2,10,50,100\}$ here, which further supports the robustness with respect to dimensionality. 

\begin{figure}[H]
\centering
\begin{minipage}{0.16\linewidth}
\includegraphics[width=1\linewidth]{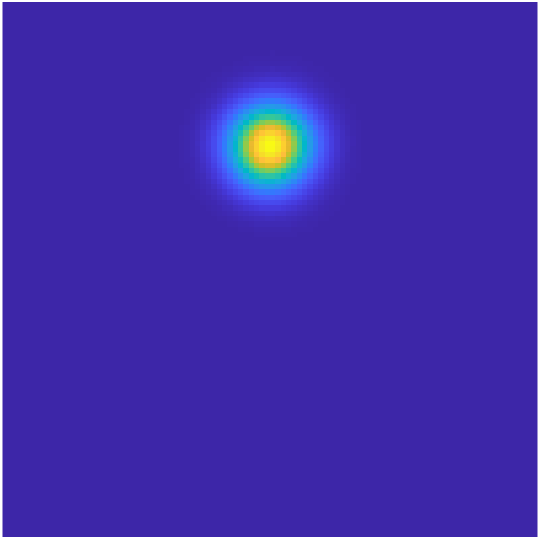}\\
\centering {$P_0$}
\end{minipage}\hfill
\begin{minipage}{0.16\linewidth}
\includegraphics[width=1\linewidth]{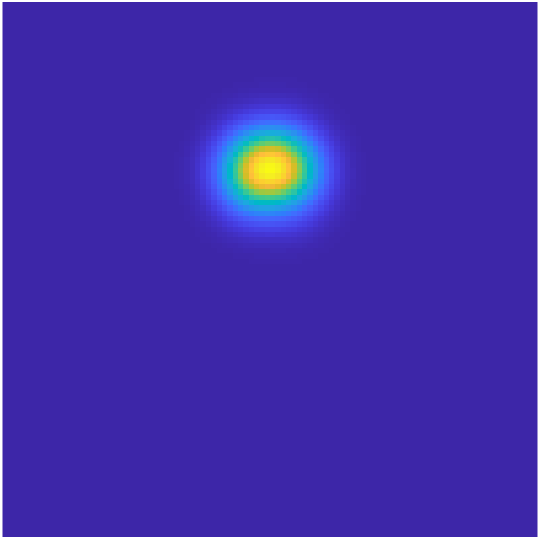}\\
\centering {$F_{1*}(P_0)$}
\end{minipage}\hfill
\begin{minipage}{0.16\linewidth}
\includegraphics[width=1\linewidth]{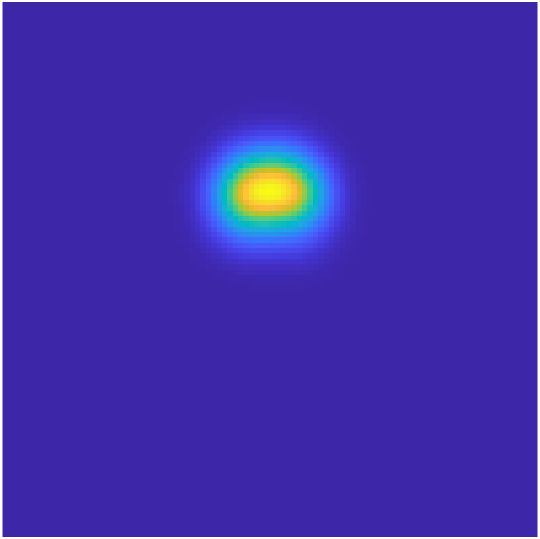}\\
\centering {$F_{2*}(P_0)$}
\end{minipage}\hfill
\begin{minipage}{0.16\linewidth}
\includegraphics[width=1\linewidth]{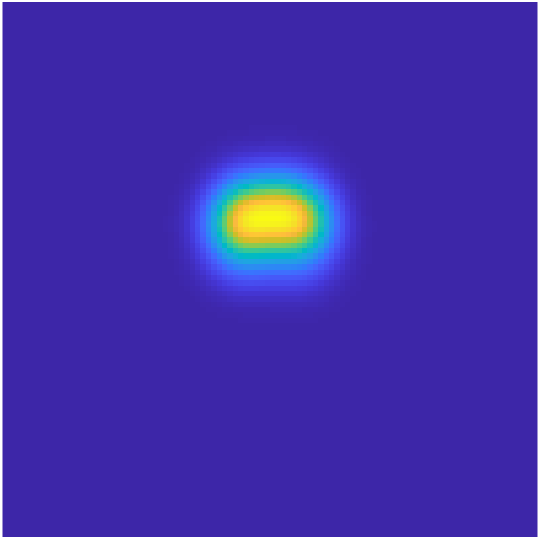}\\
\centering {$F_{3*}(P_0)$}
\end{minipage}\hfill
\begin{minipage}{0.16\linewidth}
\includegraphics[width=1\linewidth]{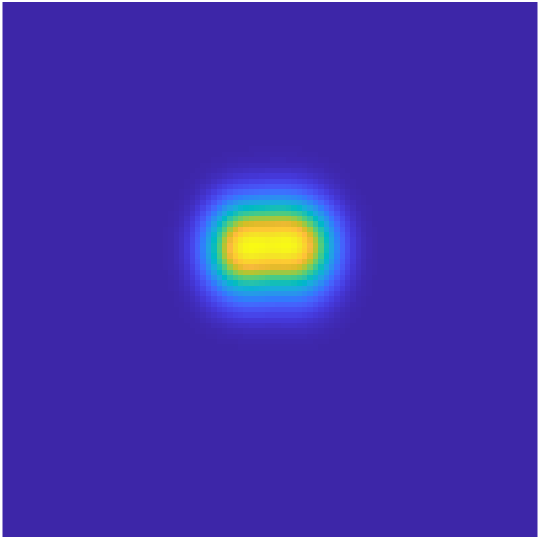}\\
\centering {$F_{4*}(P_0)$}
\end{minipage}\hfill
\begin{minipage}{0.16\linewidth}
\includegraphics[width=1\linewidth]{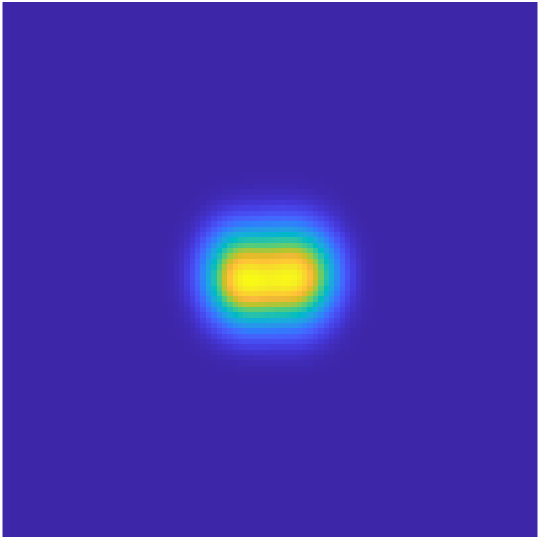}\\
\centering {$F_{5*}(P_0)$}
\end{minipage}\hfill
\begin{minipage}{0.16\linewidth}
\includegraphics[width=1\linewidth]{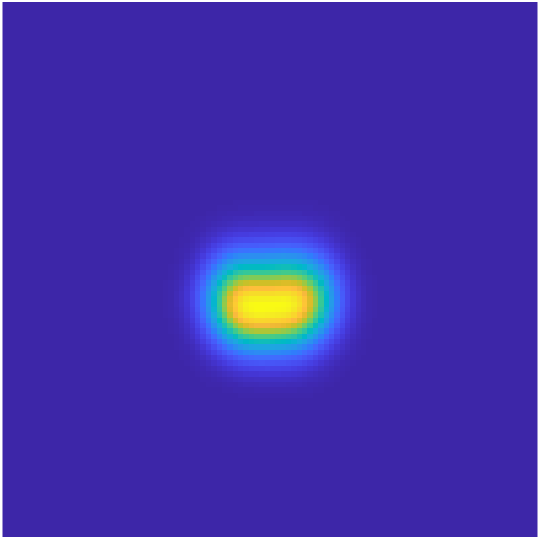}\\
\centering {$F_{6*}(P_0)$}
\end{minipage}\hfill
\begin{minipage}{0.16\linewidth}
\includegraphics[width=1\linewidth]{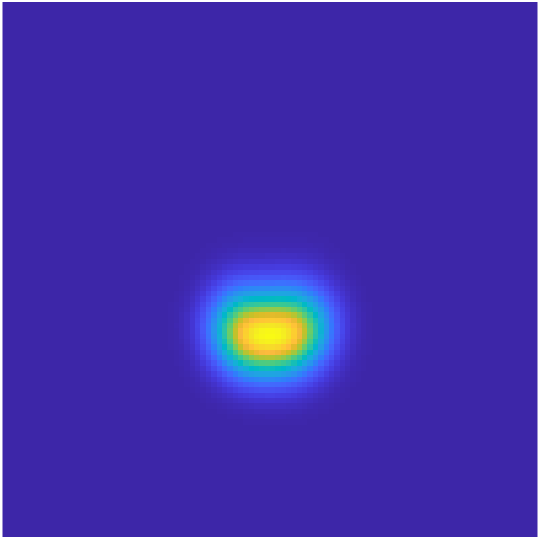}\\
\centering {$F_{7*}(P_0)$}
\end{minipage}\hfill
\begin{minipage}{0.16\linewidth}
\includegraphics[width=1\linewidth]{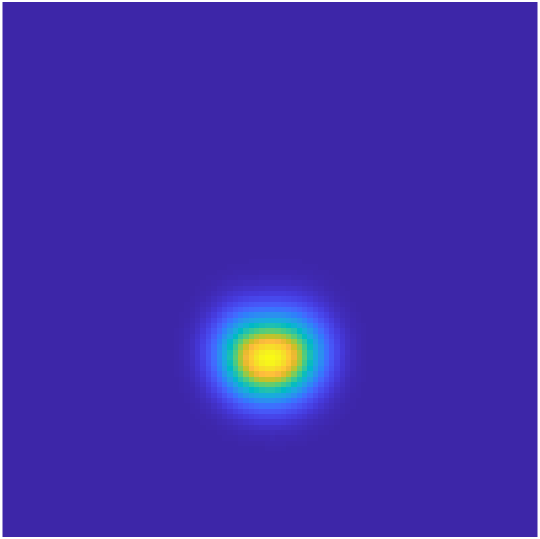}\\
\centering {$F_{8*}(P_0)$}
\end{minipage}\hfill
\begin{minipage}{0.16\linewidth}
\includegraphics[width=1\linewidth]{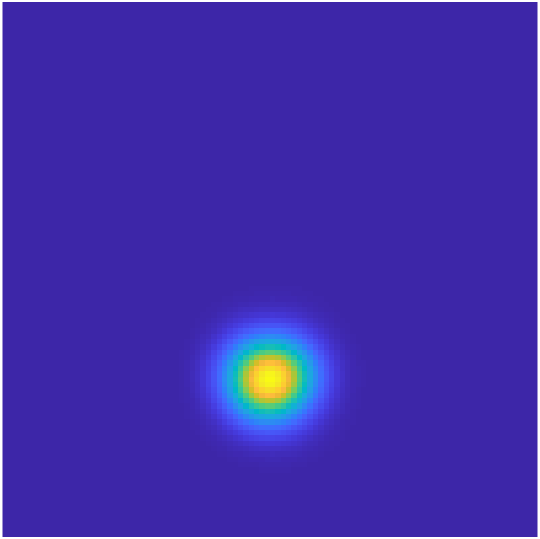}\\
\centering {$F_{9*}(P_0)$}
\end{minipage}\hfill
\begin{minipage}{0.16\linewidth}
\includegraphics[width=1\linewidth]{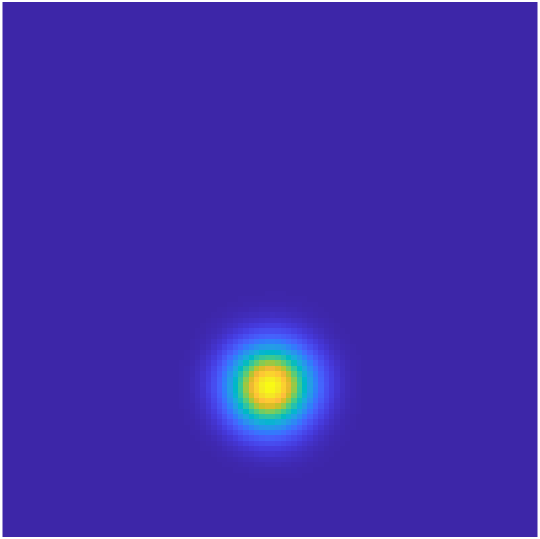}\\
\centering {$F_{10*}(P_0)$}
\end{minipage}\hfill
\begin{minipage}{0.16\linewidth}
\includegraphics[width=1\linewidth]{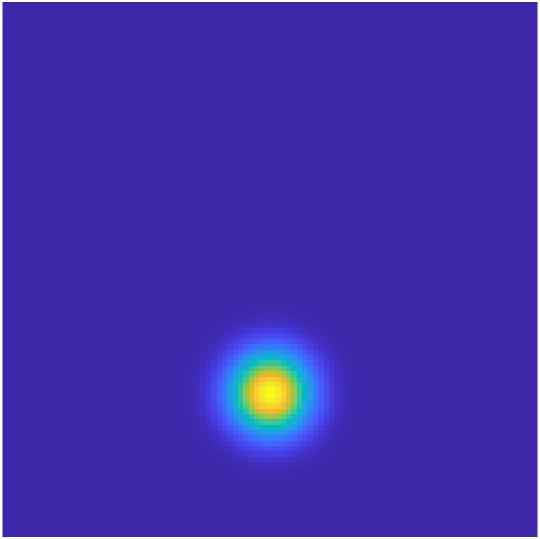}\\
\centering {$P_1$}
\end{minipage}\hfill

\caption{Evolution of the density for the crowd motion problem~\eqref{crowd_motion_MFG} in 2D.}
\label{crowd_2D_evo}
\end{figure}

\begin{figure}[H]
\centering
\begin{minipage}{0.16\linewidth}
\includegraphics[width=1\linewidth]{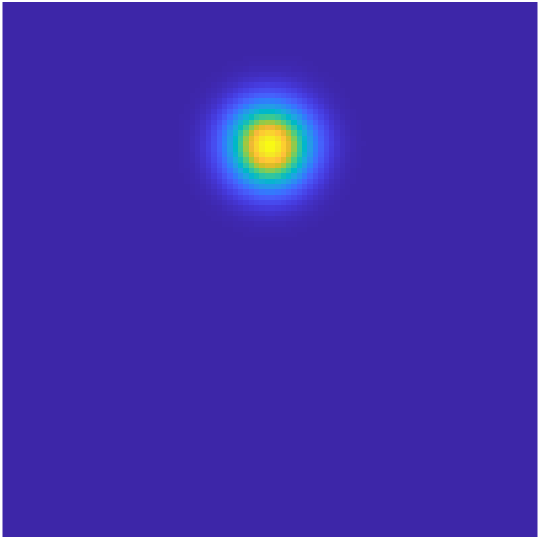}\\
\centering {$P_0$}
\end{minipage}\hfill
\begin{minipage}{0.16\linewidth}
\includegraphics[width=1\linewidth]{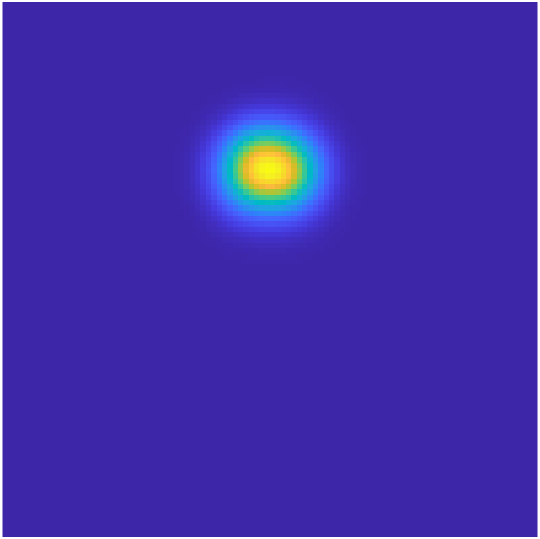}\\
\centering {$F_{1*}(P_0)$}
\end{minipage}\hfill
\begin{minipage}{0.16\linewidth}
\includegraphics[width=1\linewidth]{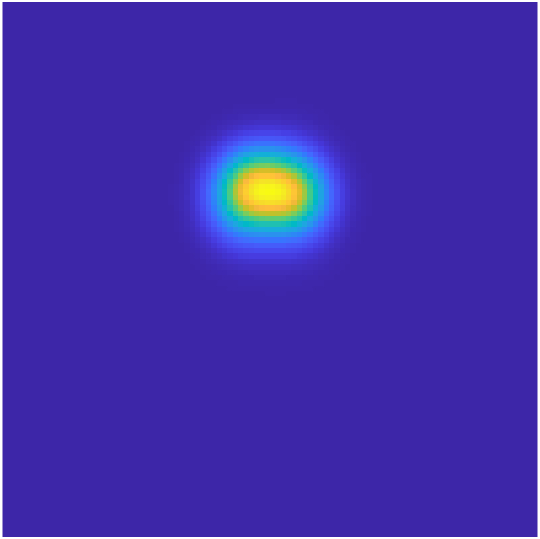}\\
\centering {$F_{2*}(P_0)$}
\end{minipage}\hfill
\begin{minipage}{0.16\linewidth}
\includegraphics[width=1\linewidth]{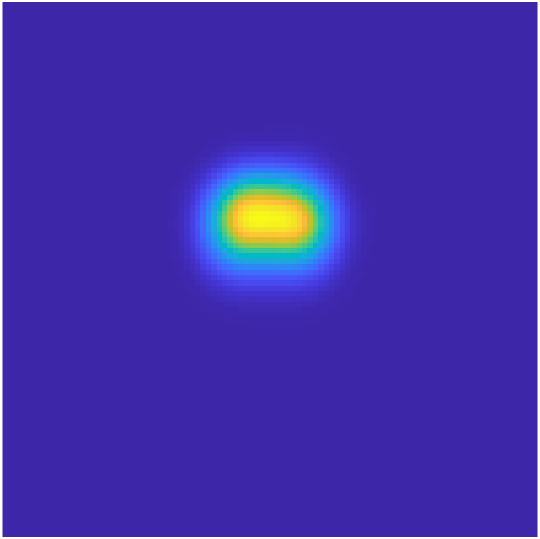}\\
\centering {$F_{3*}(P_0)$}
\end{minipage}\hfill
\begin{minipage}{0.16\linewidth}
\includegraphics[width=1\linewidth]{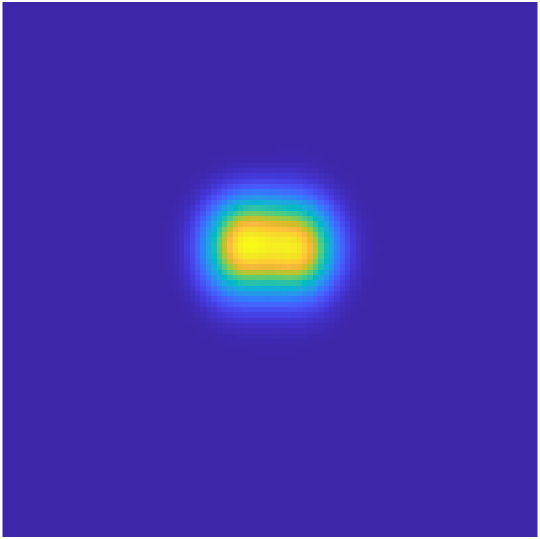}\\
\centering {$F_{4*}(P_0)$}
\end{minipage}\hfill
\begin{minipage}{0.16\linewidth}
\includegraphics[width=1\linewidth]{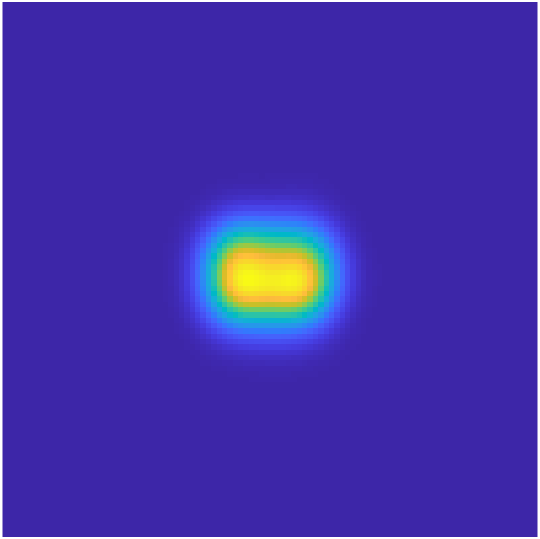}\\
\centering {$F_{5*}(P_0)$}
\end{minipage}\hfill
\begin{minipage}{0.16\linewidth}
\includegraphics[width=1\linewidth]{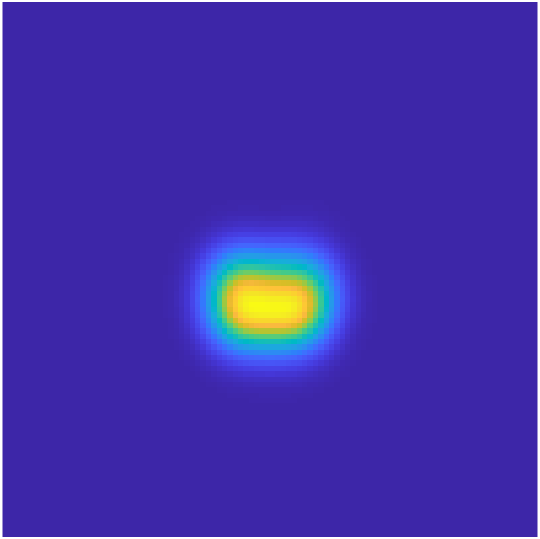}\\
\centering {$F_{6*}(P_0)$}
\end{minipage}\hfill
\begin{minipage}{0.16\linewidth}
\includegraphics[width=1\linewidth]{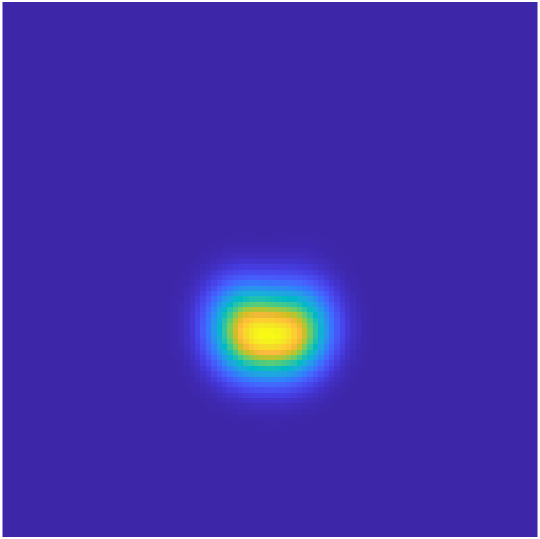}\\
\centering {$F_{7*}(P_0)$}
\end{minipage}\hfill
\begin{minipage}{0.16\linewidth}
\includegraphics[width=1\linewidth]{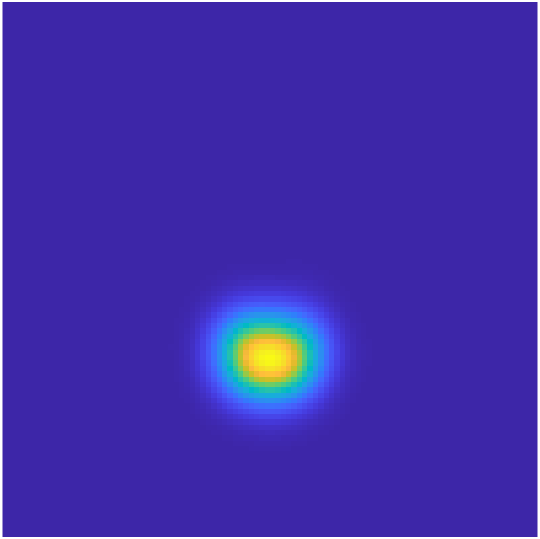}\\
\centering {$F_{8*}(P_0)$}
\end{minipage}\hfill
\begin{minipage}{0.16\linewidth}
\includegraphics[width=1\linewidth]{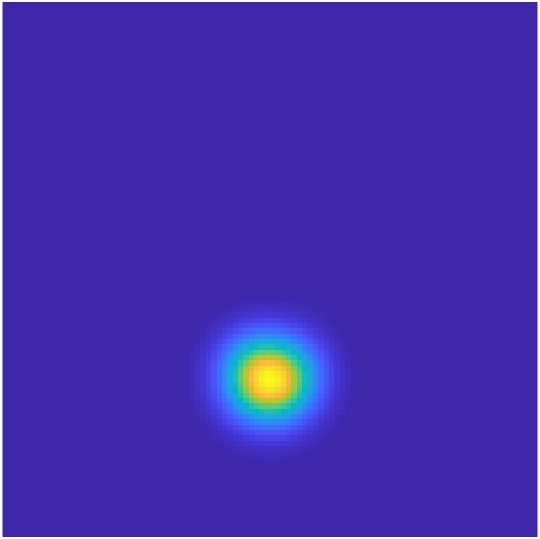}\\
\centering {$F_{9*}(P_0)$}
\end{minipage}\hfill
\begin{minipage}{0.16\linewidth}
\includegraphics[width=1\linewidth]{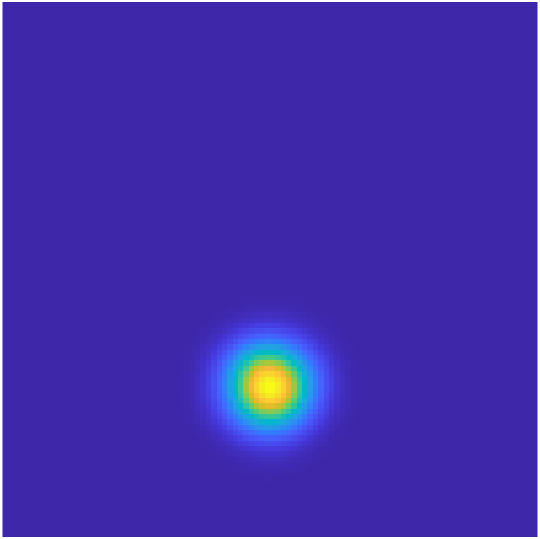}\\
\centering {$F_{10*}(P_0)$}
\end{minipage}\hfill
\begin{minipage}{0.16\linewidth}
\includegraphics[width=1\linewidth]{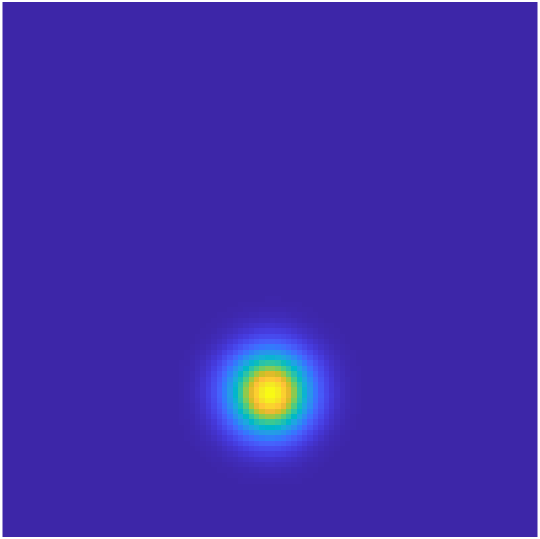}\\
\centering {$P_1$}
\end{minipage}\hfill

\caption{Evolution of the density for the crowd motion problem~\eqref{crowd_motion_MFG} in 10D.}
\label{crowd_10D_evo}
\end{figure}

\begin{figure}[H]
\centering
\begin{minipage}{0.16\linewidth}
\includegraphics[width=1\linewidth]{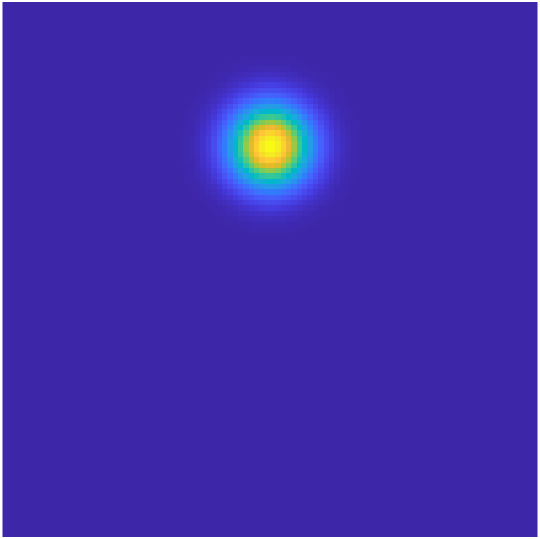}\\
\centering {$P_0$}
\end{minipage}\hfill
\begin{minipage}{0.16\linewidth}
\includegraphics[width=1\linewidth]{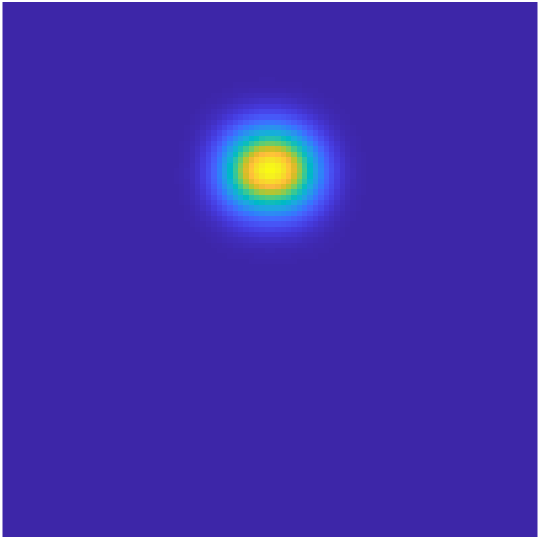}\\
\centering {$F_{1*}(P_0)$}
\end{minipage}\hfill
\begin{minipage}{0.16\linewidth}
\includegraphics[width=1\linewidth]{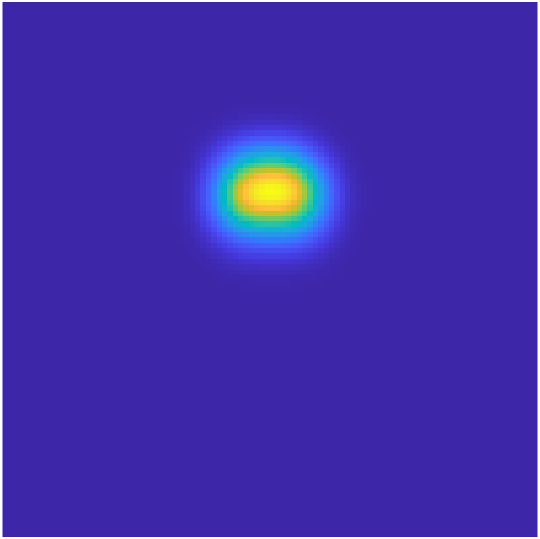}\\
\centering {$F_{2*}(P_0)$}
\end{minipage}\hfill
\begin{minipage}{0.16\linewidth}
\includegraphics[width=1\linewidth]{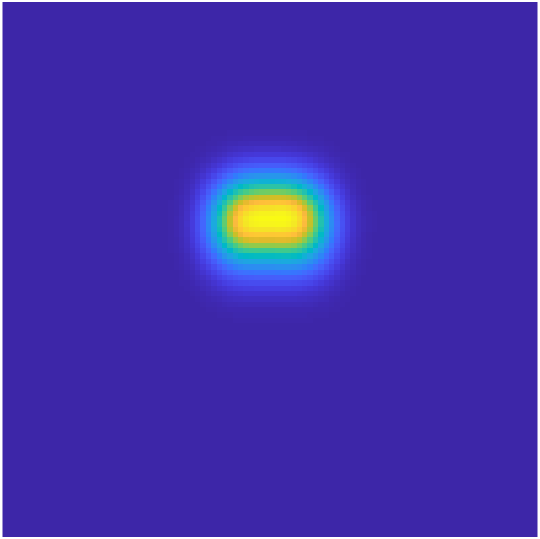}\\
\centering {$F_{3*}(P_0)$}
\end{minipage}\hfill
\begin{minipage}{0.16\linewidth}
\includegraphics[width=1\linewidth]{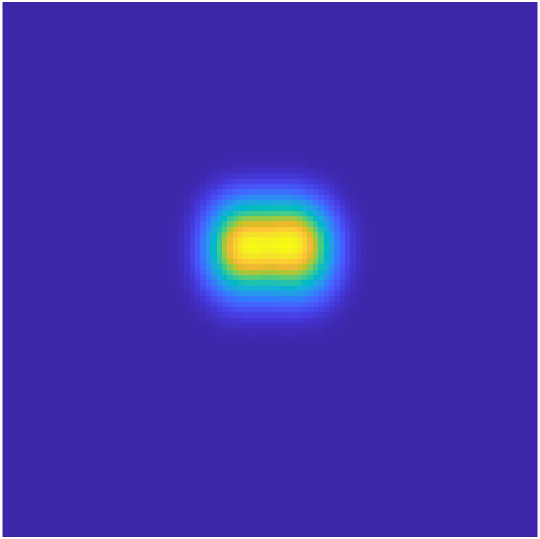}\\
\centering {$F_{4*}(P_0)$}
\end{minipage}\hfill
\begin{minipage}{0.16\linewidth}
\includegraphics[width=1\linewidth]{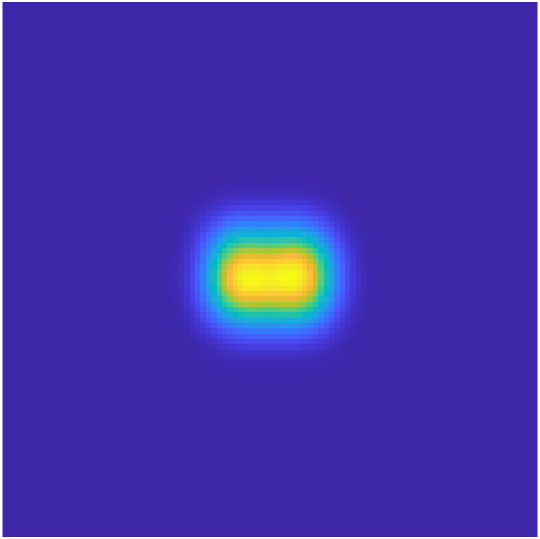}\\
\centering {$F_{5*}(P_0)$}
\end{minipage}\hfill
\begin{minipage}{0.16\linewidth}
\includegraphics[width=1\linewidth]{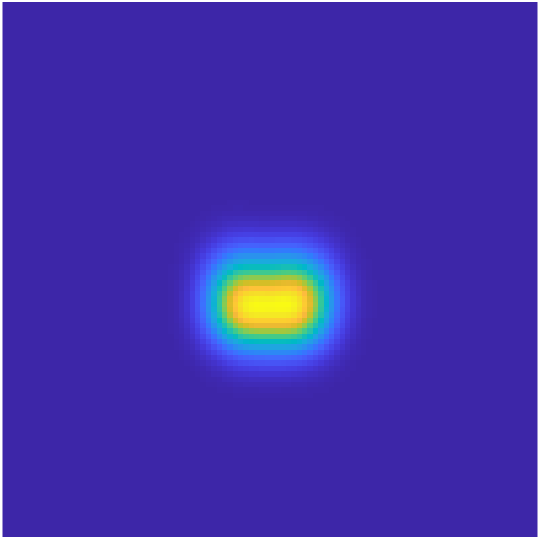}\\
\centering {$F_{6*}(P_0)$}
\end{minipage}\hfill
\begin{minipage}{0.16\linewidth}
\includegraphics[width=1\linewidth]{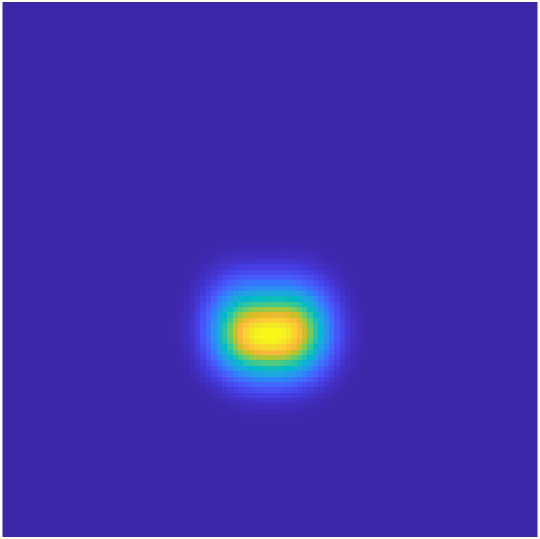}\\
\centering {$F_{7*}(P_0)$}
\end{minipage}\hfill
\begin{minipage}{0.16\linewidth}
\includegraphics[width=1\linewidth]{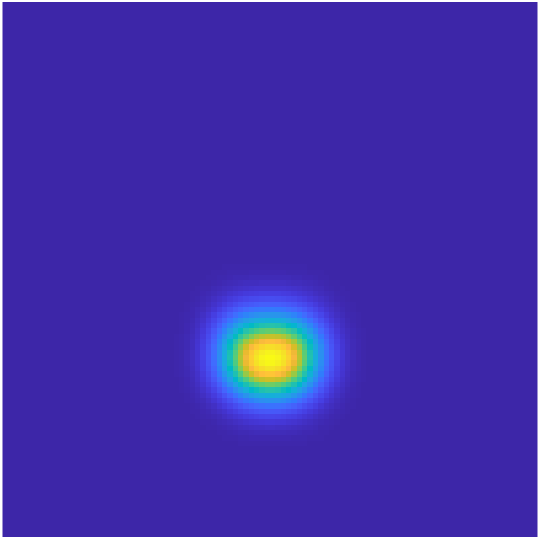}\\
\centering {$F_{8*}(P_0)$}
\end{minipage}\hfill
\begin{minipage}{0.16\linewidth}
\includegraphics[width=1\linewidth]{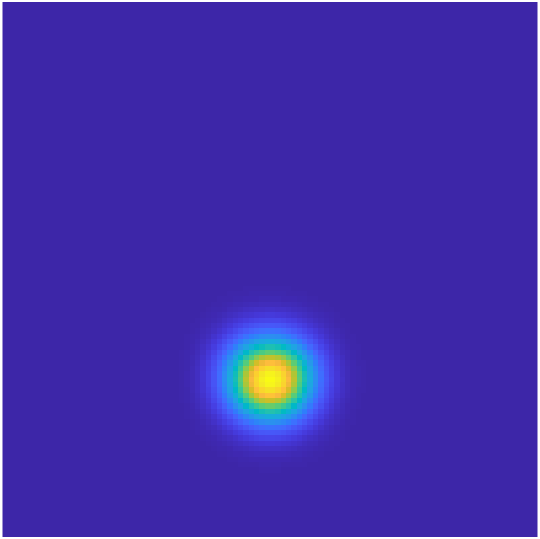}\\
\centering {$F_{9*}(P_0)$}
\end{minipage}\hfill
\begin{minipage}{0.16\linewidth}
\includegraphics[width=1\linewidth]{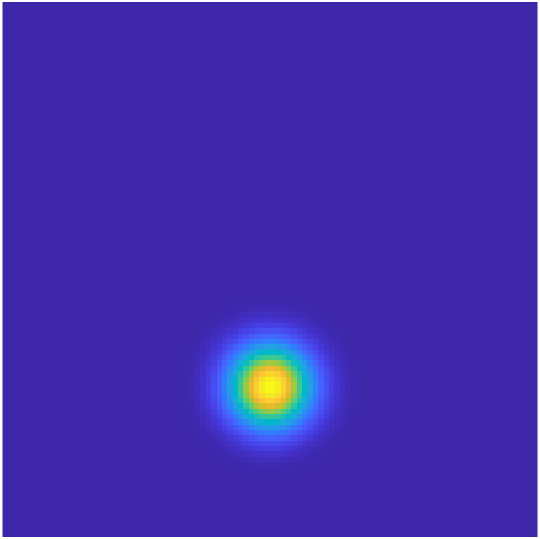}\\
\centering {$F_{10*}(P_0)$}
\end{minipage}\hfill
\begin{minipage}{0.16\linewidth}
\includegraphics[width=1\linewidth]{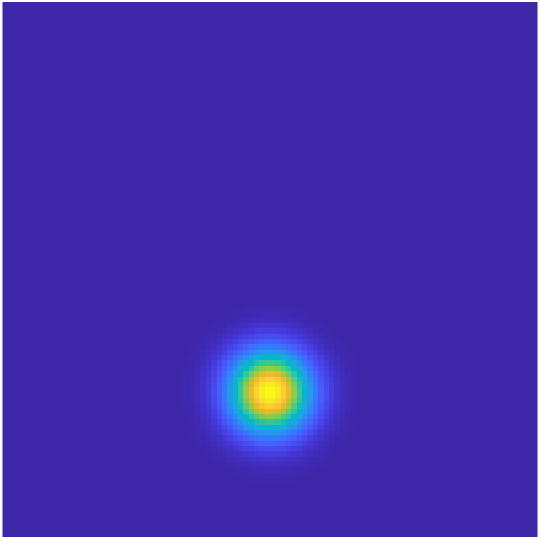}\\
\centering {$P_1$}
\end{minipage}\hfill

\caption{Evolution of the density for the crowd motion problem~\eqref{crowd_motion_MFG} in 50D.}
\label{crowd_50D_evo}
\end{figure}

\begin{figure}[H]
\centering
\begin{minipage}{0.16\linewidth}
\includegraphics[width=1\linewidth]{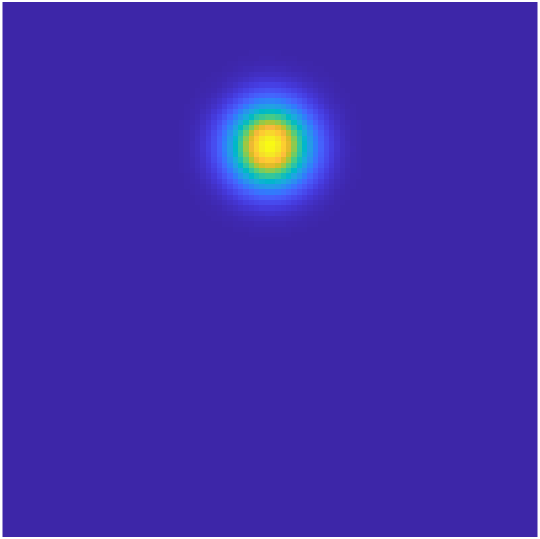}\\
\centering {$P_0$}
\end{minipage}\hfill
\begin{minipage}{0.16\linewidth}
\includegraphics[width=1\linewidth]{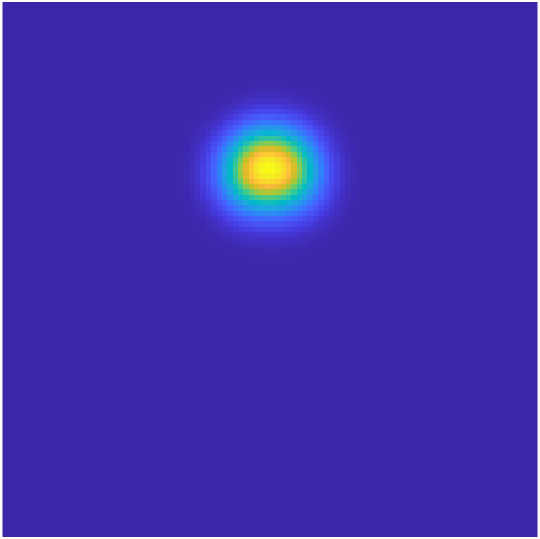}\\
\centering {$F_{1*}(P_0)$}
\end{minipage}\hfill
\begin{minipage}{0.16\linewidth}
\includegraphics[width=1\linewidth]{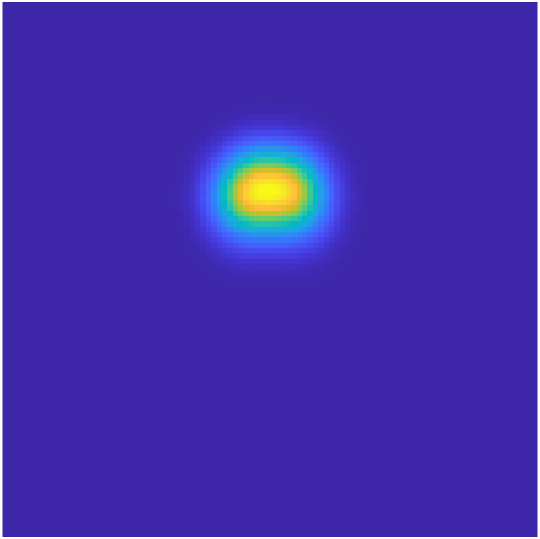}\\
\centering {$F_{2*}(P_0)$}
\end{minipage}\hfill
\begin{minipage}{0.16\linewidth}
\includegraphics[width=1\linewidth]{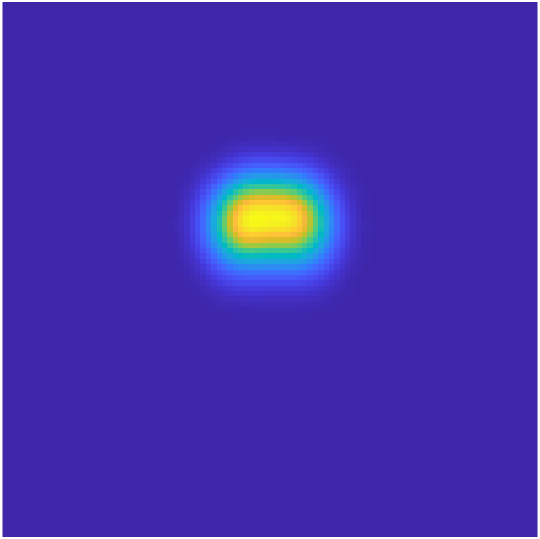}\\
\centering {$F_{3*}(P_0)$}
\end{minipage}\hfill
\begin{minipage}{0.16\linewidth}
\includegraphics[width=1\linewidth]{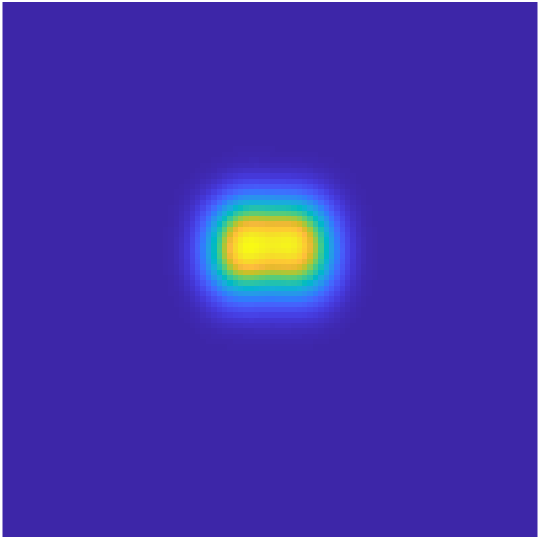}\\
\centering {$F_{4*}(P_0)$}
\end{minipage}\hfill
\begin{minipage}{0.16\linewidth}
\includegraphics[width=1\linewidth]{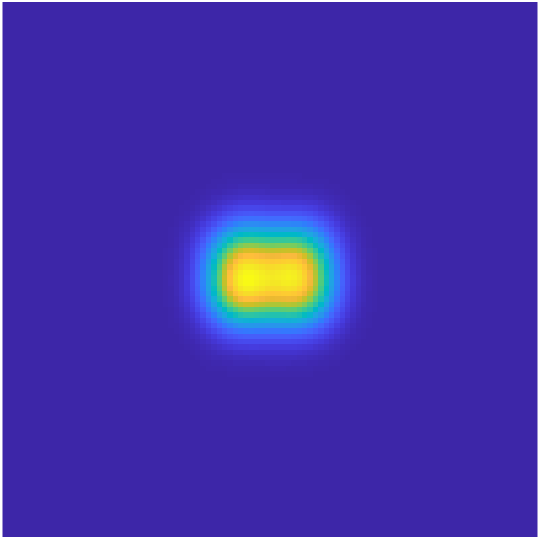}\\
\centering {$F_{5*}(P_0)$}
\end{minipage}\hfill
\begin{minipage}{0.16\linewidth}
\includegraphics[width=1\linewidth]{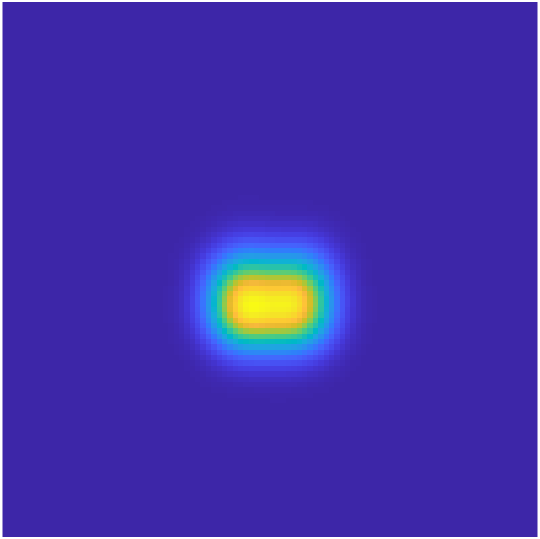}\\
\centering {$F_{6*}(P_0)$}
\end{minipage}\hfill
\begin{minipage}{0.16\linewidth}
\includegraphics[width=1\linewidth]{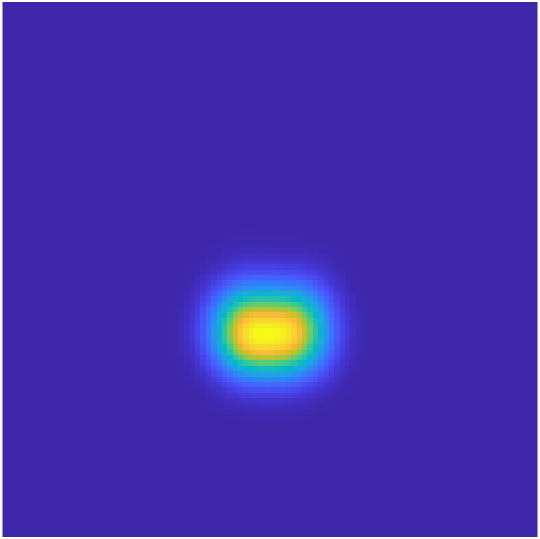}\\
\centering {$F_{7*}(P_0)$}
\end{minipage}\hfill
\begin{minipage}{0.16\linewidth}
\includegraphics[width=1\linewidth]{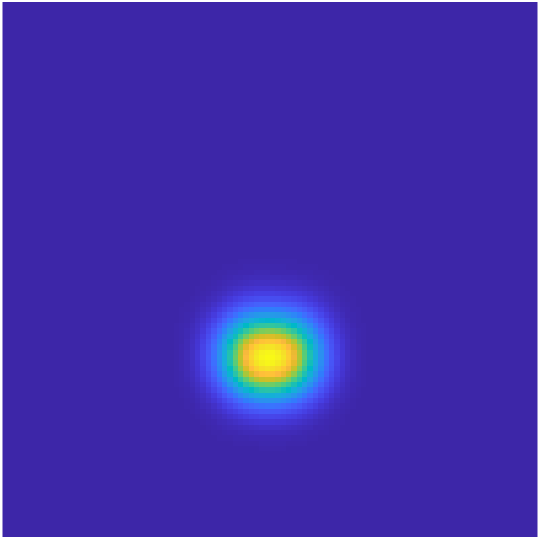}\\
\centering {$F_{8*}(P_0)$}
\end{minipage}\hfill
\begin{minipage}{0.16\linewidth}
\includegraphics[width=1\linewidth]{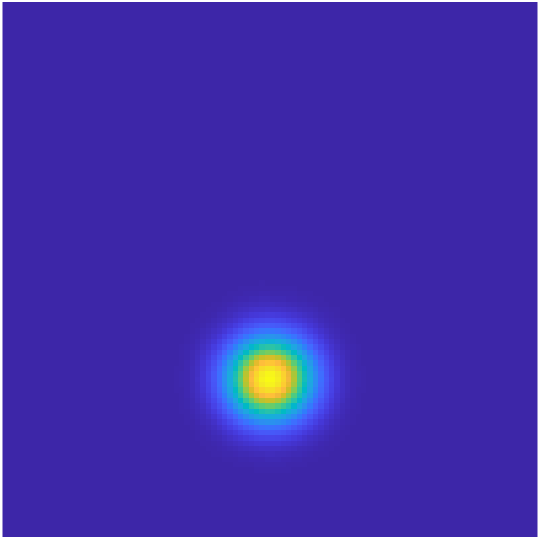}\\
\centering {$F_{9*}(P_0)$}
\end{minipage}\hfill
\begin{minipage}{0.16\linewidth}
\includegraphics[width=1\linewidth]{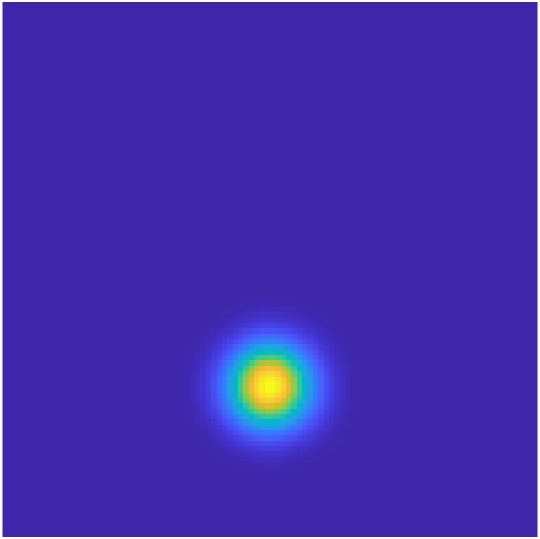}\\
\centering {$F_{10*}(P_0)$}
\end{minipage}\hfill
\begin{minipage}{0.16\linewidth}
\includegraphics[width=1\linewidth]{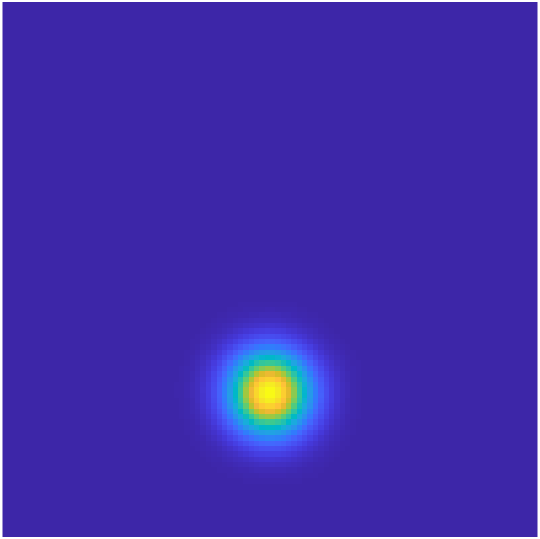}\\
\centering {$P_1$}
\end{minipage}\hfill

\caption{Evolution of the density for the crowd motion problem~\eqref{crowd_motion_MFG} in 100D.}
\label{crowd_100D_evo}
\end{figure}

\subsubsection{Different Levels of Penalty on Conflicts with the Obstacle}

In addition to the sampled trajectories provided in Figure~\ref{crowd_traj}, we show the density evolution as well as the computed cost values for different aversion preferences here. All experiments are conducted at $d=10$ with identical hyperparameter settings as the one used in Table~\ref{crowd_motion_table}, except for the weights on the cost terms. The level of aversion is manifested through the outward-curving behavior when the densities are transported near the origin.

\begin{table}[H]
\caption{Results for crowd motion at different weights assigned to $\mathcal{I}$. $\lambda_{\mathcal{I}}$: weight for $\mathcal{I}$; $L$: transport cost; $\mathcal{I}$: interaction cost; $\mathcal{M}$: terminal cost.}
\label{crowd_table_avoidance_preference}
\begin{center}
\begin{small}
\begin{sc}
\begin{tabular}{lccccccccccc}
\toprule
$\lambda_{\mathcal{I}}$ & $L$ & $\mathcal{I}$ & $\mathcal{M}$\\
\midrule
 0.2 & 32.8511 $\pm$ 0.0060  & 2.1274 $\pm$ 0.0114 & 0.0436 $\pm$ 0.0059 \\
 0.5 & 35.5121 $\pm$ 0.0134  & 1.3463 $\pm$ 0.0096 & 0.0424 $\pm$ 0.0058 \\
 1.0 & 39.9230 $\pm$ 0.0223  & 0.7076 $\pm$ 0.0053 & 0.0430 $\pm$ 0.0058 \\
\bottomrule
\end{tabular}
\end{sc}
\end{small}
\end{center}
\end{table}

\begin{figure}[H]
\centering
\begin{minipage}{0.16\linewidth}
\includegraphics[width=1\linewidth]{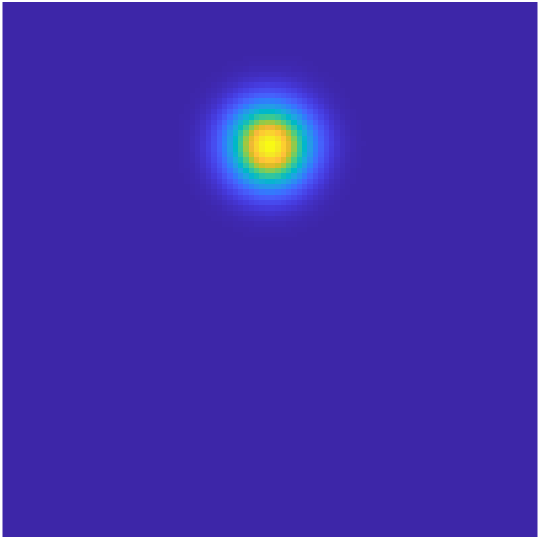}\\
\centering {$P_0$}
\end{minipage}\hfill
\begin{minipage}{0.16\linewidth}
\includegraphics[width=1\linewidth]{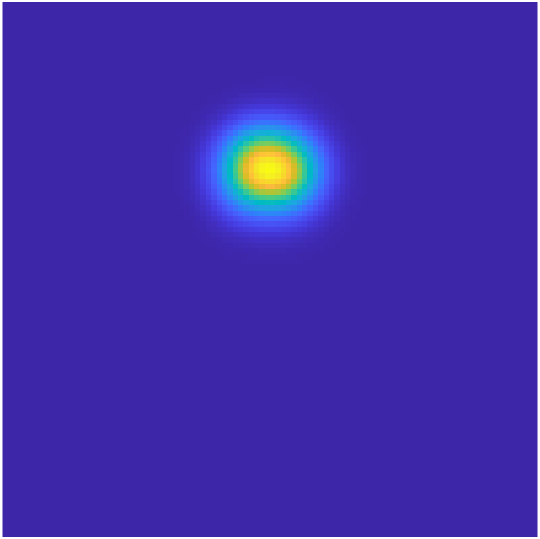}\\
\centering {$F_{1*}(P_0)$}
\end{minipage}\hfill
\begin{minipage}{0.16\linewidth}
\includegraphics[width=1\linewidth]{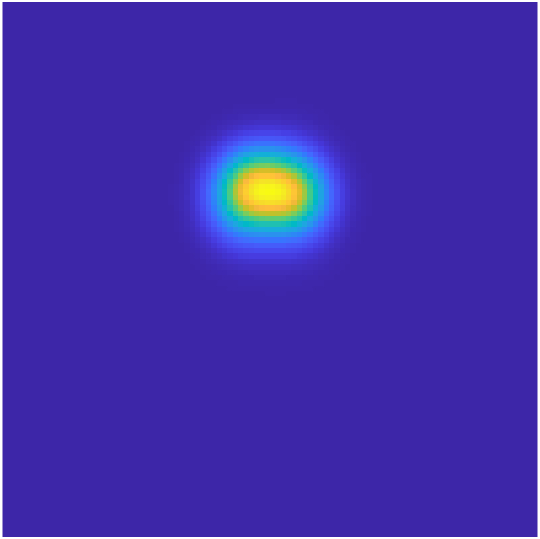}\\
\centering {$F_{2*}(P_0)$}
\end{minipage}\hfill
\begin{minipage}{0.16\linewidth}
\includegraphics[width=1\linewidth]{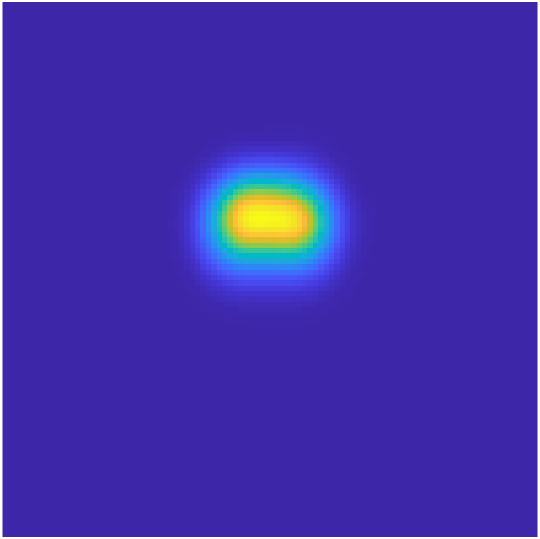}\\
\centering {$F_{3*}(P_0)$}
\end{minipage}\hfill
\begin{minipage}{0.16\linewidth}
\includegraphics[width=1\linewidth]{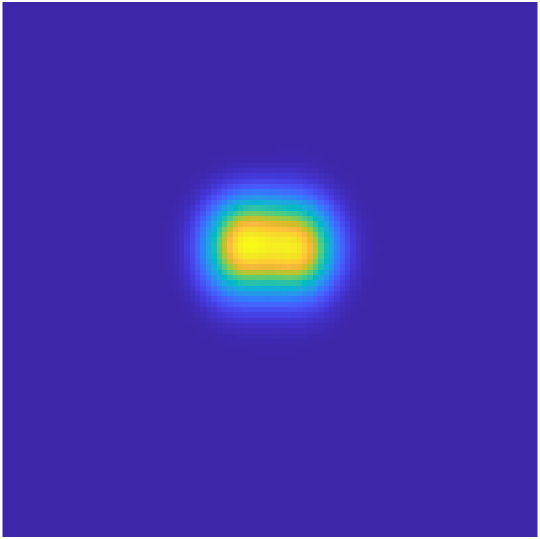}\\
\centering {$F_{4*}(P_0)$}
\end{minipage}\hfill
\begin{minipage}{0.16\linewidth}
\includegraphics[width=1\linewidth]{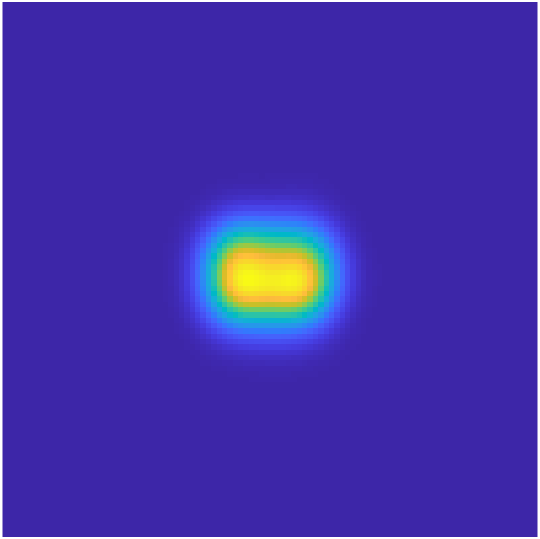}\\
\centering {$F_{5*}(P_0)$}
\end{minipage}\hfill
\begin{minipage}{0.16\linewidth}
\includegraphics[width=1\linewidth]{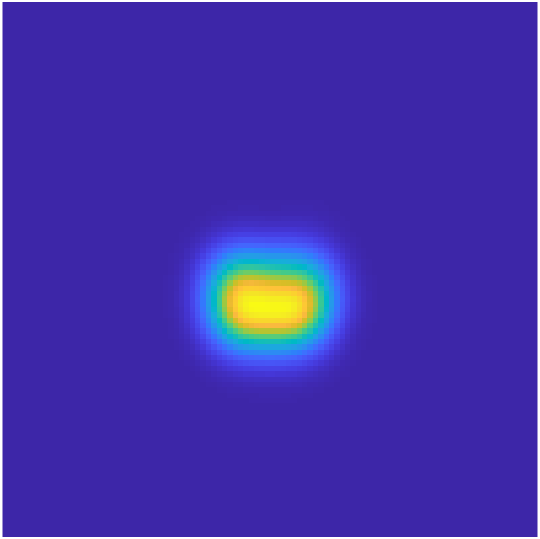}\\
\centering {$F_{6*}(P_0)$}
\end{minipage}\hfill
\begin{minipage}{0.16\linewidth}
\includegraphics[width=1\linewidth]{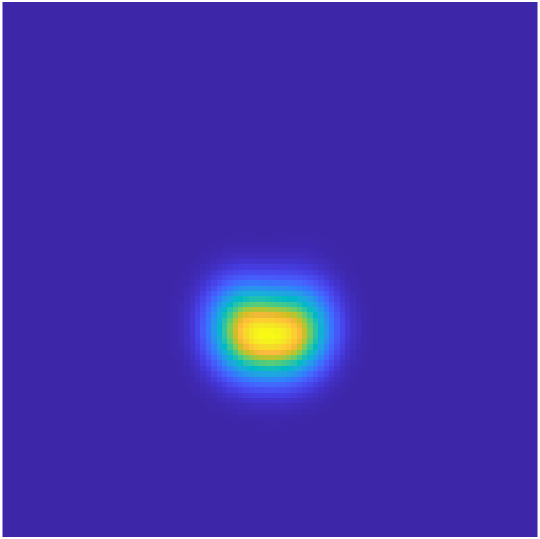}\\
\centering {$F_{7*}(P_0)$}
\end{minipage}\hfill
\begin{minipage}{0.16\linewidth}
\includegraphics[width=1\linewidth]{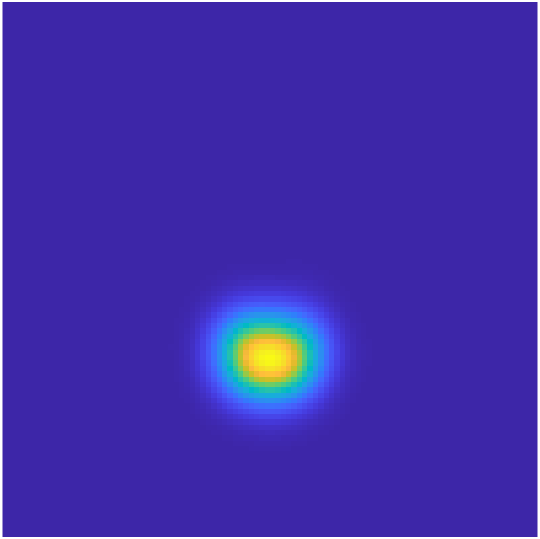}\\
\centering {$F_{8*}(P_0)$}
\end{minipage}\hfill
\begin{minipage}{0.16\linewidth}
\includegraphics[width=1\linewidth]{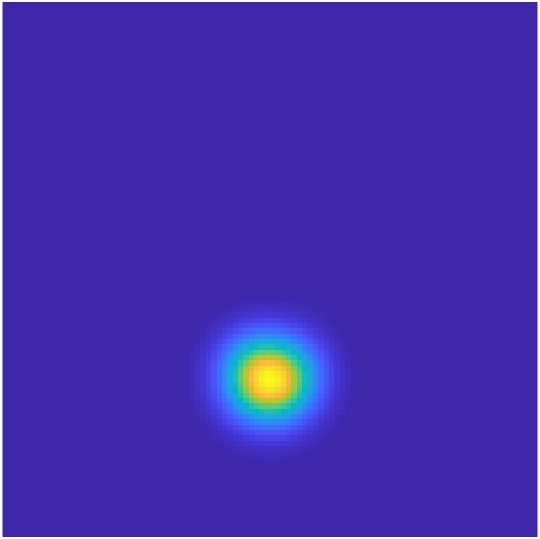}\\
\centering {$F_{9*}(P_0)$}
\end{minipage}\hfill
\begin{minipage}{0.16\linewidth}
\includegraphics[width=1\linewidth]{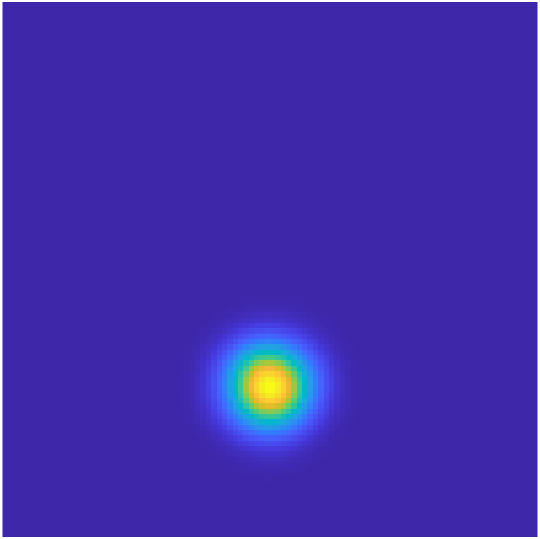}\\
\centering {$F_{10*}(P_0)$}
\end{minipage}\hfill
\begin{minipage}{0.16\linewidth}
\includegraphics[width=1\linewidth]{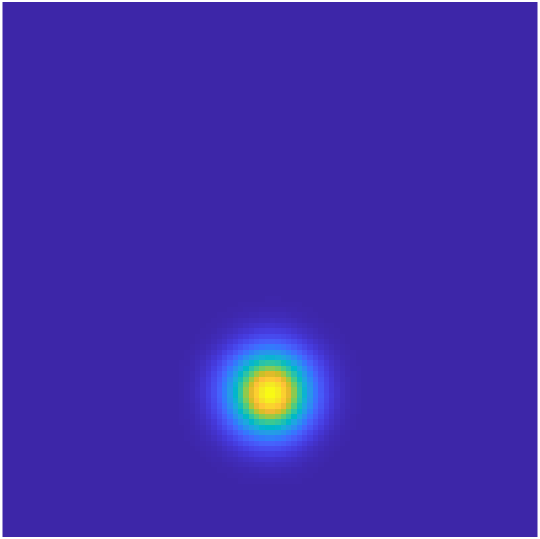}\\
\centering {$P_1$}
\end{minipage}\hfill

\caption{Evolution of the density for minor ($\lambda_{\mathcal{I}}=0.2$) penalty on conflicts with the obstacle.}
\label{crowd_10D_minor_evo}
\end{figure}

\begin{figure}[H]
\centering
\begin{minipage}{0.16\linewidth}
\includegraphics[width=1\linewidth]{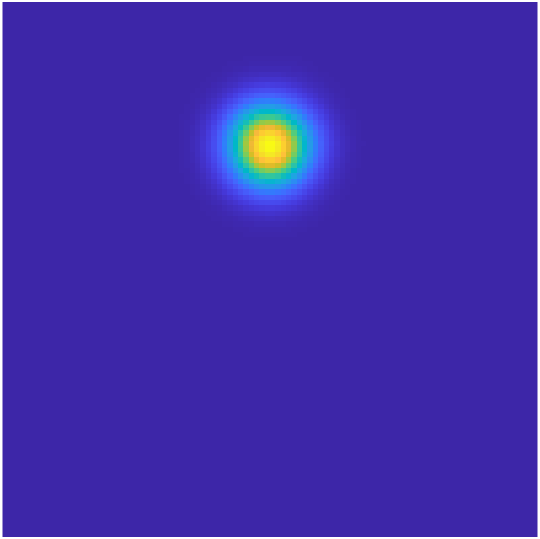}\\
\centering {$P_0$}
\end{minipage}\hfill
\begin{minipage}{0.16\linewidth}
\includegraphics[width=1\linewidth]{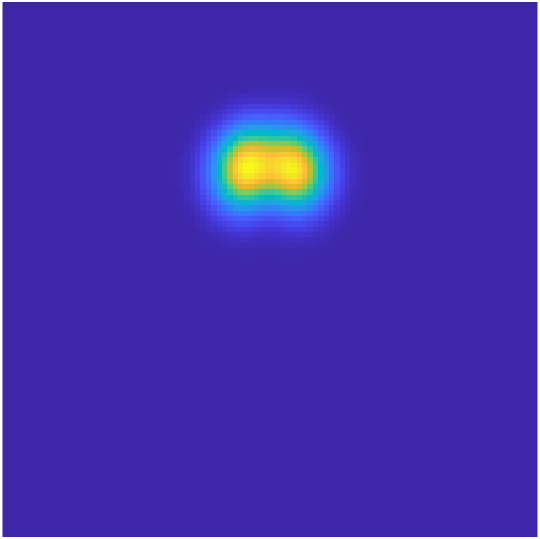}\\
\centering {$F_{1*}(P_0)$}
\end{minipage}\hfill
\begin{minipage}{0.16\linewidth}
\includegraphics[width=1\linewidth]{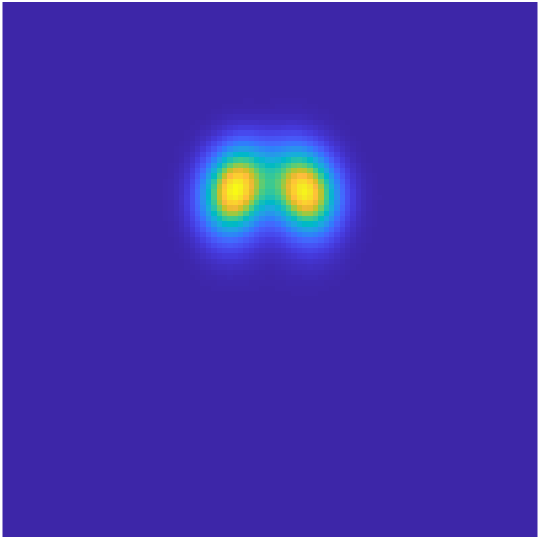}\\
\centering {$F_{2*}(P_0)$}
\end{minipage}\hfill
\begin{minipage}{0.16\linewidth}
\includegraphics[width=1\linewidth]{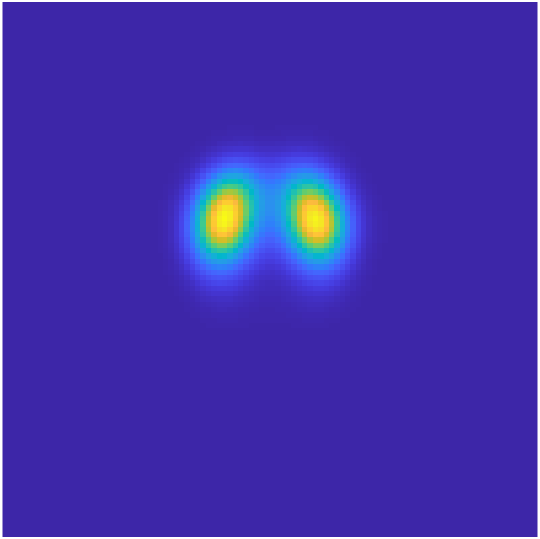}\\
\centering {$F_{3*}(P_0)$}
\end{minipage}\hfill
\begin{minipage}{0.16\linewidth}
\includegraphics[width=1\linewidth]{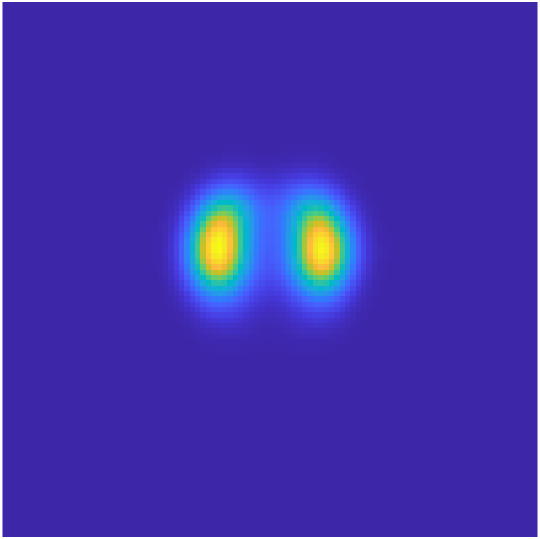}\\
\centering {$F_{4*}(P_0)$}
\end{minipage}\hfill
\begin{minipage}{0.16\linewidth}
\includegraphics[width=1\linewidth]{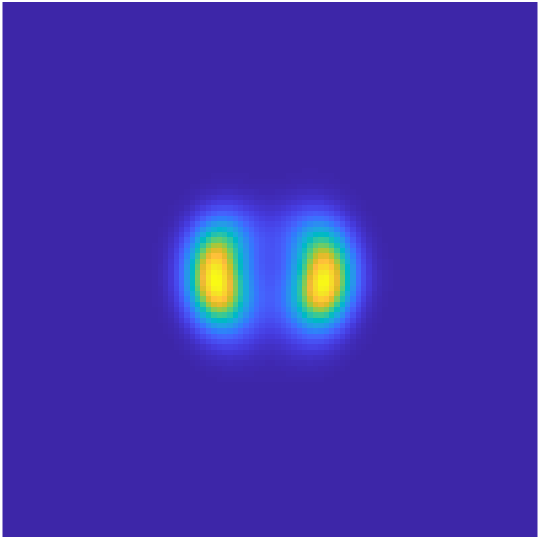}\\
\centering {$F_{5*}(P_0)$}
\end{minipage}\hfill
\begin{minipage}{0.16\linewidth}
\includegraphics[width=1\linewidth]{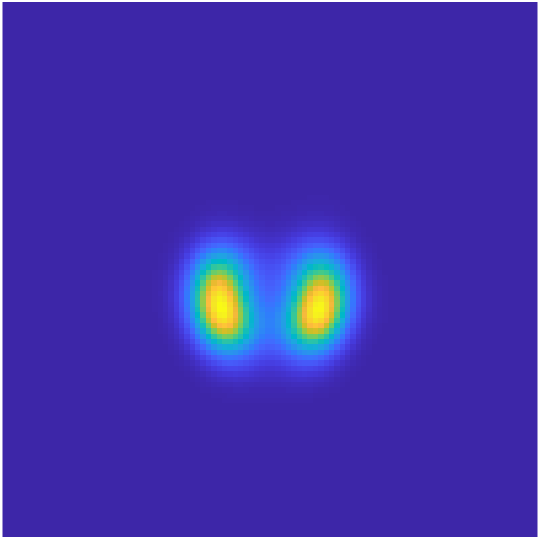}\\
\centering {$F_{6*}(P_0)$}
\end{minipage}\hfill
\begin{minipage}{0.16\linewidth}
\includegraphics[width=1\linewidth]{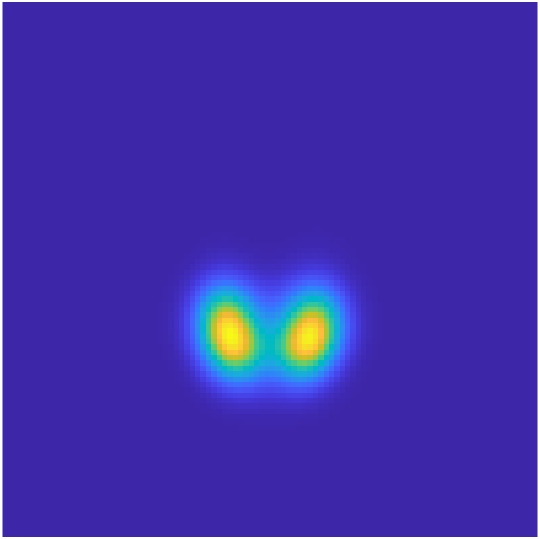}\\
\centering {$F_{7*}(P_0)$}
\end{minipage}\hfill
\begin{minipage}{0.16\linewidth}
\includegraphics[width=1\linewidth]{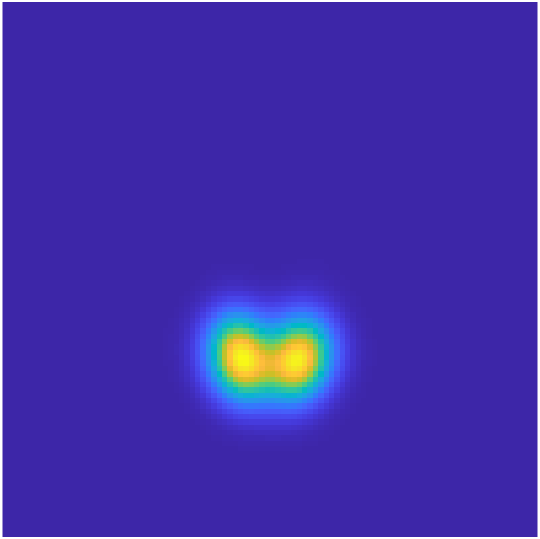}\\
\centering {$F_{8*}(P_0)$}
\end{minipage}\hfill
\begin{minipage}{0.16\linewidth}
\includegraphics[width=1\linewidth]{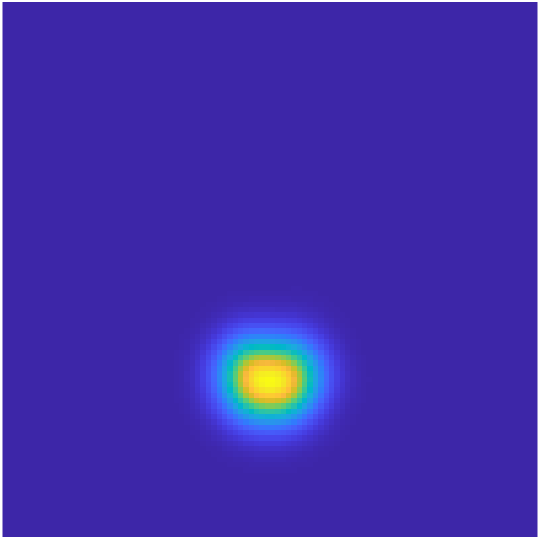}\\
\centering {$F_{9*}(P_0)$}
\end{minipage}\hfill
\begin{minipage}{0.16\linewidth}
\includegraphics[width=1\linewidth]{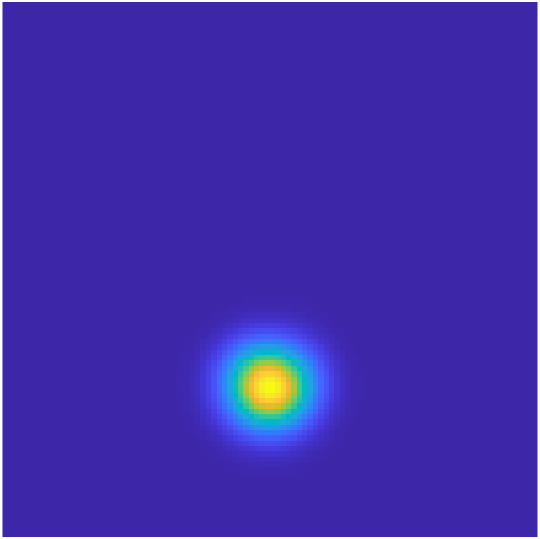}\\
\centering {$F_{10*}(P_0)$}
\end{minipage}\hfill
\begin{minipage}{0.16\linewidth}
\includegraphics[width=1\linewidth]{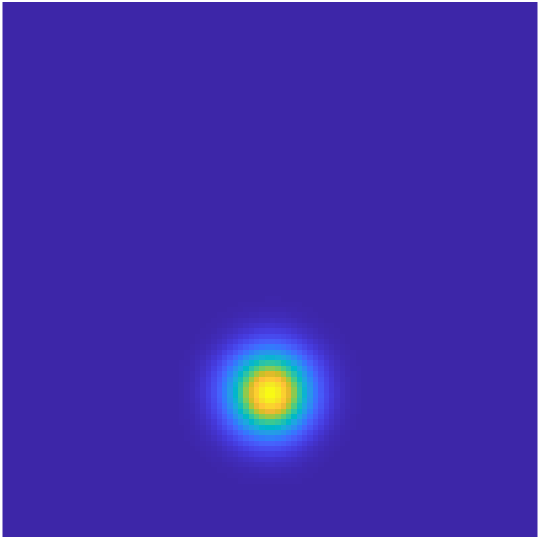}\\
\centering {$F_{10*}(P_0)$}
\end{minipage}\hfill

\caption{Evolution of the density for moderate ($\lambda_{\mathcal{I}}=0.5$) penalty on conflicts with the obstacle.}
\label{crowd_10D_morderate_evo}
\end{figure}

\begin{figure}[H]
\centering
\begin{minipage}{0.16\linewidth}
\includegraphics[width=1\linewidth]{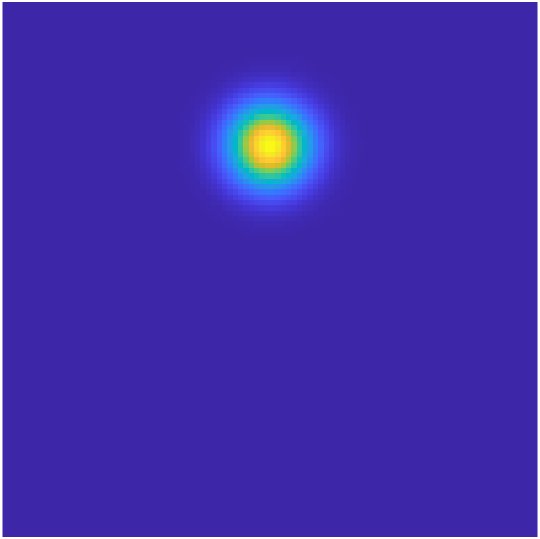}\\
\centering {$P_0$}
\end{minipage}\hfill
\begin{minipage}{0.16\linewidth}
\includegraphics[width=1\linewidth]{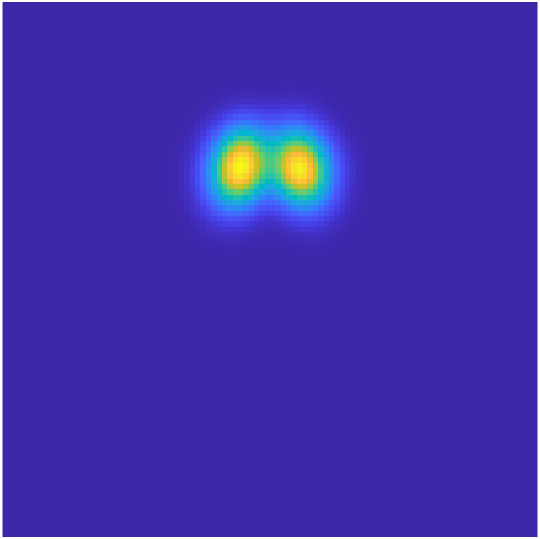}\\
\centering {$F_{1*}(P_0)$}
\end{minipage}\hfill
\begin{minipage}{0.16\linewidth}
\includegraphics[width=1\linewidth]{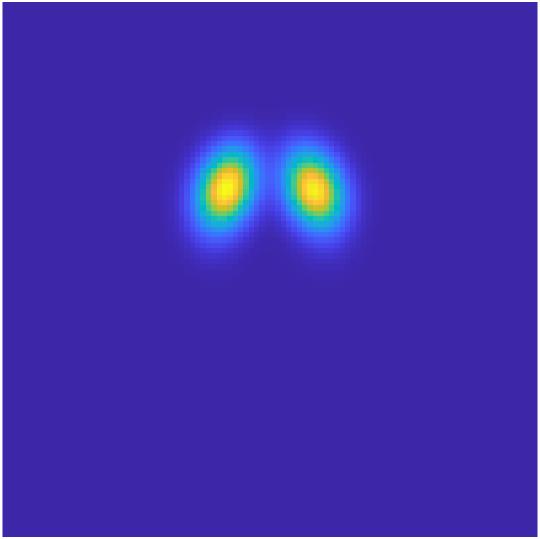}\\
\centering {$F_{2*}(P_0)$}
\end{minipage}\hfill
\begin{minipage}{0.16\linewidth}
\includegraphics[width=1\linewidth]{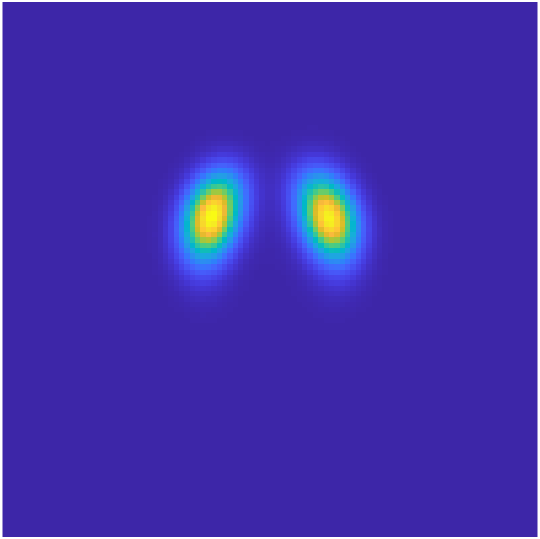}\\
\centering {$F_{3*}(P_0)$}
\end{minipage}\hfill
\begin{minipage}{0.16\linewidth}
\includegraphics[width=1\linewidth]{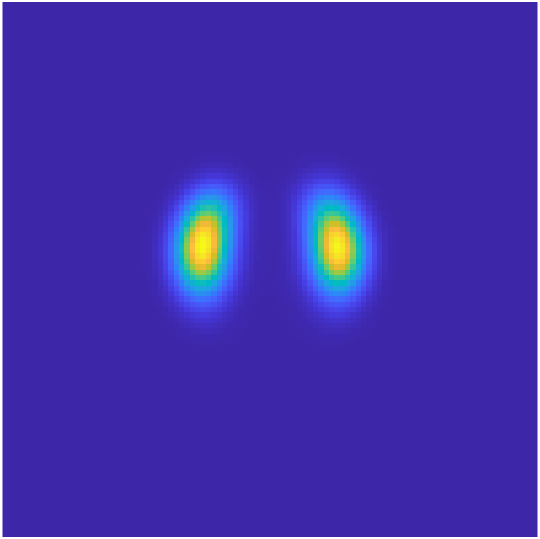}\\
\centering {$F_{4*}(P_0)$}
\end{minipage}\hfill
\begin{minipage}{0.16\linewidth}
\includegraphics[width=1\linewidth]{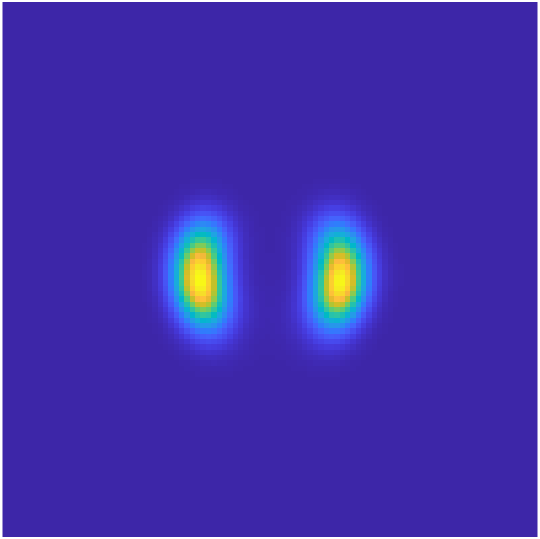}\\
\centering {$F_{5*}(P_0)$}
\end{minipage}\hfill
\begin{minipage}{0.16\linewidth}
\includegraphics[width=1\linewidth]{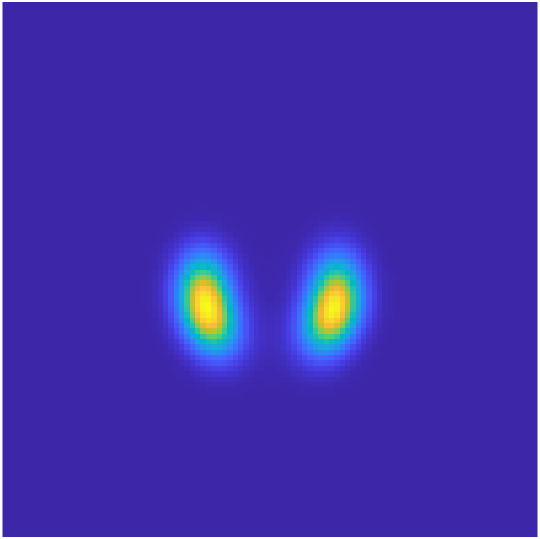}\\
\centering {$F_{6*}(P_0)$}
\end{minipage}\hfill
\begin{minipage}{0.16\linewidth}
\includegraphics[width=1\linewidth]{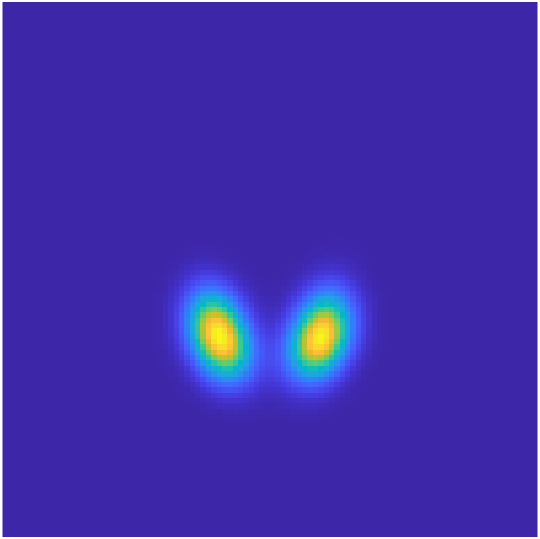}\\
\centering {$F_{7*}(P_0)$}
\end{minipage}\hfill
\begin{minipage}{0.16\linewidth}
\includegraphics[width=1\linewidth]{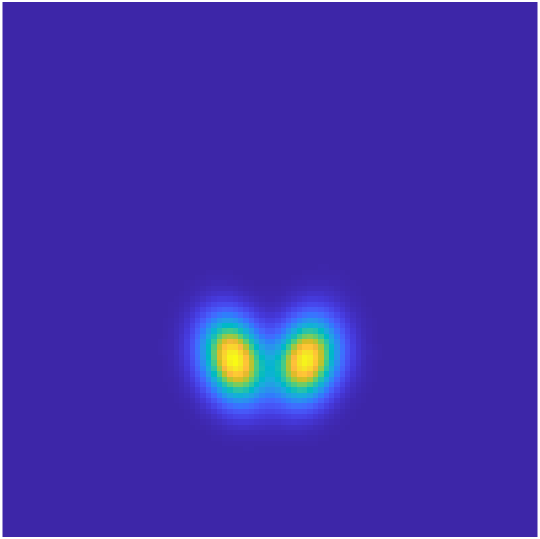}\\
\centering {$F_{8*}(P_0)$}
\end{minipage}\hfill
\begin{minipage}{0.16\linewidth}
\includegraphics[width=1\linewidth]{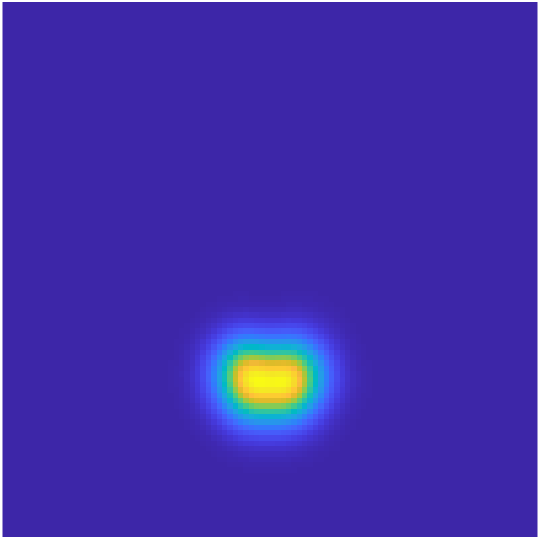}\\
\centering {$F_{9*}(P_0)$}
\end{minipage}\hfill
\begin{minipage}{0.16\linewidth}
\includegraphics[width=1\linewidth]{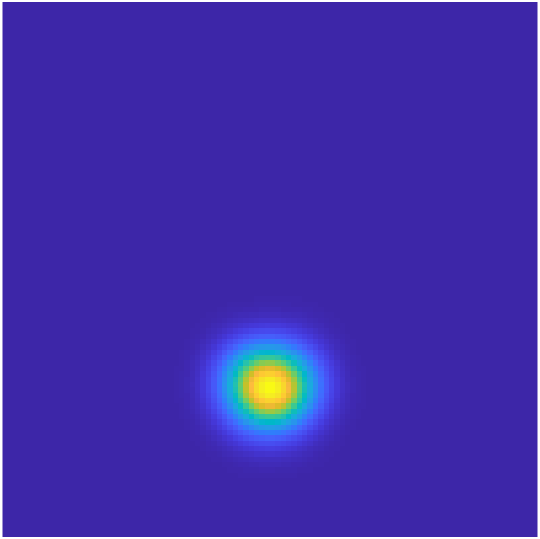}\\
\centering {$F_{10*}(P_0)$}
\end{minipage}\hfill
\begin{minipage}{0.16\linewidth}
\includegraphics[width=1\linewidth]{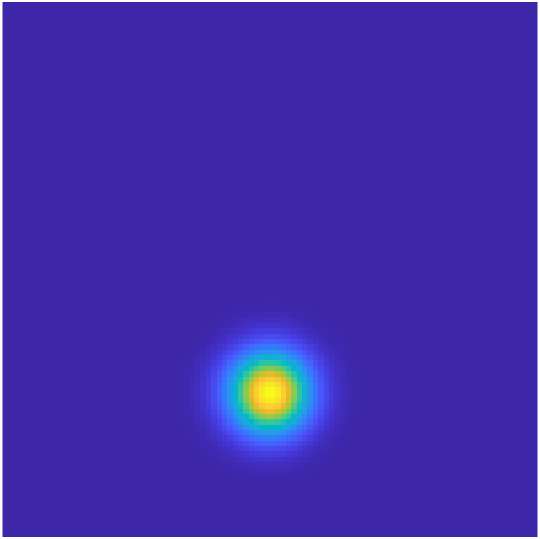}\\
\centering {$P_1$}
\end{minipage}\hfill

\caption{Evolution of the density for strong ($\lambda_{\mathcal{I}}=1$) penalty on conflicts with the obstacle.}
\label{crowd_10D_strong_evo}
\end{figure}







\section{NF Experiments}

\subsection{More Synthetic Data}

\begin{figure}[H]
\centering
\begin{minipage}{0.142\linewidth}
\includegraphics[width=1\linewidth]{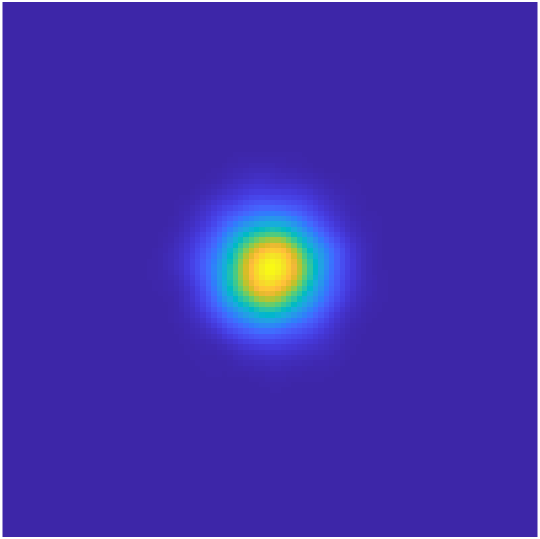}\\
\end{minipage}\hfill
\begin{minipage}{0.142\linewidth}
\includegraphics[width=1\linewidth]{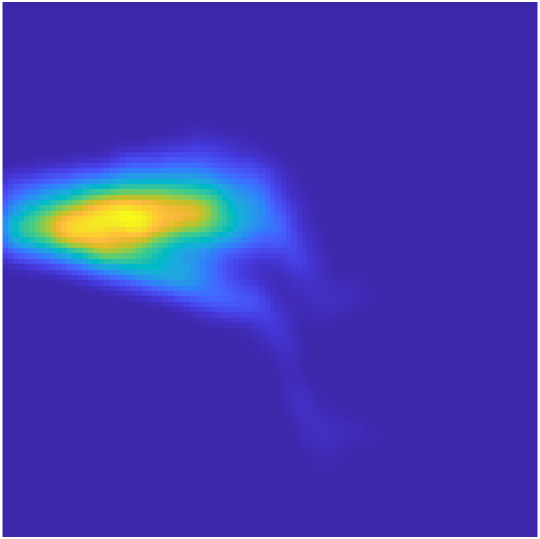}\\
\end{minipage}\hfill
\begin{minipage}{0.142\linewidth}
\includegraphics[width=1\linewidth]{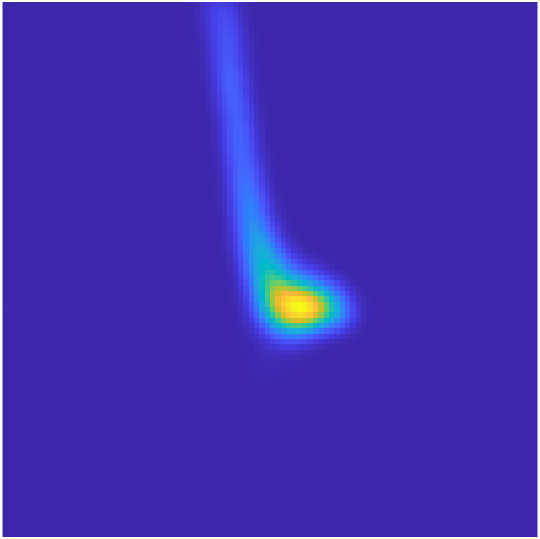}\\
\end{minipage}\hfill
\begin{minipage}{0.142\linewidth}
\includegraphics[width=1\linewidth]{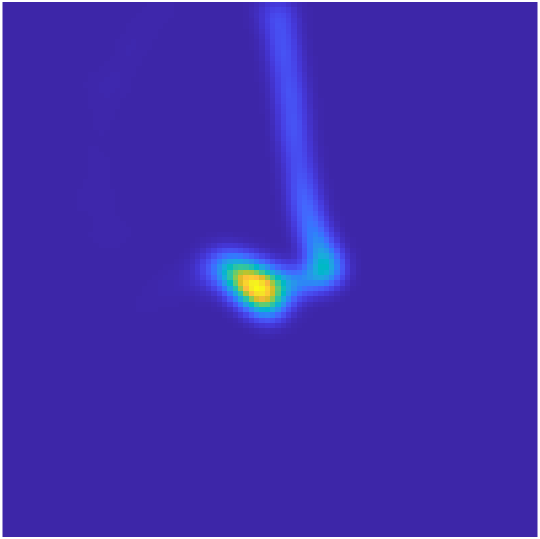}\\
\end{minipage}\hfill
\begin{minipage}{0.142\linewidth}
\includegraphics[width=1\linewidth]{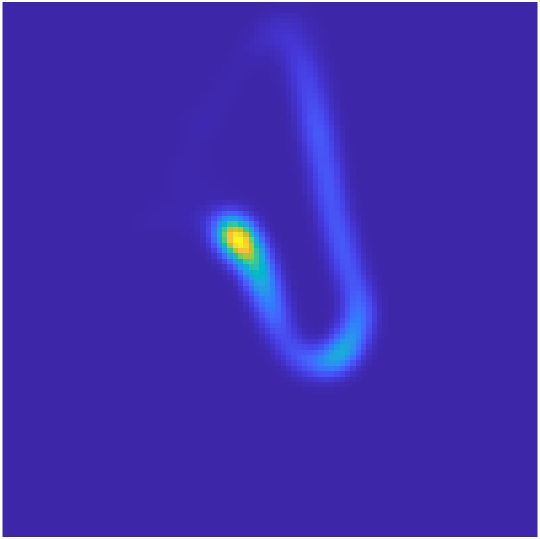}\\
\end{minipage}\hfill
\begin{minipage}{0.142\linewidth}
\includegraphics[width=1\linewidth]{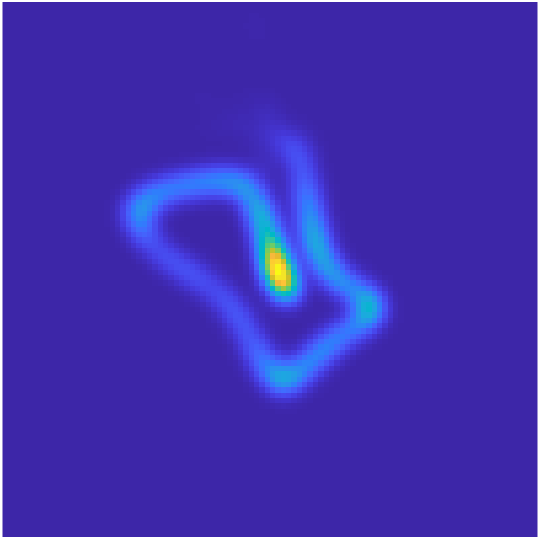}\\
\end{minipage}\hfill
\begin{minipage}{0.142\linewidth}
\includegraphics[width=1\linewidth]{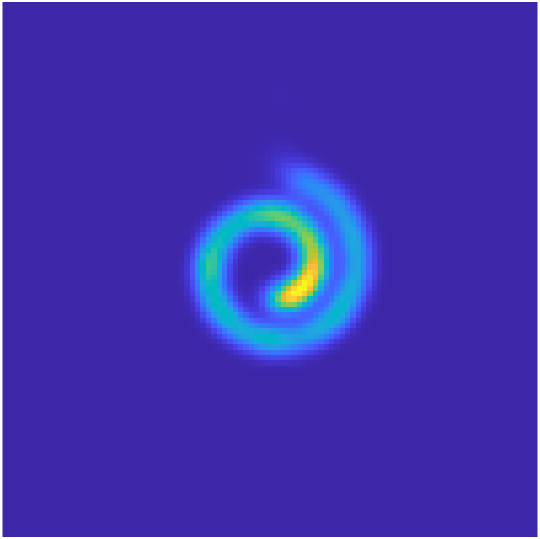}\\
\end{minipage}\hfill

\begin{minipage}{0.142\linewidth}
\includegraphics[width=1\linewidth]{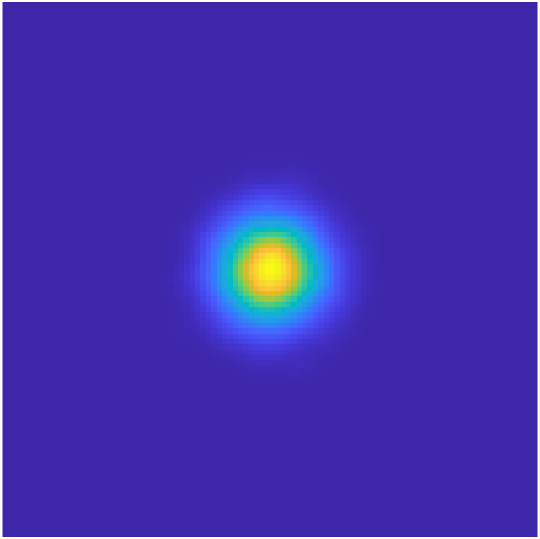}\\
\end{minipage}\hfill
\begin{minipage}{0.142\linewidth}
\includegraphics[width=1\linewidth]{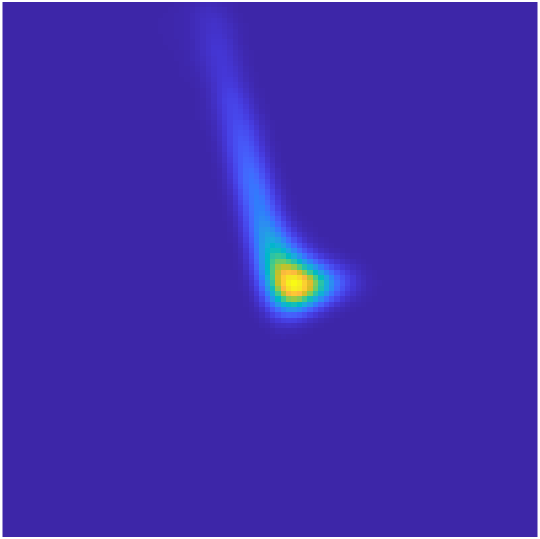}\\
\end{minipage}\hfill
\begin{minipage}{0.142\linewidth}
\includegraphics[width=1\linewidth]{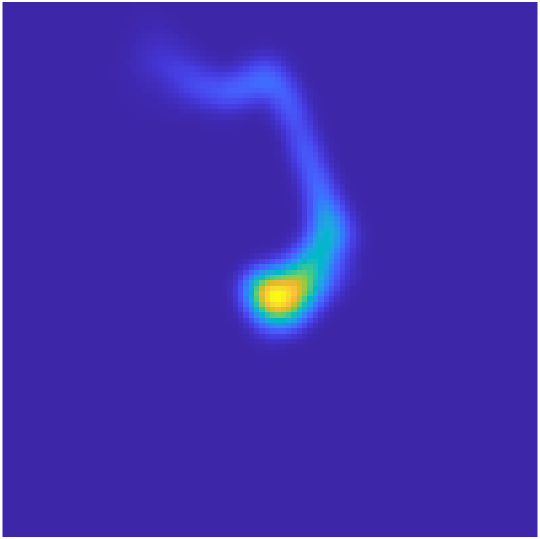}\\
\end{minipage}\hfill
\begin{minipage}{0.142\linewidth}
\includegraphics[width=1\linewidth]{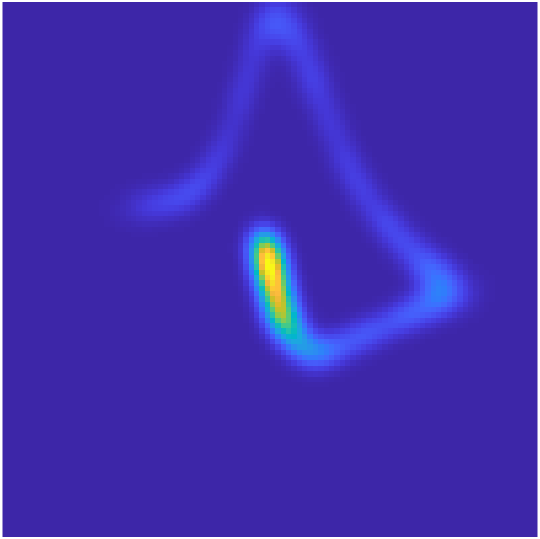}\\
\end{minipage}\hfill
\begin{minipage}{0.142\linewidth}
\includegraphics[width=1\linewidth]{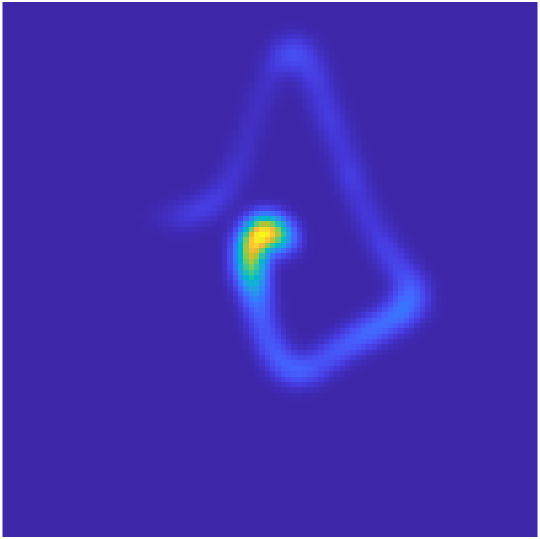}\\
\end{minipage}\hfill
\begin{minipage}{0.142\linewidth}
\includegraphics[width=1\linewidth]{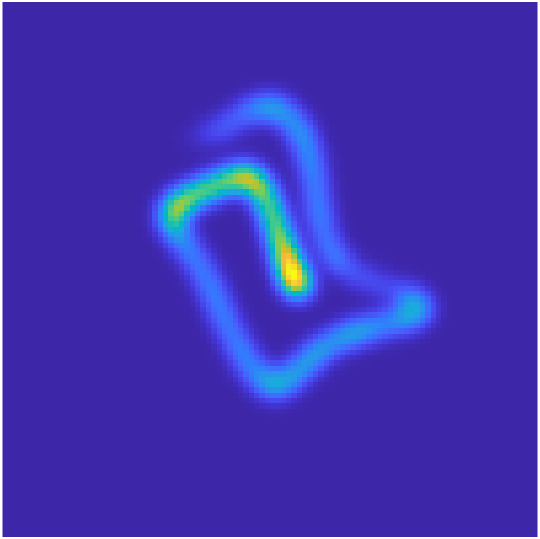}\\
\end{minipage}\hfill
\begin{minipage}{0.142\linewidth}
\includegraphics[width=1\linewidth]{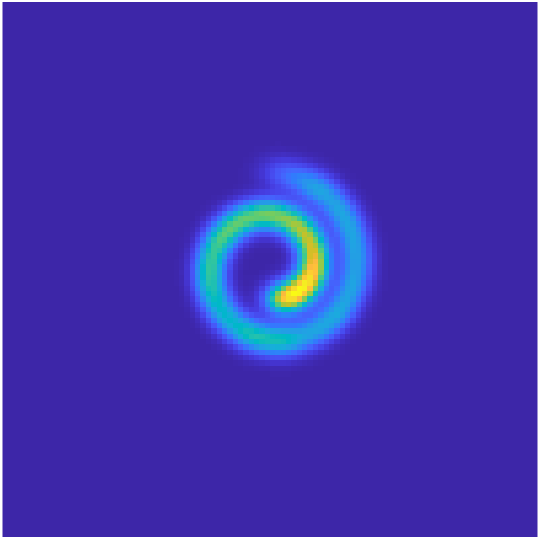}\\
\end{minipage}\hfill

\begin{minipage}{0.49\linewidth}
\includegraphics[width=1\linewidth]{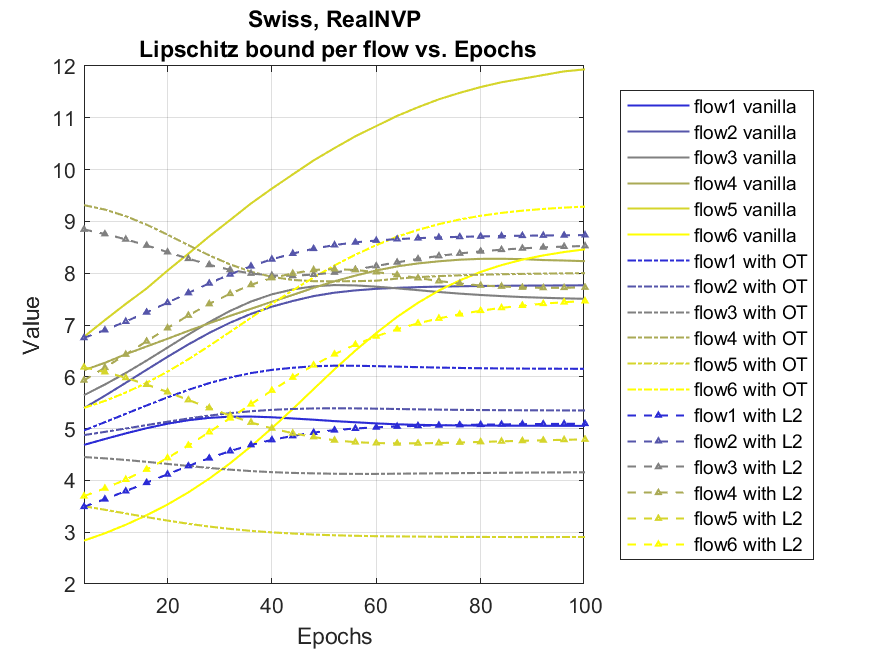}\\
\end{minipage}\hfill
\begin{minipage}{0.49\linewidth}
\includegraphics[width=1\linewidth]{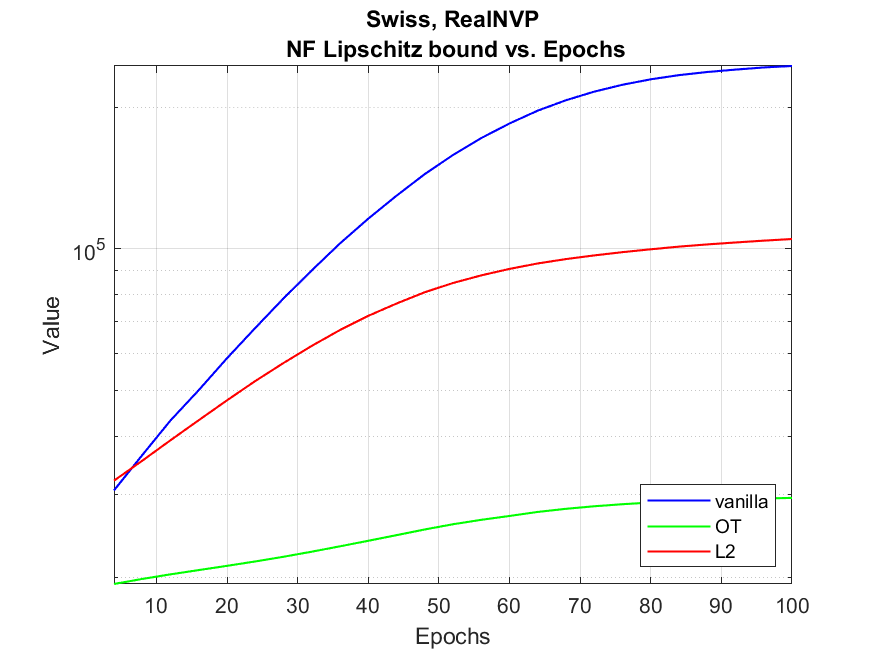}\\
\end{minipage}\hfill

\caption{Top row: output of each intermediate RealNVP flow between a single Gaussian and the Swiss roll-shape density, trained without using the transport cost. Middle row: same but trained with using the transport cost. Bottom left: the Lipschitz bound for each flow over the training epochs. Bottom right: the Lipschitz bound for the entire flow. The weights of the transport regularization are chosen so that the negative log-likelihood is not severely obstructed.}
\label{RNVP_Swiss}
\end{figure}

\begin{figure}[H]
\centering
\begin{minipage}{0.142\linewidth}
\includegraphics[width=1\linewidth]{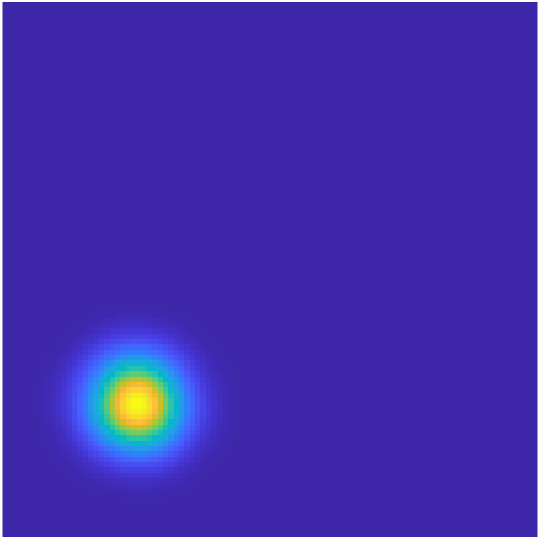}\\
\end{minipage}\hfill
\begin{minipage}{0.142\linewidth}
\includegraphics[width=1\linewidth]{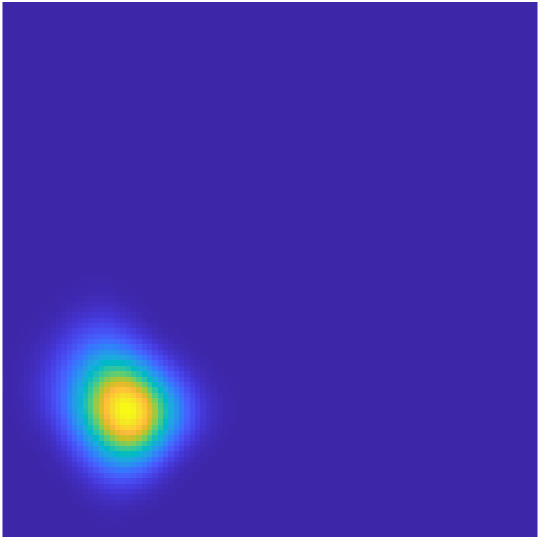}\\
\end{minipage}\hfill
\begin{minipage}{0.142\linewidth}
\includegraphics[width=1\linewidth]{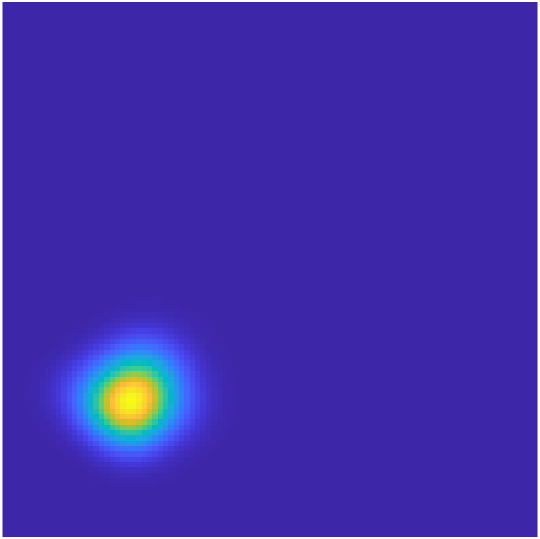}\\
\end{minipage}\hfill
\begin{minipage}{0.142\linewidth}
\includegraphics[width=1\linewidth]{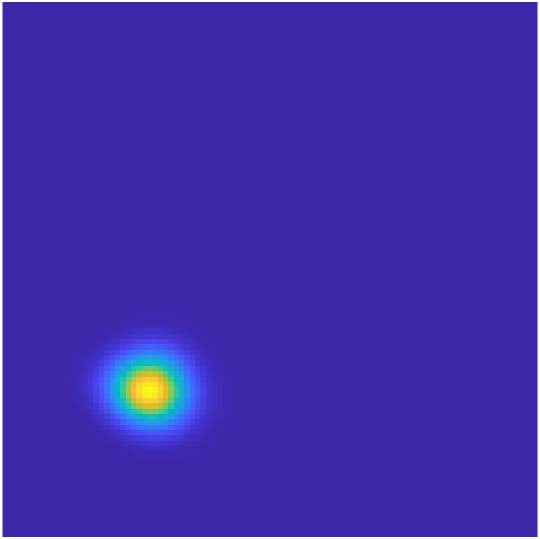}\\
\end{minipage}\hfill
\begin{minipage}{0.142\linewidth}
\includegraphics[width=1\linewidth]{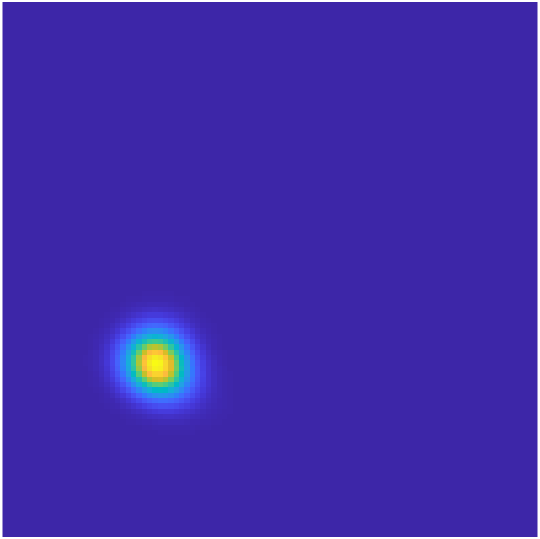}\\
\end{minipage}\hfill
\begin{minipage}{0.142\linewidth}
\includegraphics[width=1\linewidth]{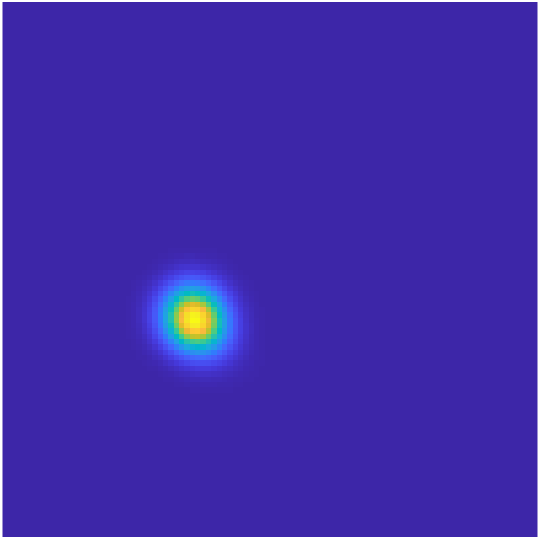}\\
\end{minipage}\hfill
\begin{minipage}{0.142\linewidth}
\includegraphics[width=1\linewidth]{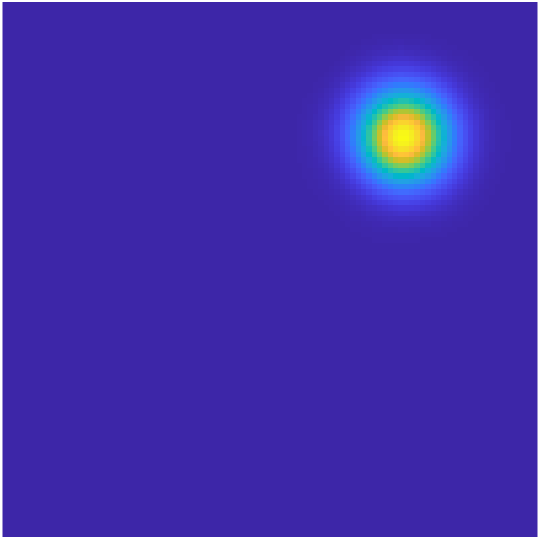}\\
\end{minipage}\hfill

\begin{minipage}{0.142\linewidth}
\includegraphics[width=1\linewidth]{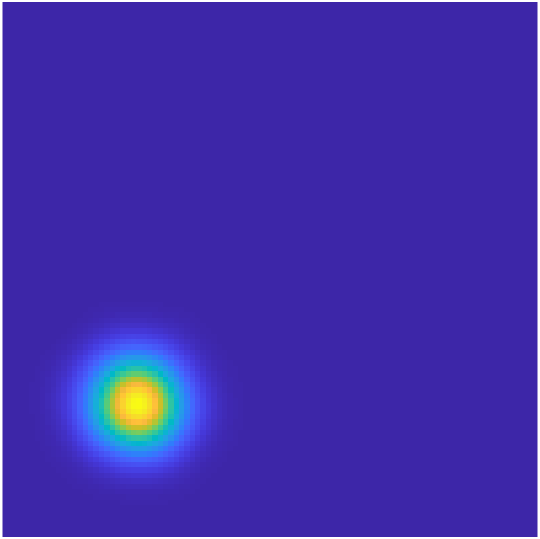}\\
\end{minipage}\hfill
\begin{minipage}{0.142\linewidth}
\includegraphics[width=1\linewidth]{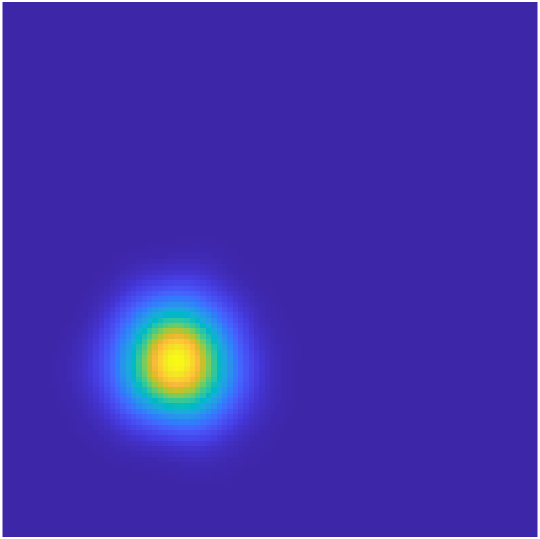}\\
\end{minipage}\hfill
\begin{minipage}{0.142\linewidth}
\includegraphics[width=1\linewidth]{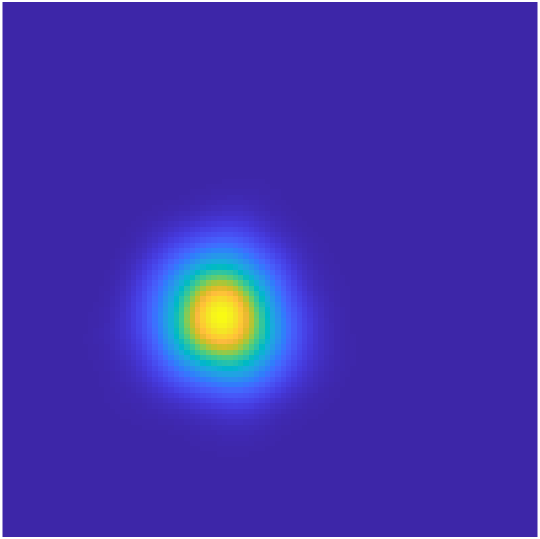}\\
\end{minipage}\hfill
\begin{minipage}{0.142\linewidth}
\includegraphics[width=1\linewidth]{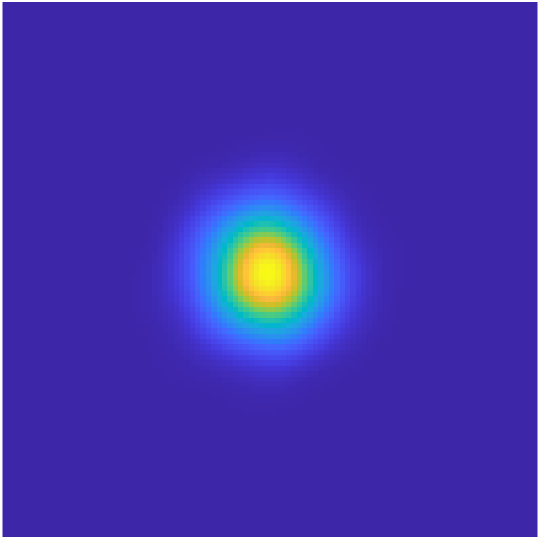}\\
\end{minipage}\hfill
\begin{minipage}{0.142\linewidth}
\includegraphics[width=1\linewidth]{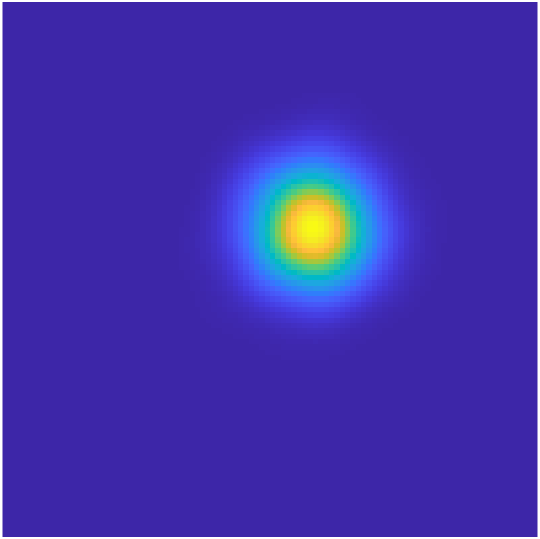}\\
\end{minipage}\hfill
\begin{minipage}{0.142\linewidth}
\includegraphics[width=1\linewidth]{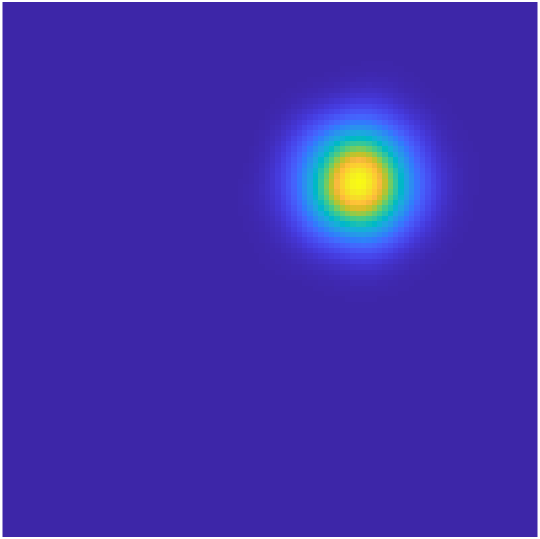}\\
\end{minipage}\hfill
\begin{minipage}{0.142\linewidth}
\includegraphics[width=1\linewidth]{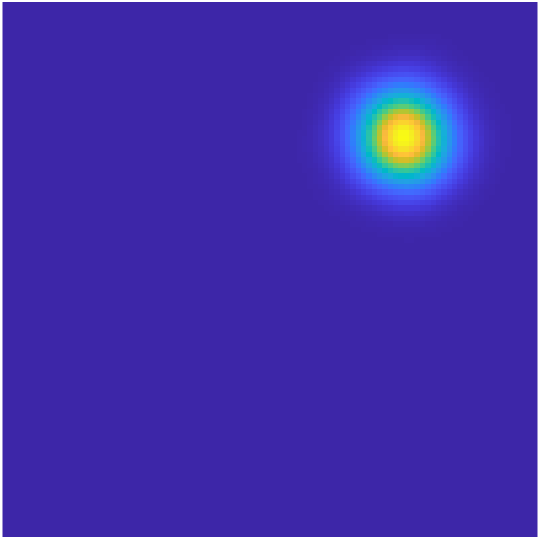}\\
\end{minipage}\hfill

\begin{minipage}{0.142\linewidth}
\includegraphics[width=1\linewidth]{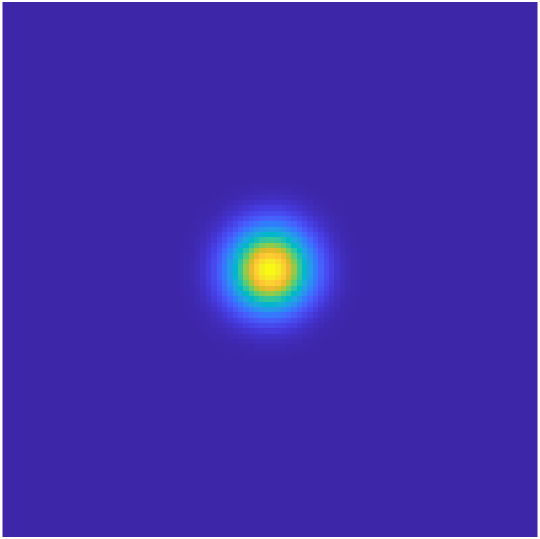}\\
\end{minipage}\hfill
\begin{minipage}{0.142\linewidth}
\includegraphics[width=1\linewidth]{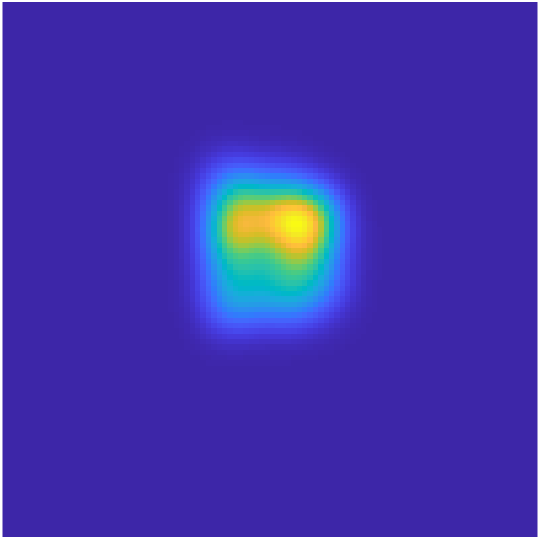}\\
\end{minipage}\hfill
\begin{minipage}{0.142\linewidth}
\includegraphics[width=1\linewidth]{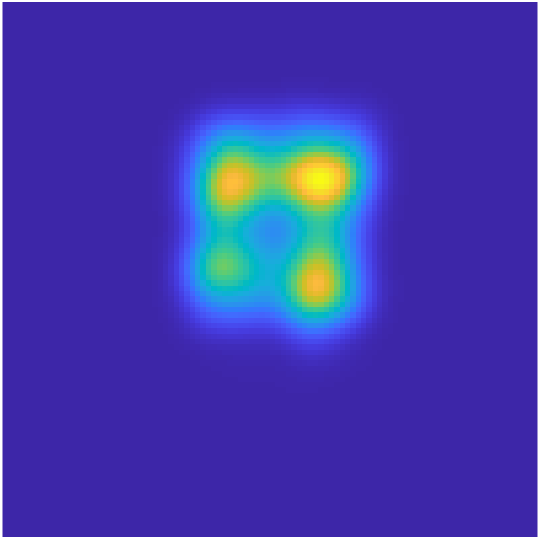}\\
\end{minipage}\hfill
\begin{minipage}{0.142\linewidth}
\includegraphics[width=1\linewidth]{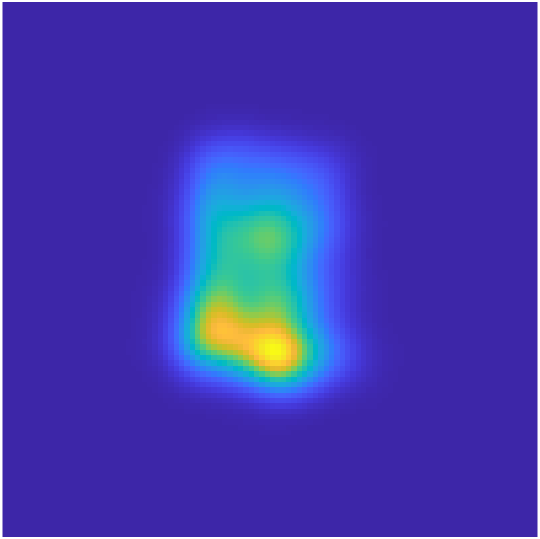}\\
\end{minipage}\hfill
\begin{minipage}{0.142\linewidth}
\includegraphics[width=1\linewidth]{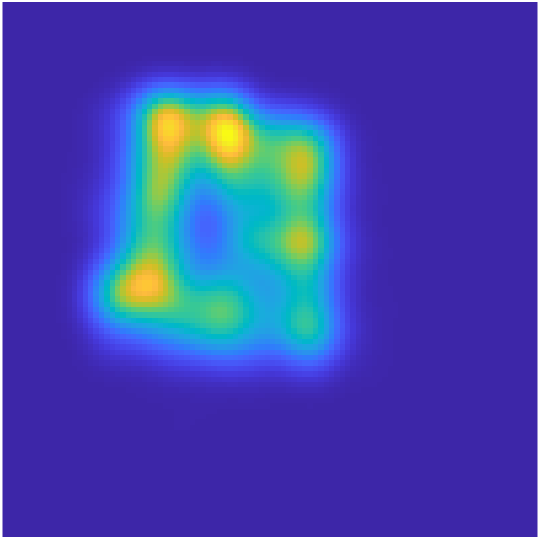}\\
\end{minipage}\hfill
\begin{minipage}{0.142\linewidth}
\includegraphics[width=1\linewidth]{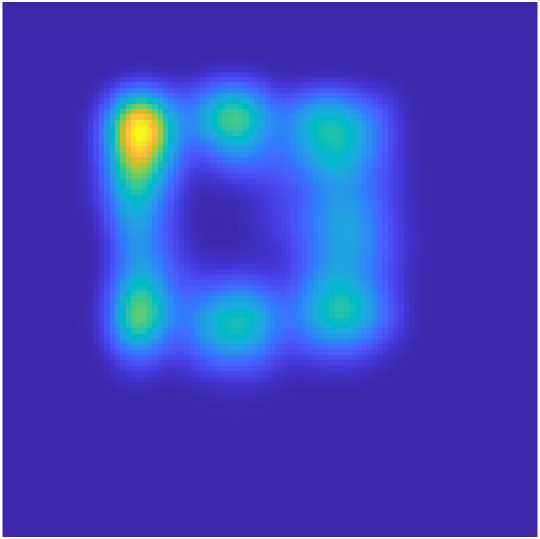}\\
\end{minipage}\hfill
\begin{minipage}{0.142\linewidth}
\includegraphics[width=1\linewidth]{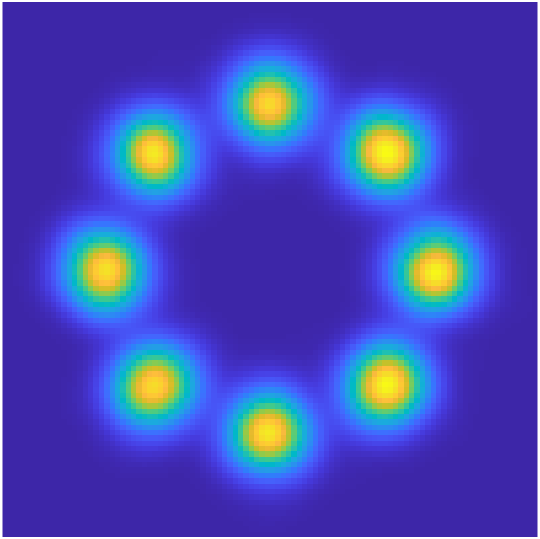}\\
\end{minipage}\hfill

\begin{minipage}{0.142\linewidth}
\includegraphics[width=1\linewidth]{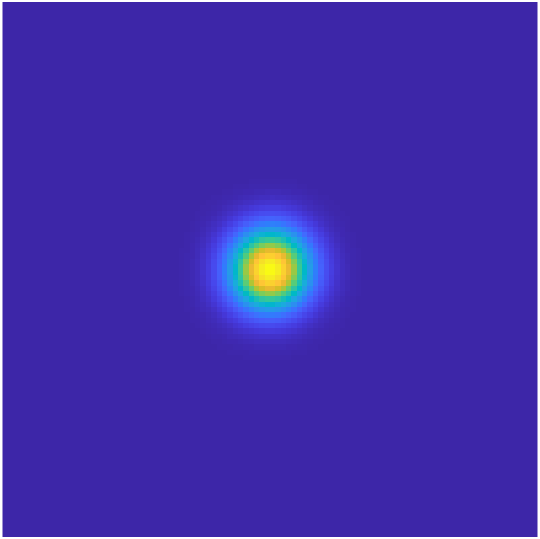}\\
\end{minipage}\hfill
\begin{minipage}{0.142\linewidth}
\includegraphics[width=1\linewidth]{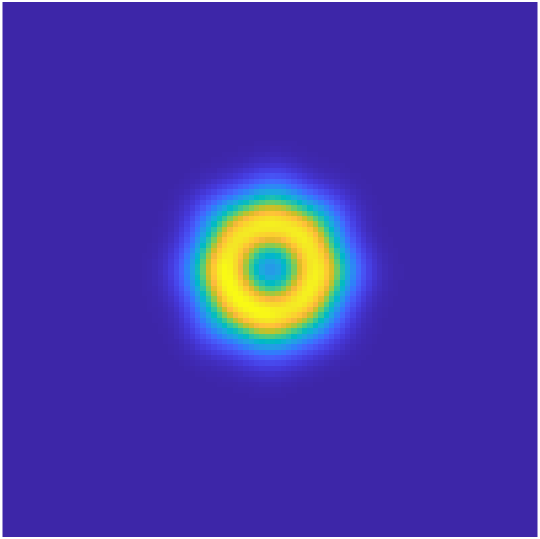}\\
\end{minipage}\hfill
\begin{minipage}{0.142\linewidth}
\includegraphics[width=1\linewidth]{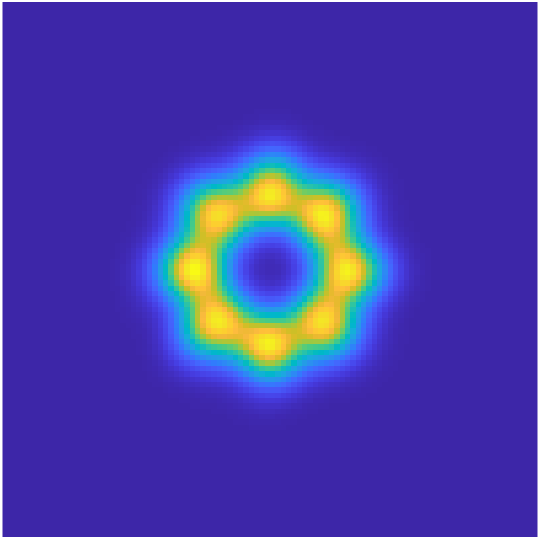}\\
\end{minipage}\hfill
\begin{minipage}{0.142\linewidth}
\includegraphics[width=1\linewidth]{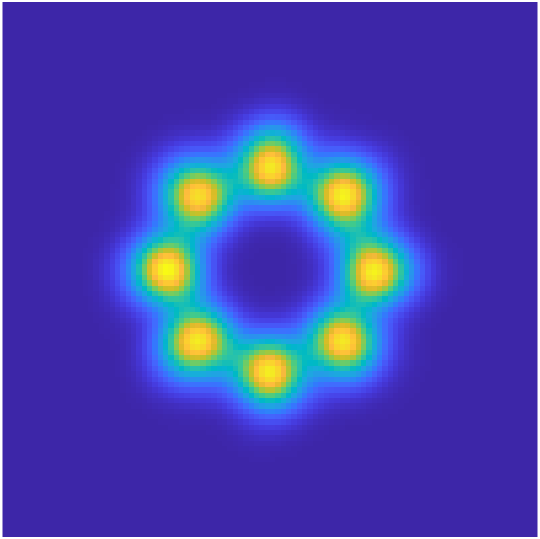}\\
\end{minipage}\hfill
\begin{minipage}{0.142\linewidth}
\includegraphics[width=1\linewidth]{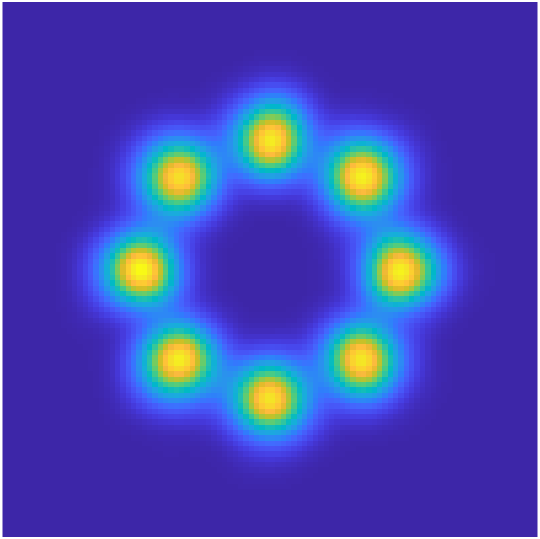}\\
\end{minipage}\hfill
\begin{minipage}{0.142\linewidth}
\includegraphics[width=1\linewidth]{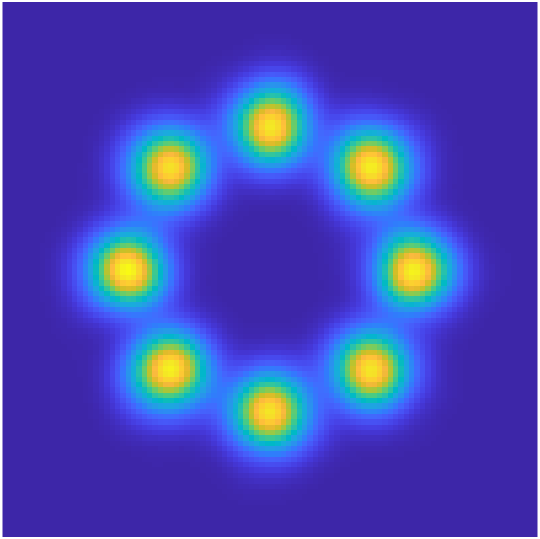}\\
\end{minipage}\hfill
\begin{minipage}{0.142\linewidth}
\includegraphics[width=1\linewidth]{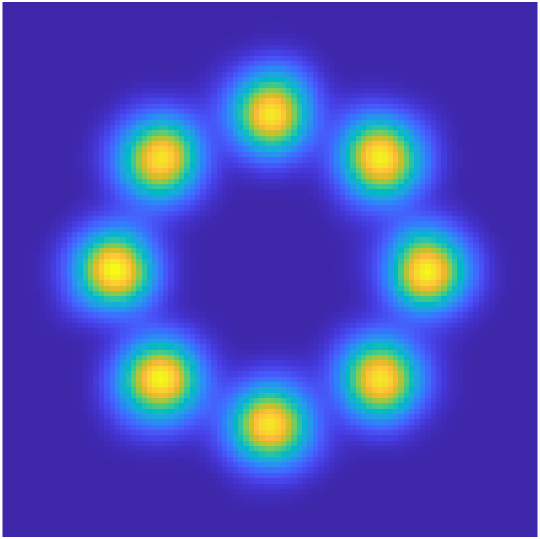}\\
\end{minipage}\hfill

\caption{Rows 1 and 2: output of each intermediate RealNVP flow between two Gaussians, trained without and with using the transport cost, respectively. Rows 3 and 4: output of each intermediate NSF-CL flow between a Gaussian and a Gaussian mixture, trained without and with using the transport cost, respectively.}
\label{syn_fig_RealNVP}
\end{figure}

\begin{figure}[H]
\centering
\begin{minipage}{0.142\linewidth}
\includegraphics[width=1\linewidth]{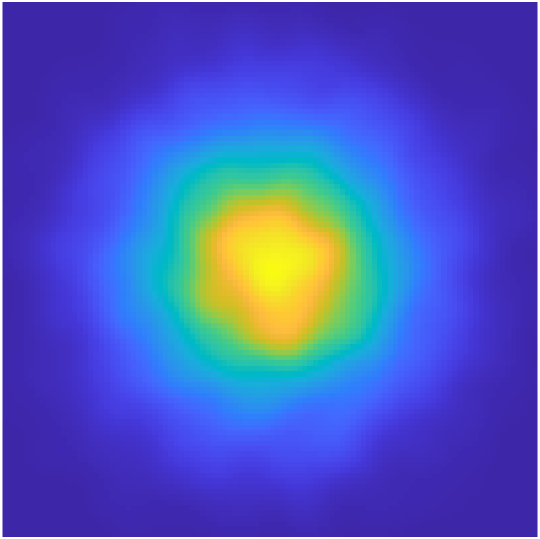}\\
\end{minipage}\hfill
\begin{minipage}{0.142\linewidth}
\includegraphics[width=1\linewidth]{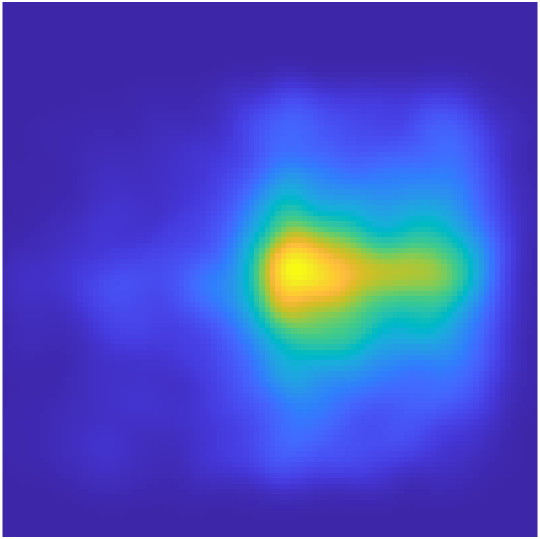}\\
\end{minipage}\hfill
\begin{minipage}{0.142\linewidth}
\includegraphics[width=1\linewidth]{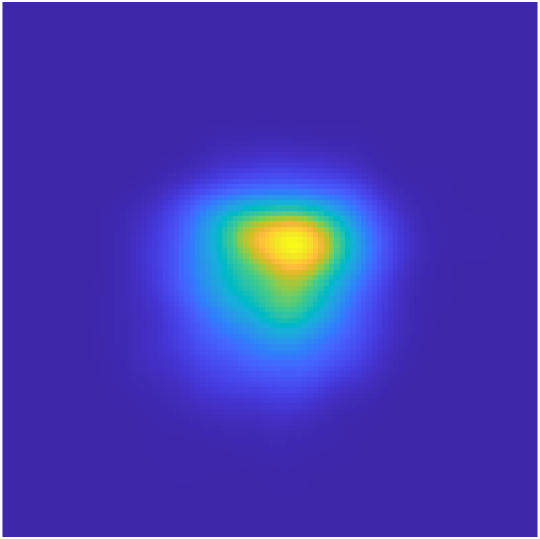}\\
\end{minipage}\hfill
\begin{minipage}{0.142\linewidth}
\includegraphics[width=1\linewidth]{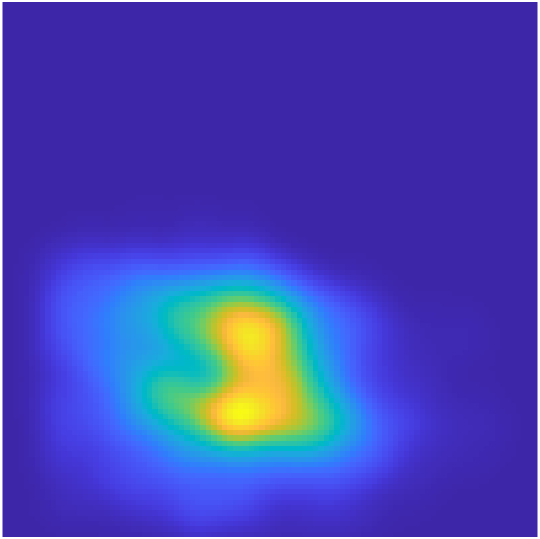}\\
\end{minipage}\hfill
\begin{minipage}{0.142\linewidth}
\includegraphics[width=1\linewidth]{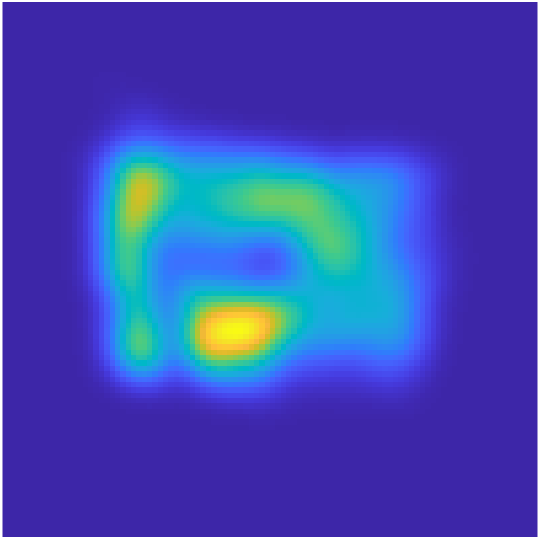}\\
\end{minipage}\hfill
\begin{minipage}{0.142\linewidth}
\includegraphics[width=1\linewidth]{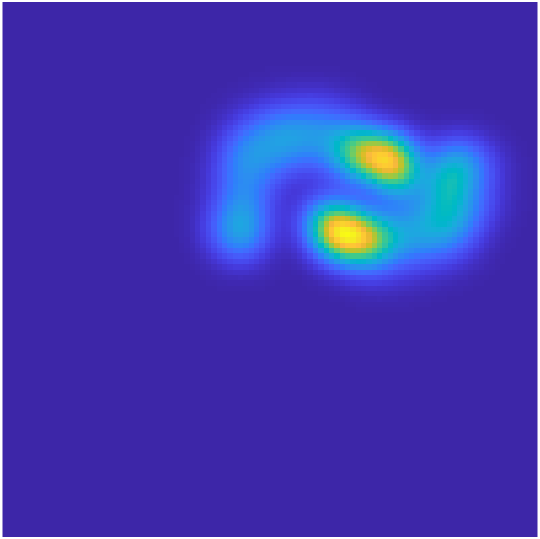}\\
\end{minipage}\hfill
\begin{minipage}{0.142\linewidth}
\includegraphics[width=1\linewidth]{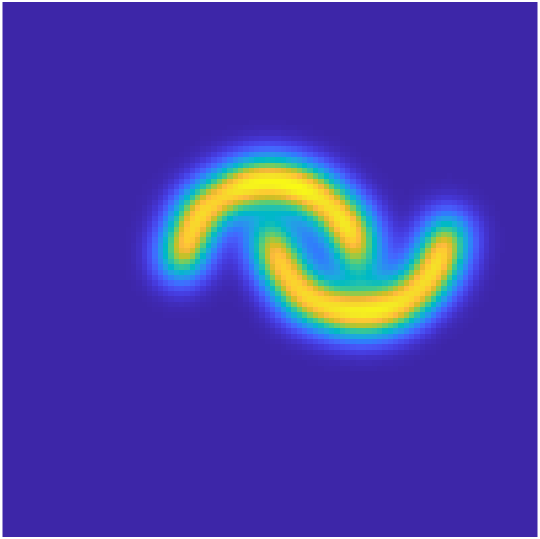}\\
\end{minipage}\hfill

\begin{minipage}{0.142\linewidth}
\includegraphics[width=1\linewidth]{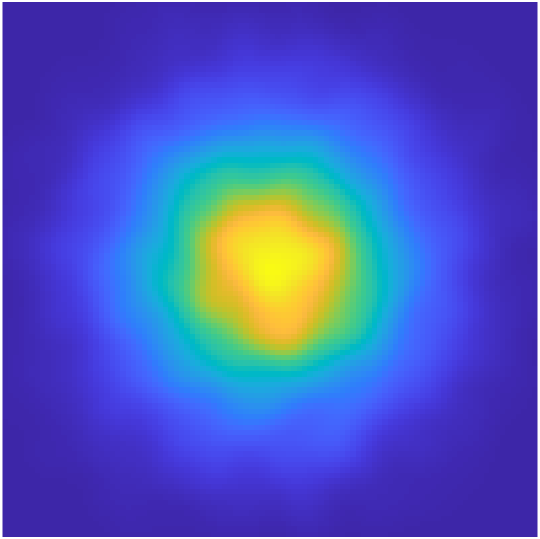}\\
\end{minipage}\hfill
\begin{minipage}{0.142\linewidth}
\includegraphics[width=1\linewidth]{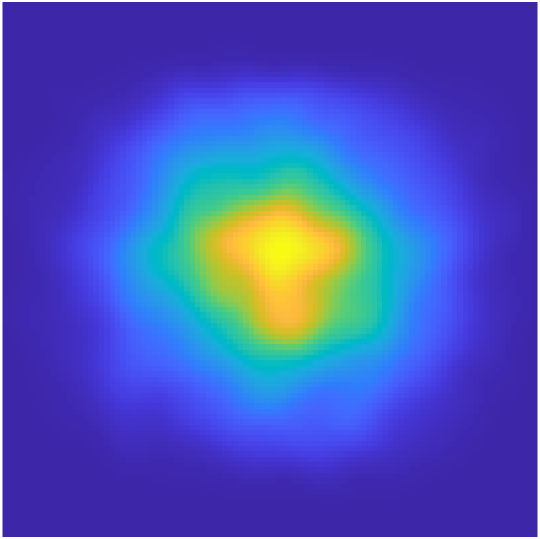}\\
\end{minipage}\hfill
\begin{minipage}{0.142\linewidth}
\includegraphics[width=1\linewidth]{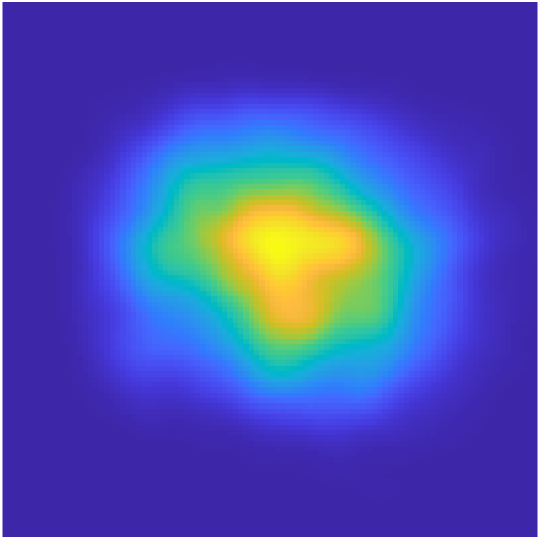}\\
\end{minipage}\hfill
\begin{minipage}{0.142\linewidth}
\includegraphics[width=1\linewidth]{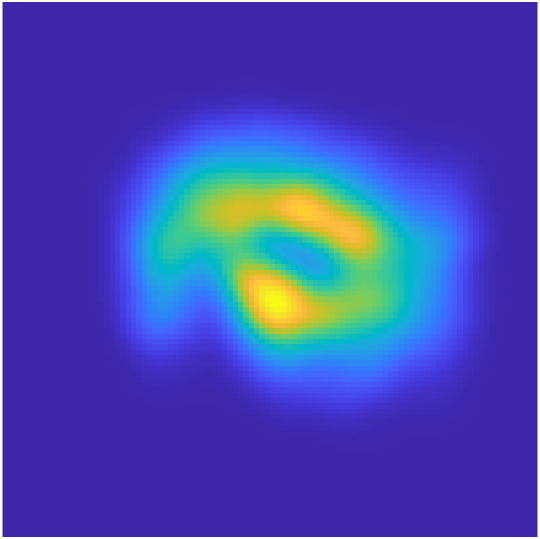}\\
\end{minipage}\hfill
\begin{minipage}{0.142\linewidth}
\includegraphics[width=1\linewidth]{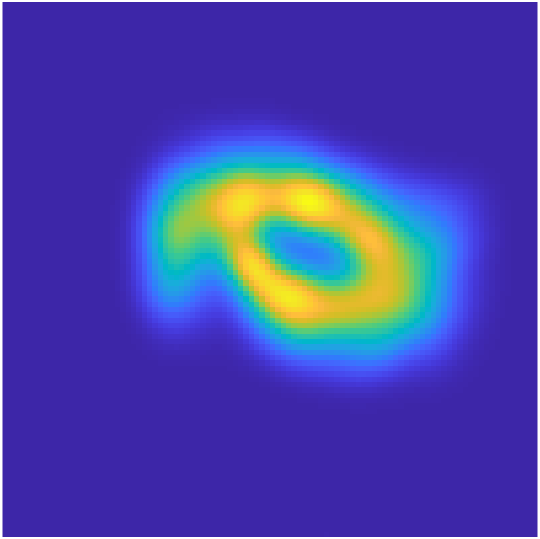}\\
\end{minipage}\hfill
\begin{minipage}{0.142\linewidth}
\includegraphics[width=1\linewidth]{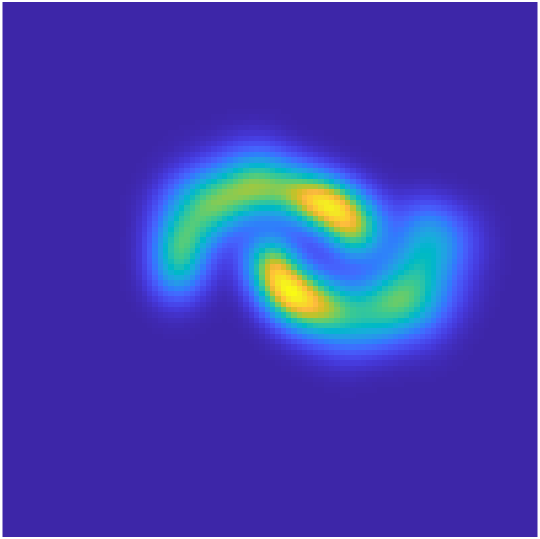}\\
\end{minipage}\hfill
\begin{minipage}{0.142\linewidth}
\includegraphics[width=1\linewidth]{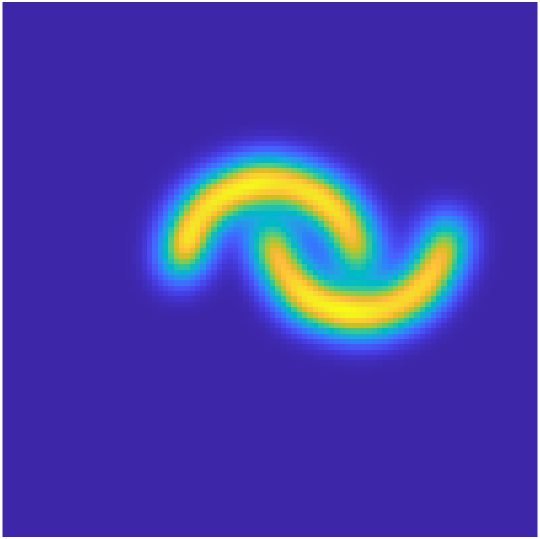}\\
\end{minipage}\hfill

\begin{minipage}{0.142\linewidth}
\includegraphics[width=1\linewidth]{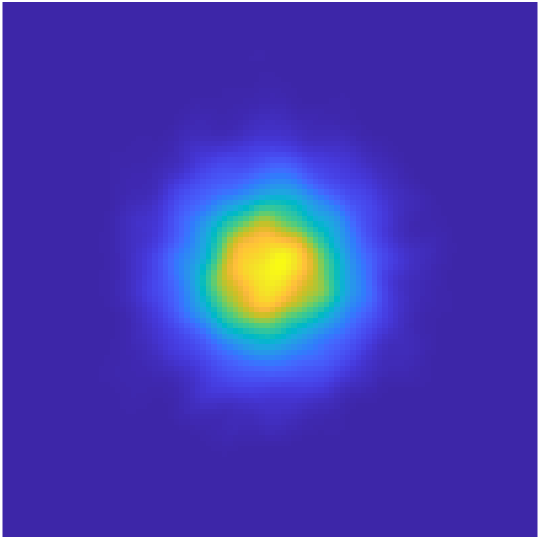}\\
\end{minipage}\hfill
\begin{minipage}{0.142\linewidth}
\includegraphics[width=1\linewidth]{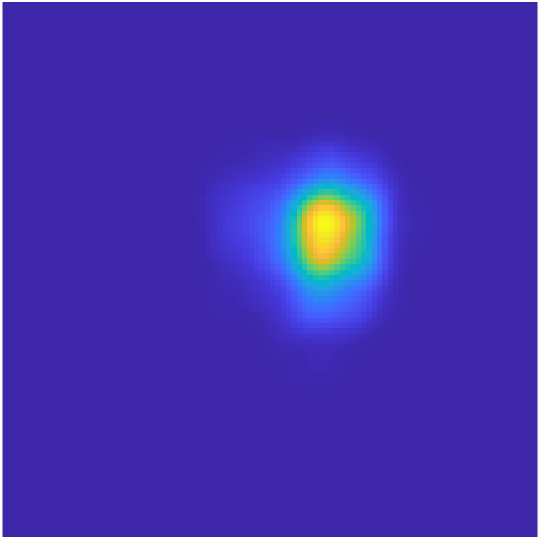}\\
\end{minipage}\hfill
\begin{minipage}{0.142\linewidth}
\includegraphics[width=1\linewidth]{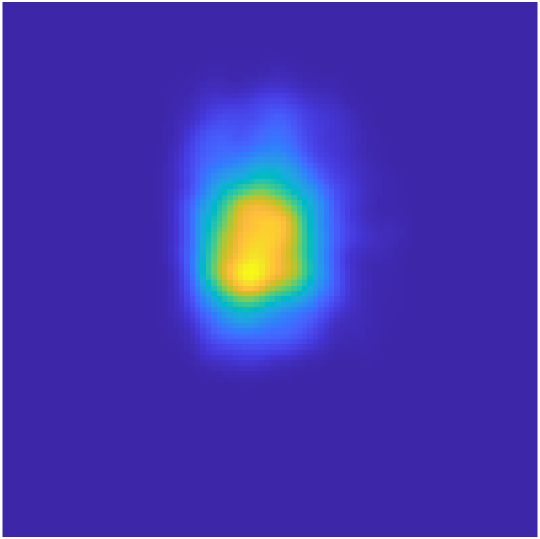}\\
\end{minipage}\hfill
\begin{minipage}{0.142\linewidth}
\includegraphics[width=1\linewidth]{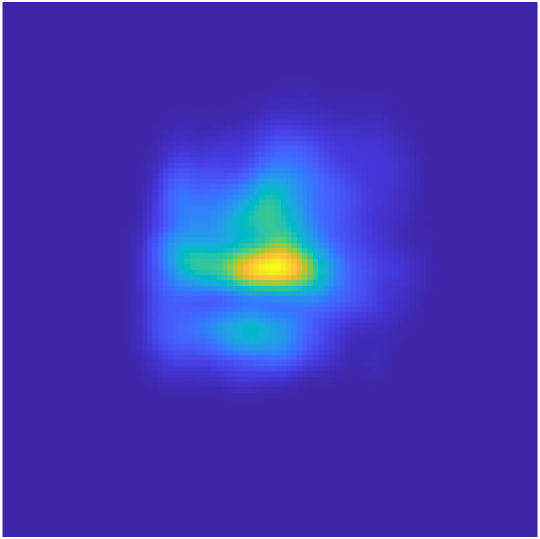}\\
\end{minipage}\hfill
\begin{minipage}{0.142\linewidth}
\includegraphics[width=1\linewidth]{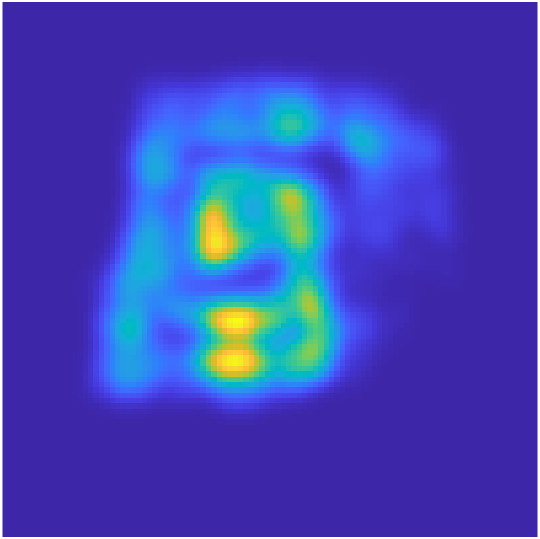}\\
\end{minipage}\hfill
\begin{minipage}{0.142\linewidth}
\includegraphics[width=1\linewidth]{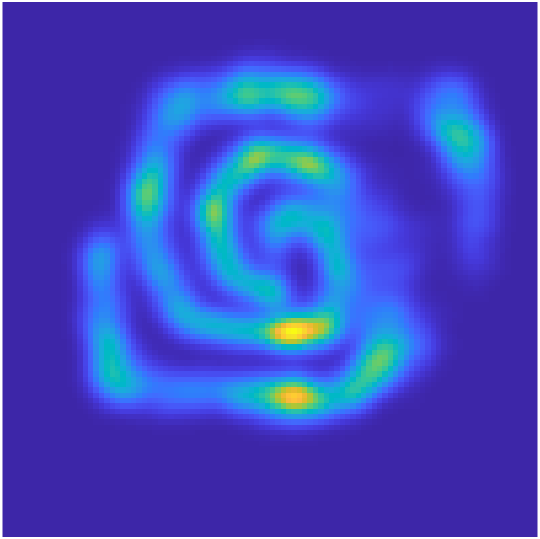}\\
\end{minipage}\hfill
\begin{minipage}{0.142\linewidth}
\includegraphics[width=1\linewidth]{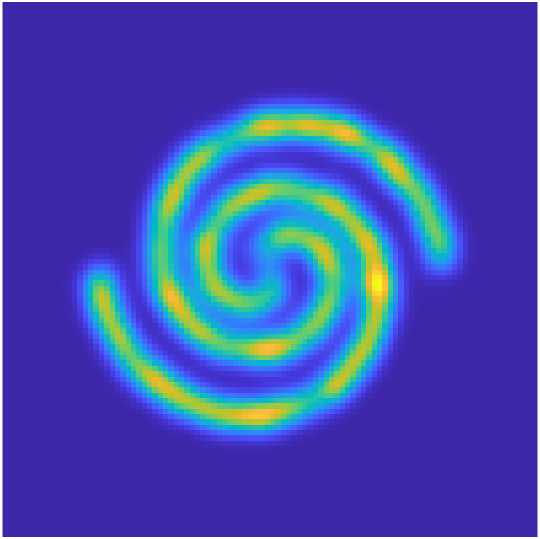}\\
\end{minipage}\hfill

\begin{minipage}{0.142\linewidth}
\includegraphics[width=1\linewidth]{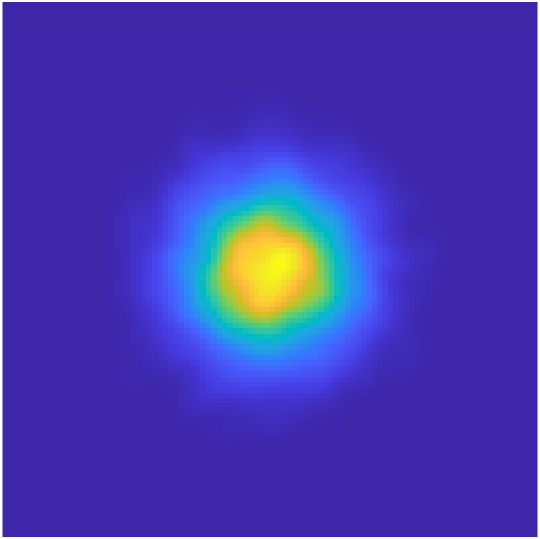}\\
\end{minipage}\hfill
\begin{minipage}{0.142\linewidth}
\includegraphics[width=1\linewidth]{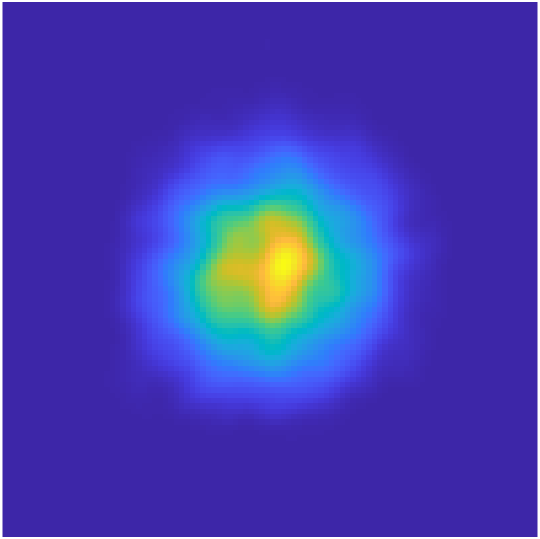}\\
\end{minipage}\hfill
\begin{minipage}{0.142\linewidth}
\includegraphics[width=1\linewidth]{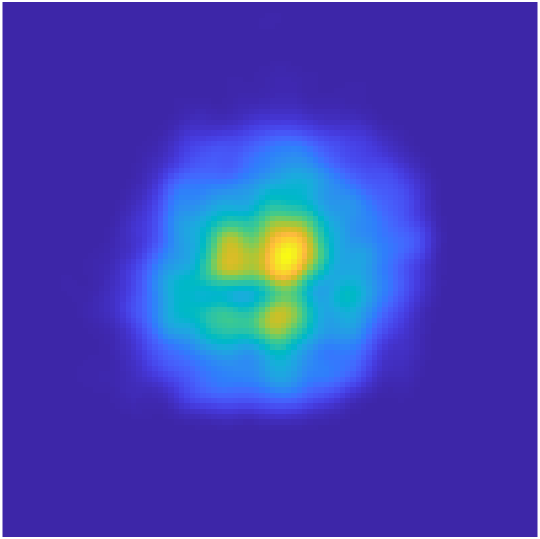}\\
\end{minipage}\hfill
\begin{minipage}{0.142\linewidth}
\includegraphics[width=1\linewidth]{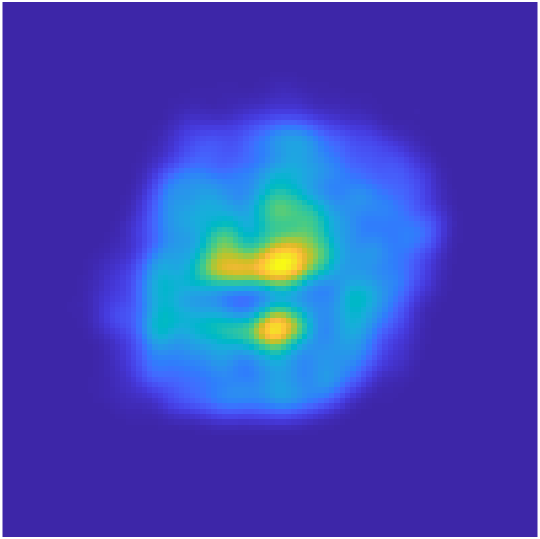}\\
\end{minipage}\hfill
\begin{minipage}{0.142\linewidth}
\includegraphics[width=1\linewidth]{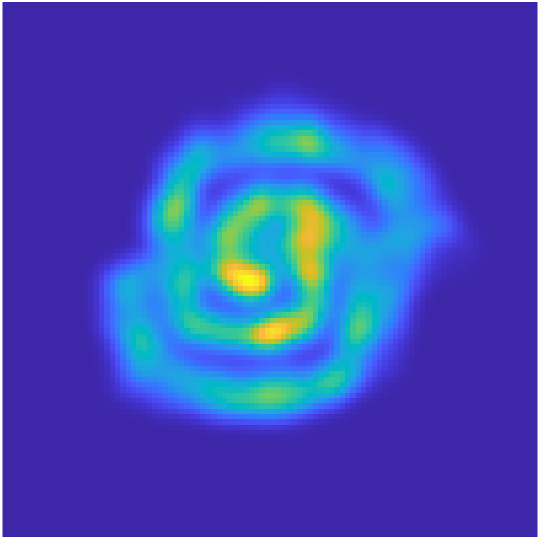}\\
\end{minipage}\hfill
\begin{minipage}{0.142\linewidth}
\includegraphics[width=1\linewidth]{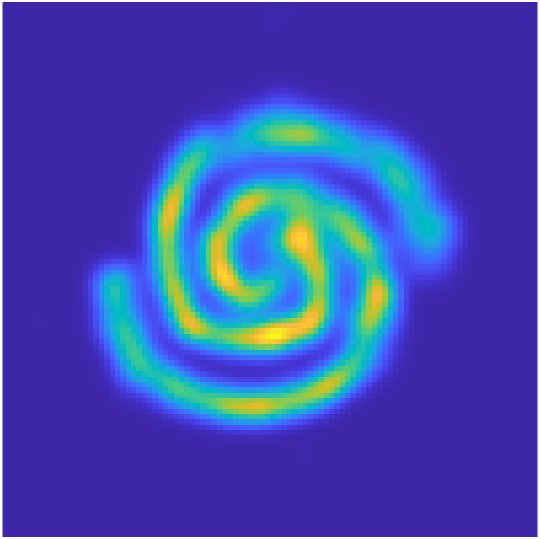}\\
\end{minipage}\hfill
\begin{minipage}{0.142\linewidth}
\includegraphics[width=1\linewidth]{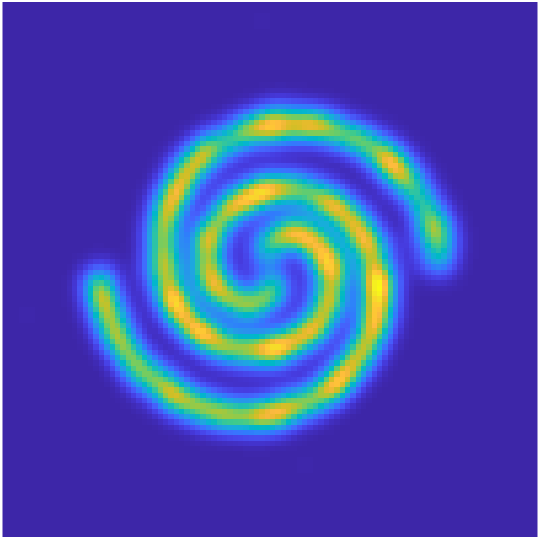}\\
\end{minipage}\hfill

\caption{Rows 1 and 2: output of each intermediate RealNVP flow between two Gaussians, trained without and with using the transport cost, respectively. Rows 3 and 4: output of each intermediate NSF-CL flow between a Gaussian and the spiral density, trained without and with using the transport cost, respectively.}
\label{syn_fig_NSF}
\end{figure}

\subsection{Tabular Datasets}

For the RealNVP flows, we used a replementation provided by~\cite{RNVP_github}, which achieved better log-likelihoods than the results reported in the original paper~\cite{RealNVP}. All RealNVP models use 6 flows with hidden dimensions of 256 in the $s,t$ networks. For NSF, we used the original implementation~\cite{NSF} with the reported hyperparameters in their paper. The likelihood results for the other models are taken from Table 3 of~\cite{NF_survey}. Our Lipschitz bound for a flow $f$ is defined as $\max_{\vx\in P_1} \|\nabla f(\vx)\|_2$, the estimated maximum gradient spectral norm over the training set $P_1$. The Lipschitz bound for the entire NF is the product of the bounds of the individual flows.

For a controlled study, all hyperparameters for the standard RealNVP (resp. NSF) and the RealNVP (resp. NSF) trained with transport costs are identical. The normalized weights $\frac{\lambda_L}{\lambda_{\mathcal{M}}}$ for the transport costs used to produce the results in Table~\ref{Tabular_table} are summarized in Table~\ref{Tabular_OT_weights}.

The dimension-wise permutations have shown to improve the flexibility of flow transforms both theoretically~\cite{universality_CL} and empirically~\cite{glow}. When transport costs are computed, we omit the costs of permutation transforms, since the underlying object stays the same.

\begin{figure}[H]
\centering
\begin{minipage}{0.45\linewidth}
\includegraphics[width=1\linewidth]{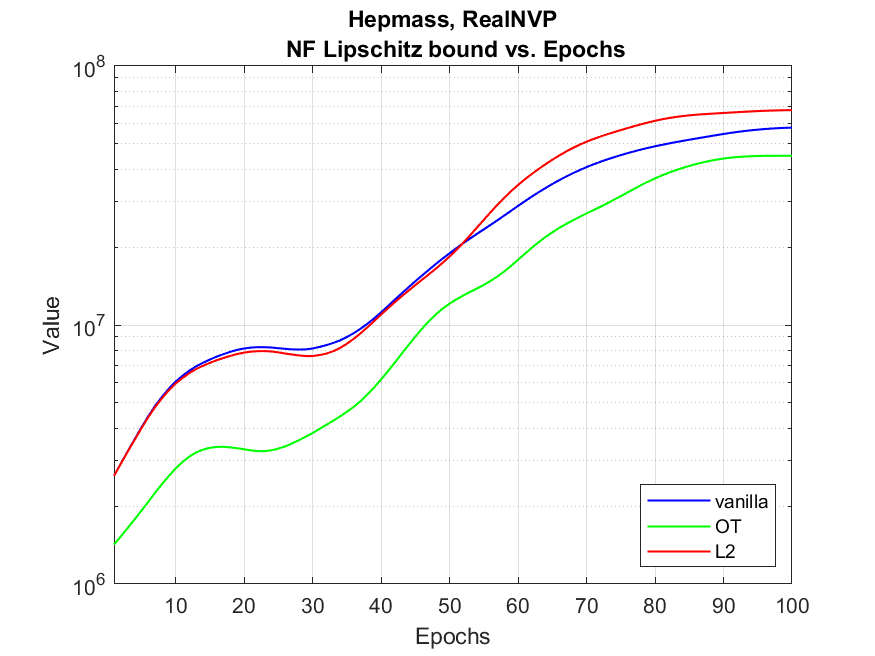}\\
\end{minipage}\hfill
\begin{minipage}{0.45\linewidth}
\includegraphics[width=1\linewidth]{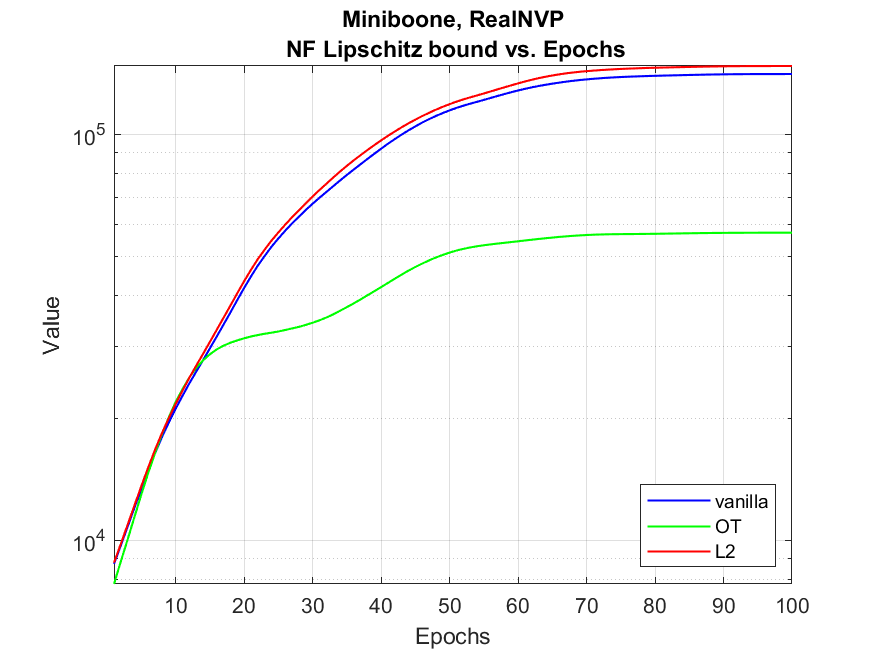}\\
\end{minipage}\hfill
%
\caption{The Lipschitz bound over trianing epochs, for different tabular datasets (left to right: Hepmass and Minibone).}
\label{RNVP_tabular_Lip_bound2}
\end{figure}

\begin{table}[H]
\caption{Normalized weights associated with the transport cost ($\frac{\lambda_L}{\lambda_{\mathcal{M}}}$) used to produce results in Table~\ref{Tabular_table}.}
\label{Tabular_OT_weights}
\begin{center}
\begin{small}
\begin{sc}
\begin{tabular}{lccccc}
\toprule
 & Power & Gas & Hepmass & Miniboone & BSDS300\\
\midrule
RNVP OT Weight & 5e-6 & 1e-6 & 5e-5 & 5e-3 & 1e-5
\\
NSF OT Weight & 1e-5 & 5e-5 & 5e-5 & 5e-3 & 5e-2\\
\bottomrule
\end{tabular}
\end{sc}
\end{small}
\end{center}
\vspace{-.2cm}
\end{table}

In the main text, we showed the marginal distributions on the first two dimensions for RealNVP trained on Miniboone. The remaining dimensions follow a similar trend. We additionally show those for dimensions 3,4 and for dimensions 5,6 here.

\begin{figure}[H]
\centering

\begin{minipage}{0.142\linewidth}
\includegraphics[width=1\linewidth]{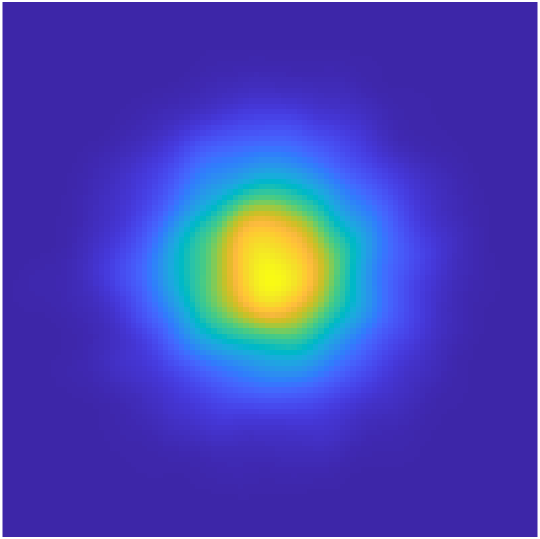}\\
\end{minipage}\hfill
\begin{minipage}{0.142\linewidth}
\includegraphics[width=1\linewidth]{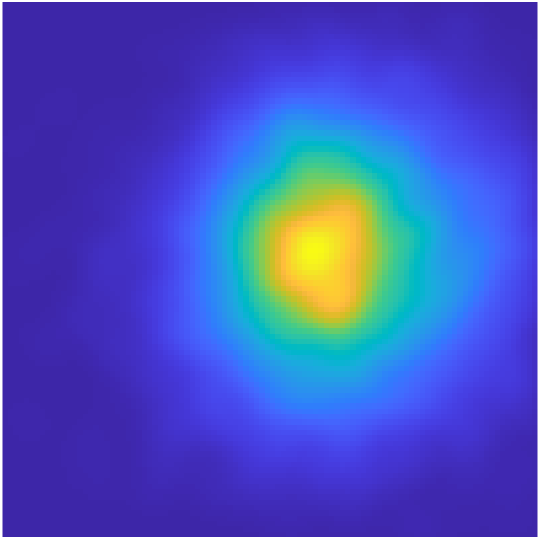}\\
\end{minipage}\hfill
\begin{minipage}{0.142\linewidth}
\includegraphics[width=1\linewidth]{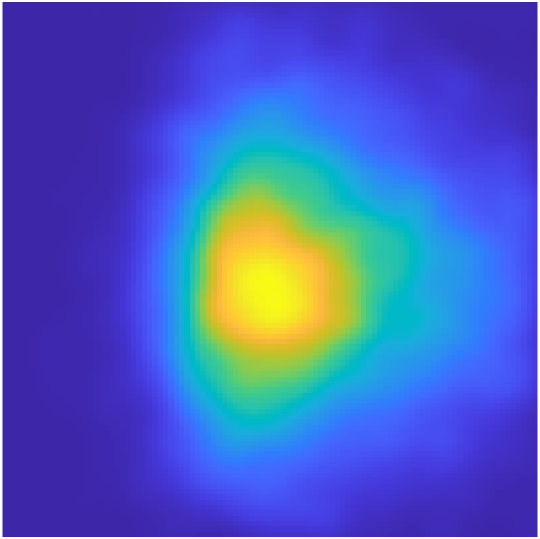}\\
\end{minipage}\hfill
\begin{minipage}{0.142\linewidth}
\includegraphics[width=1\linewidth]{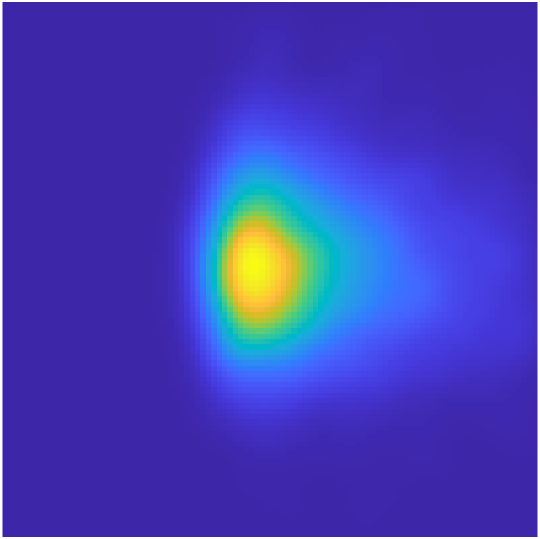}\\
\end{minipage}\hfill
\begin{minipage}{0.142\linewidth}
\includegraphics[width=1\linewidth]{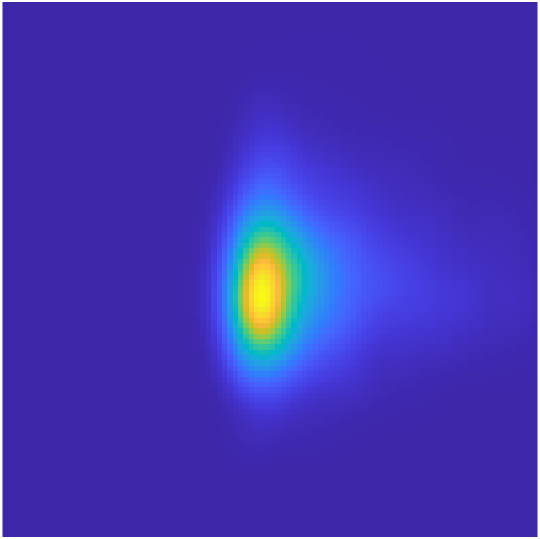}\\
\end{minipage}\hfill
\begin{minipage}{0.142\linewidth}
\includegraphics[width=1\linewidth]{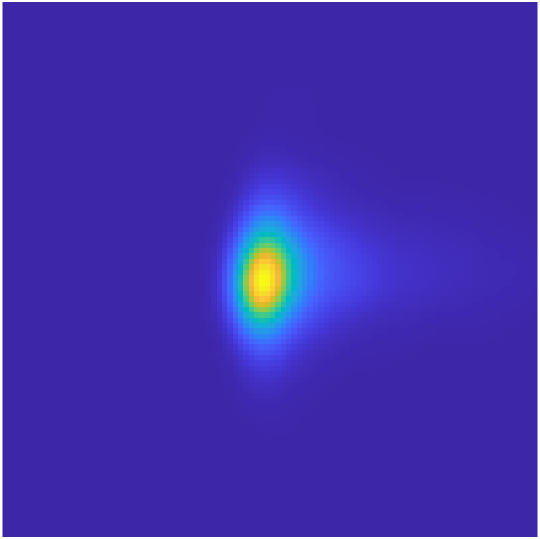}\\
\end{minipage}\hfill
\begin{minipage}{0.142\linewidth}
\includegraphics[width=1\linewidth]{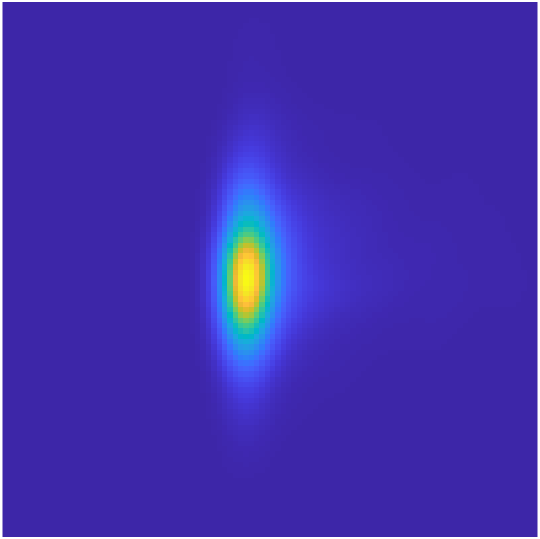}\\
\end{minipage}\hfill

\begin{minipage}{0.142\linewidth}
\includegraphics[width=1\linewidth]{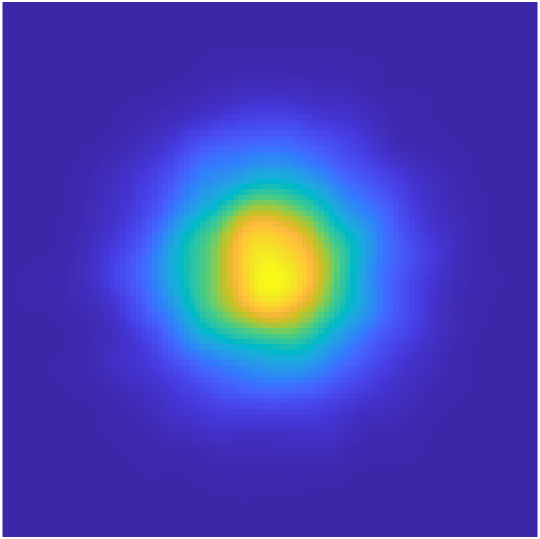}\\
\end{minipage}\hfill
\begin{minipage}{0.142\linewidth}
\includegraphics[width=1\linewidth]{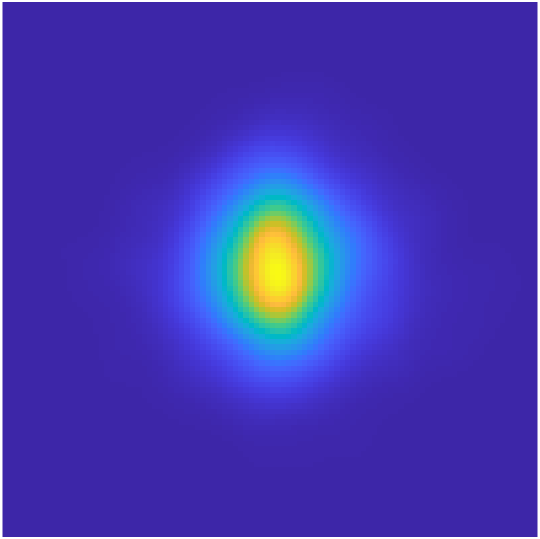}\\
\end{minipage}\hfill
\begin{minipage}{0.142\linewidth}
\includegraphics[width=1\linewidth]{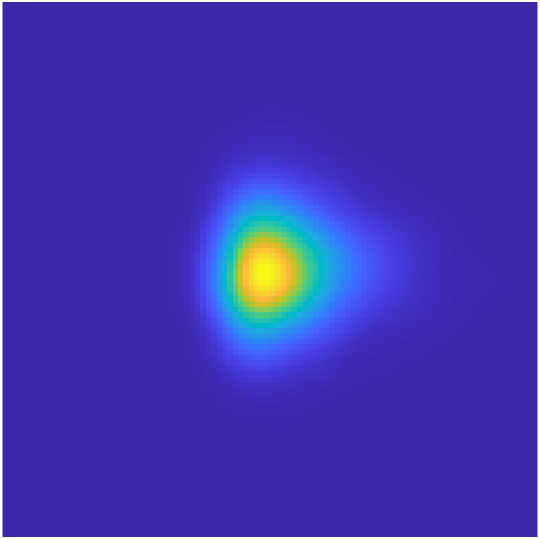}\\
\end{minipage}\hfill
\begin{minipage}{0.142\linewidth}
\includegraphics[width=1\linewidth]{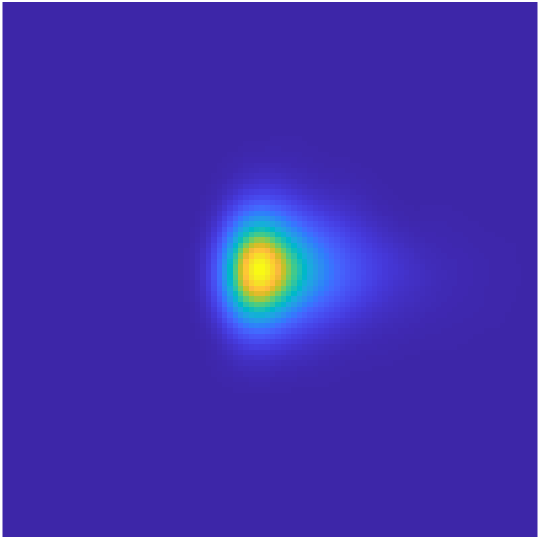}\\
\end{minipage}\hfill
\begin{minipage}{0.142\linewidth}
\includegraphics[width=1\linewidth]{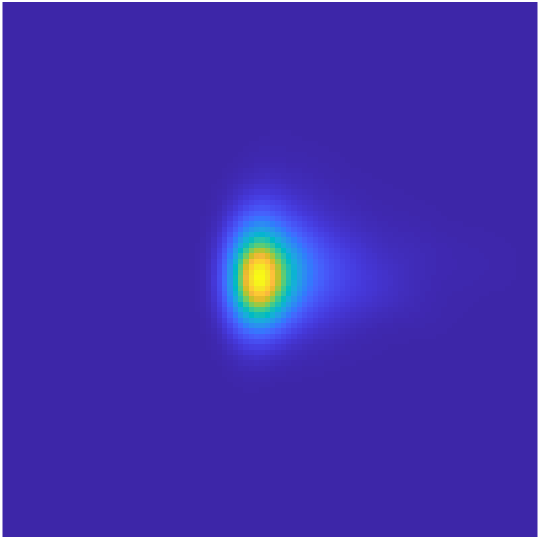}\\
\end{minipage}\hfill
\begin{minipage}{0.142\linewidth}
\includegraphics[width=1\linewidth]{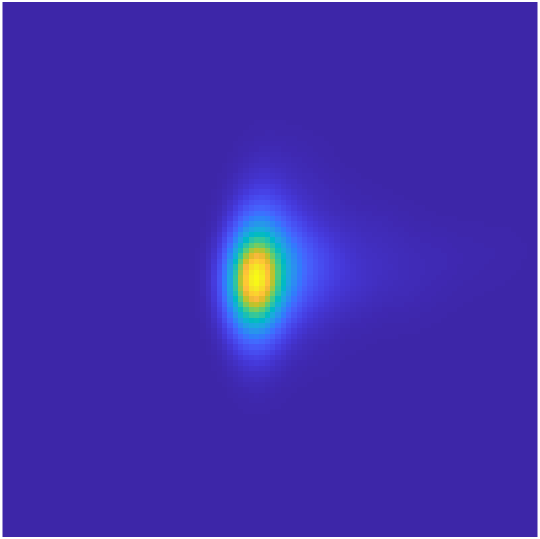}\\
\end{minipage}\hfill
\begin{minipage}{0.142\linewidth}
\includegraphics[width=1\linewidth]{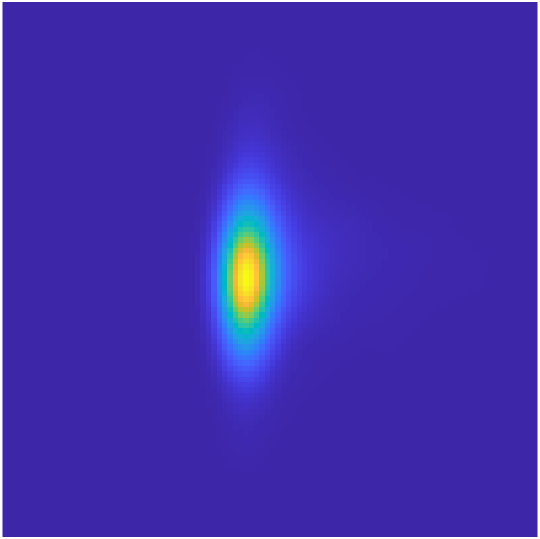}\\
\end{minipage}\hfill

\caption{Output of each intermediate RealNVP flow, projected onto dimensions 3 and 4, between a Gaussian and the Miniboone dataset. Top: trained without using the transport cost; bottom: with the cost.}
\label{RNVP_miniboone_34}
\end{figure}

\begin{figure}[H]
\centering
\begin{minipage}{0.142\linewidth}
\includegraphics[width=1\linewidth]{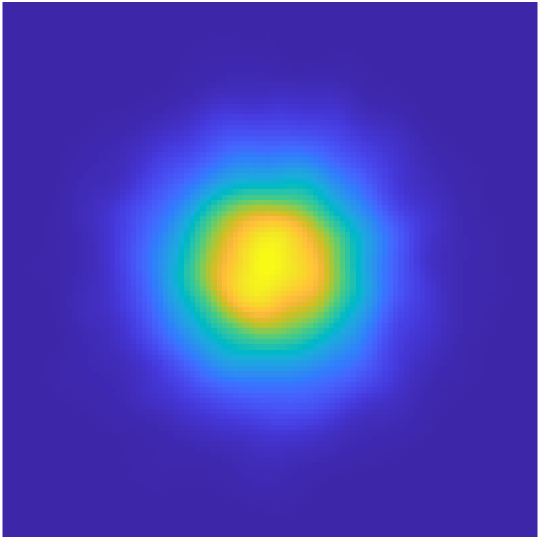}\\
\end{minipage}\hfill
\begin{minipage}{0.142\linewidth}
\includegraphics[width=1\linewidth]{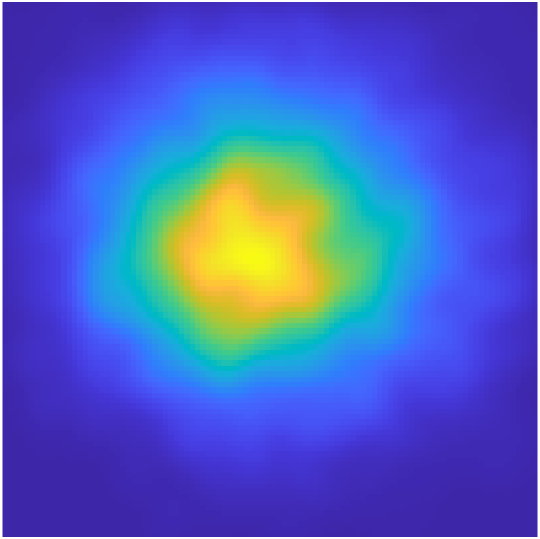}\\
\end{minipage}\hfill
\begin{minipage}{0.142\linewidth}
\includegraphics[width=1\linewidth]{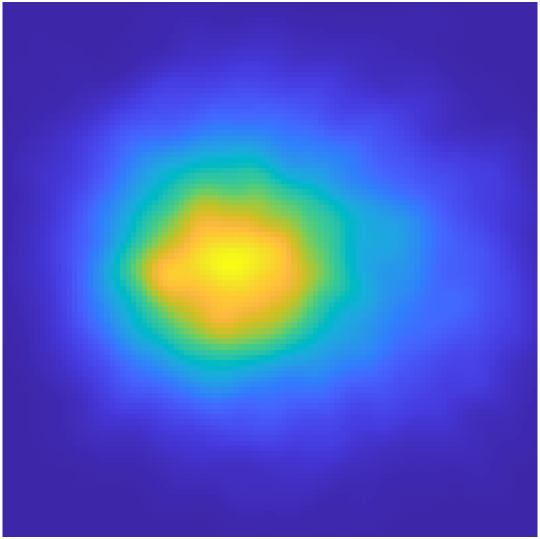}\\
\end{minipage}\hfill
\begin{minipage}{0.142\linewidth}
\includegraphics[width=1\linewidth]{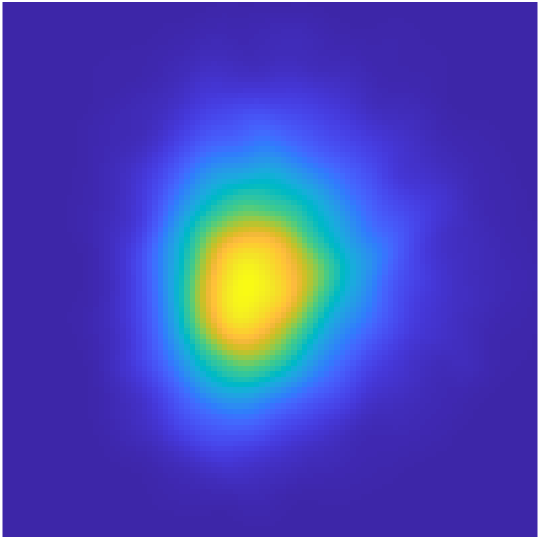}\\
\end{minipage}\hfill
\begin{minipage}{0.142\linewidth}
\includegraphics[width=1\linewidth]{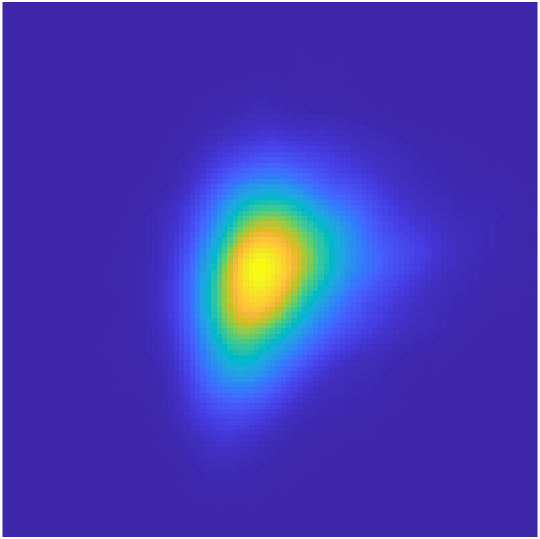}\\
\end{minipage}\hfill
\begin{minipage}{0.142\linewidth}
\includegraphics[width=1\linewidth]{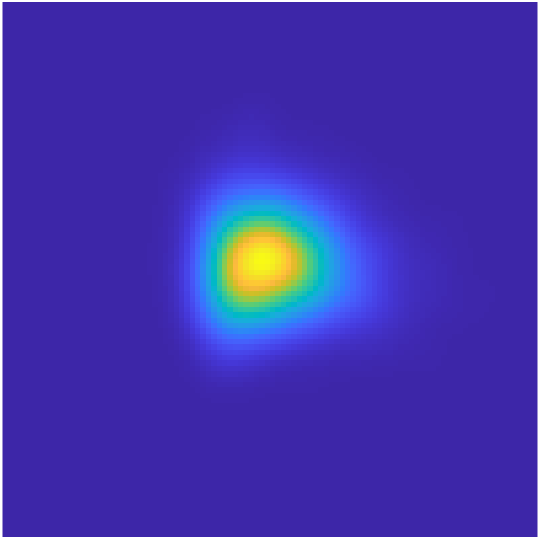}\\
\end{minipage}\hfill
\begin{minipage}{0.142\linewidth}
\includegraphics[width=1\linewidth]{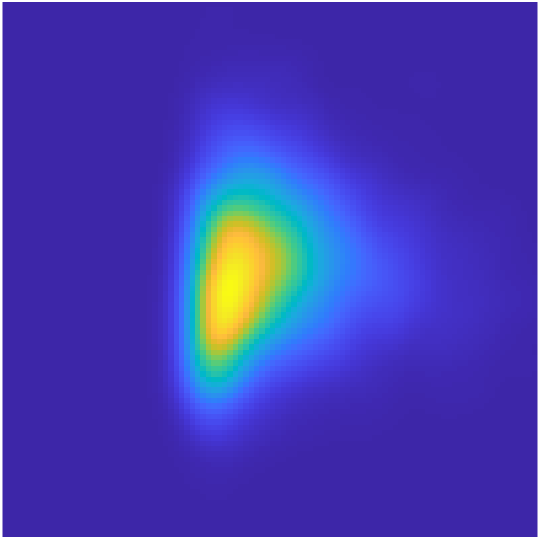}\\
\end{minipage}\hfill

\begin{minipage}{0.142\linewidth}
\includegraphics[width=1\linewidth]{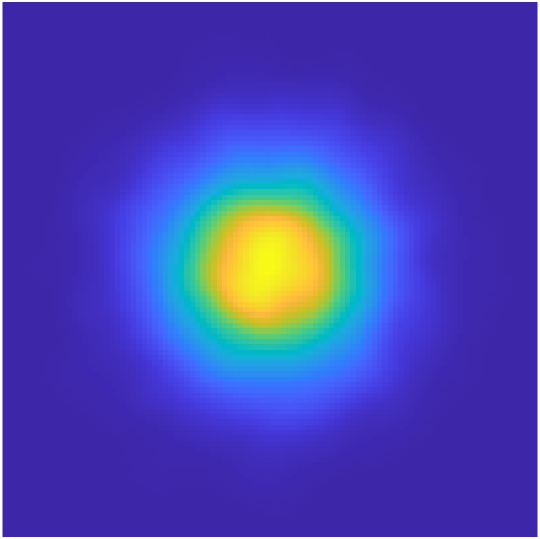}\\
\end{minipage}\hfill
\begin{minipage}{0.142\linewidth}
\includegraphics[width=1\linewidth]{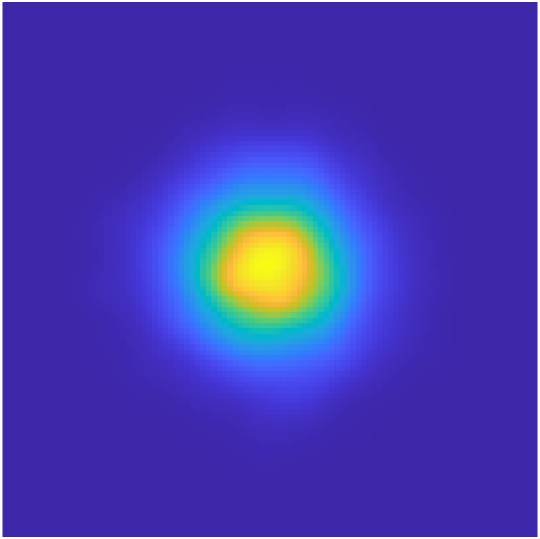}\\
\end{minipage}\hfill
\begin{minipage}{0.142\linewidth}
\includegraphics[width=1\linewidth]{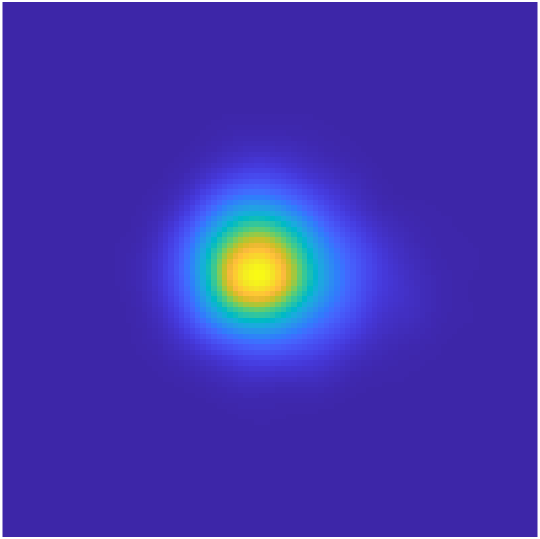}\\
\end{minipage}\hfill
\begin{minipage}{0.142\linewidth}
\includegraphics[width=1\linewidth]{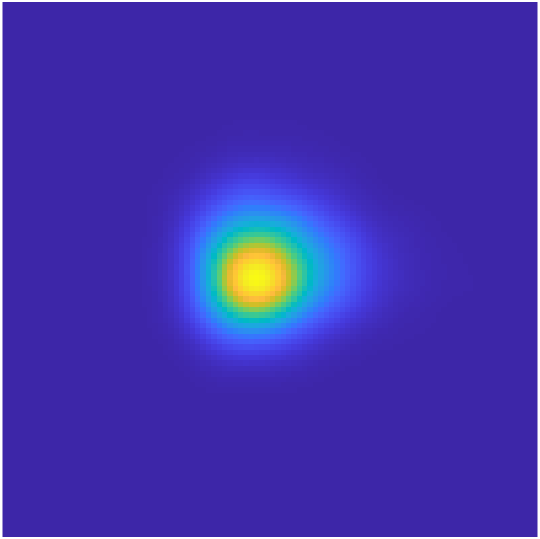}\\
\end{minipage}\hfill
\begin{minipage}{0.142\linewidth}
\includegraphics[width=1\linewidth]{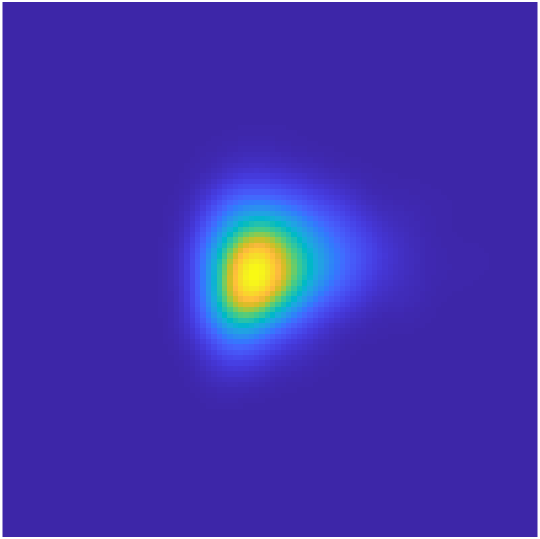}\\
\end{minipage}\hfill
\begin{minipage}{0.142\linewidth}
\includegraphics[width=1\linewidth]{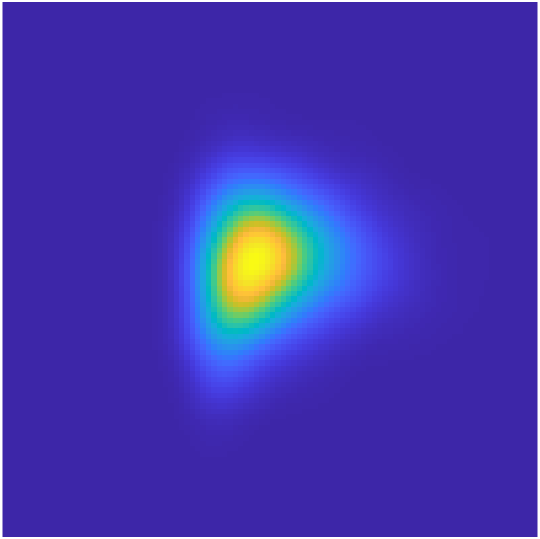}\\
\end{minipage}\hfill
\begin{minipage}{0.142\linewidth}
\includegraphics[width=1\linewidth]{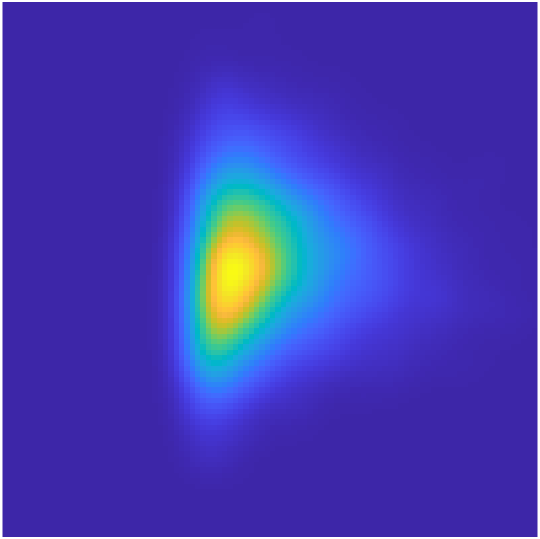}\\
\end{minipage}\hfill

\caption{Output of each intermediate RealNVP flow, projected onto dimensions 5 and 6, between a Gaussian and the Miniboone dataset. Top: trained without using the transport cost; bottom: with the cost.}
\label{RNVP_miniboone_56}
\end{figure}

\subsection{Image Datasets}

Here, we document the details of the image datasets as well as the Glow model. The SVHN dataset consists of color images of house numbers~\cite{SVHN}, and EuroSAT is composed of satellite images for various landscapes~\cite{EuroSAT}. In Table~\ref{Image_dataset_details}, we provide the meta-data for the image datasets used.

\begin{table}[H]
\caption{Meta-data for the image datasets used in Table~\ref{glow_table}.}
\label{Image_dataset_details}
\begin{center}
\begin{small}
\begin{sc}
\begin{tabular}{lccccc}
\toprule
 & MNIST & FMNIST & CIFAR-10 & SVHN & EUROSAT\\
\midrule
Number of Images & 60000 & 60000 & 50000 & 73257  & 27000 \\
Image Dimensions & $28 \times 28$ & $28 \times 28$ & $32 \times 32$ & $32 \times 32$ & $64 \times 64$\\
Number of Channels & 1 & 1 & 3 & 3 & 3\\
\bottomrule
\end{tabular}
\end{sc}
\end{small}
\end{center}
\vspace{-.2cm}
\end{table}

The Glow implementation we used is adapted from~\cite{glow_github}, which uses a simple CNN with 512 channels interlaced with actnorm as the conditioner in its affine coupling transforms. Each flow in Glow consists of an actnorm layer, an affine coupling, and a 1x1 convolution~\cite{glow}. The full model employs 32 flow transforms per level, and the number of levels depends on the image dataset specified in Table~\ref{glow_model_detail}. At the end of each level, the images are squeezed from the shape $c\times n \times n$ to $4c\times \frac{n}{2} \times \frac{n}{2}$, where $c, n$ are the number of channels and the side dimension, respectively. Afterwards, half of the channels go through an additional affine transformation and are not transformed further in the remaining flows. We optimize the models with Adam~\cite{ADAM} at the highest learning rate when training does not diverge. Further information about the Glow model is summarized in Table~\ref{glow_model_detail}.

\begin{table}[h]
\caption{Hyperparameters used to produce the results in Table~\ref{glow_table}.}
\label{glow_model_detail}
\begin{center}
\begin{small}
\begin{sc}
\begin{tabular}{lccccc}
\toprule
 & MNIST & FMNIST & CIFAR-10 & SVHN & EUROSAT\\
\midrule
Glow OT Weight & 1e-6 & 1e-6 & 5e-7 & 3e-7 & 3e-7\\
Number of levels & 2 & 2 & 4 & 4 & 4\\
\bottomrule
\end{tabular}
\end{sc}
\end{small}
\end{center}
\vspace{-.2cm}
\end{table}





\end{document}